\documentclass[11pt]{article}
\usepackage{mathtools,amsmath,amsthm,amssymb,cases}
\usepackage{epsfig}
\usepackage{leftidx}
\usepackage{graphicx}
\usepackage{amsfonts}
\usepackage{color,epstopdf}
\usepackage{cite}
\usepackage{graphicx,wrapfig,tabularx}
\usepackage{flafter}
\usepackage{fancyhdr}
\usepackage{stmaryrd}
\usepackage{graphicx,subfigure}
\usepackage{multicol,multirow,booktabs}
\usepackage{booktabs,threeparttable}
\usepackage[center]{caption2}
\usepackage{epstopdf}
\usepackage [latin1]{inputenc}
\usepackage{multirow}
\usepackage{enumerate}
\usepackage{enumitem}
\usepackage{algpseudocode,algorithm,algorithmicx}

\newtheorem{theorem}{Theorem}[section]
\newtheorem{lemma}{Lemma}[section]

\newtheorem{corollary}{Corollary}[section]
\newtheorem{example}{Example}

\newcommand{\zd}{\,\mathrm{d}}

\newcommand{\abs}[1]{\left|#1\right|}

\newcommand{\bra}[1]{\left(#1\right)}
\newcommand{\brab}[1]{\big(#1\big)}
\newcommand{\braB}[1]{\Big(#1\Big)}
\newcommand{\brat}[1]{(#1)}
\newcommand{\kbra}[1]{\left[#1\right]}
\newcommand{\kbrab}[1]{\big[#1\big]}

\newcommand{\kbrat}[1]{[#1]}
\newcommand{\myinner}[1]{\left\langle#1\right\rangle}
\newcommand{\myinnert}[1]{\langle#1\rangle}
\newcommand{\myinnerb}[1]{\big\langle#1\big\rangle}
\newcommand{\myinnerB}[1]{\Big\langle#1\Big\rangle}
\newcommand{\mynorm}[1]{\left\|#1\right\|}
\newcommand{\mynormb}[1]{\big\|#1\big\|}

\title{Average energy dissipation rates of additive implicit-explicit 
  Runge-Kutta methods for gradient flow problems}
\author{Hong-lin Liao
  \thanks{ORCID 0000-0003-0777-6832. School of Mathematics,
    Nanjing University of Aeronautics and Astronautics,
    Nanjing 211106, China; Key Laboratory of Mathematical Modeling
    and High Performance Computing of Air Vehicles (NUAA), MIIT, Nanjing 211106, China. Emails: liaohl@nuaa.edu.cn and liaohl@csrc.ac.cn.
    This author's work is supported by NSF of China under grant number 12071216.}
  \quad Xuping Wang\thanks{School of Mathematics, Nanjing University of Aeronautics and Astronautics, Nanjing 211106, China. Email: wangxp@nuaa.edu.cn.}
  \quad Cao Wen\thanks{School of Mathematics, Nanjing University of Aeronautics and Astronautics, Nanjing 211106, China. Email: 0131122255@mail.imu.edu.cn.}
}
\date{\today}

\begin{document}
  
\maketitle
  
\begin{abstract}
  A unified theoretical framework is suggested to examine the energy dissipation properties at all stages of additive implicit-explicit Runge-Kutta (IERK) methods up to fourth-order accuracy for gradient flow problems. We construct some parameterized IERK methods by applying the so-called first same as last method, that is, the diagonally implicit Runge-Kutta method with the explicit first stage and stiffly-accurate assumption for the linear stiff term, and applying the explicit Runge-Kutta method for the nonlinear term. The main part of the novel framework is to construct the differential forms and the associated differentiation  matrices of IERK methods by using the difference coefficients of method and the so-called discrete orthogonal convolution kernels. As the main result, we prove that an IERK method can preserve the original energy dissipation law unconditionally if the associated differentiation  matrix is positive semi-definite. The recent indicator, namely average energy dissipation rate, is also adopted for these multi-stage methods to evaluate the overall energy dissipation rate of an IERK method such that one can choose proper parameters in some parameterized IERK methods. It is found that the selection of method parameters in the IERK methods is at least as important as the selection of different IERK methods. Extensive numerical experiments are also included to support our theory.
  \\[1ex]   
  \\[1ex]
\end{abstract}

\section{Introduction}
\setcounter{equation}{0}
We propose a unified theoretical framework to examine the energy dissipation properties at all stages of additive implicit-explicit Runge-Kutta (in short, IERK) methods for solving the following semi-discrete semilinear parabolic problem, cf. \cite{FuYang:2022JCP,FuTangYang:2024,LiaoWang:2024arxiv,StuartHumphries:1998},
\begin{align}\label{problem: autonomous}
  u_h'(t)=M_h\kbra{L_hu_h(t)-g_h(u_h(t))},\quad u_h(t_0)=u_h^0,
\end{align}
where $L_h$ is a symmetric, positive definite matrix resulting from certain spatial discretization of stiff term, typically the Laplacian operator $-\Delta$ with periodic boundary conditions, $M_h$ is a symmetric, negative definite matrix (having the same size as $L_h$) resulting from certain spatial discretization of the mobility operator, and $g_h$ represents a nonlinear but non-stiff term. 
Without losing the generality,
the Fourier pseudo-spectral method
is assumed to approximate spatial operators and we define the discrete $L^2$ inner product $\myinner{u,v}:=v^Tu$ and the $L^2$ norm $\mynorm{v}:=\sqrt{\myinner{v,v}}$. Assume that there exists a non-negative Lyapunov function $G_h$ such that $g_h(v)=-\frac{\delta}{\delta v} G_h(v).$
Then the problem \eqref{problem: autonomous} can be formulated into  the following gradient flow system, see
\cite{AkrivisLiLi:2019,DuJuLiQiao:2021SIREV,GongZhaoWang:2020HEQRK,ShenXuYang:2018,
	ShinLeeLee:2017CMA,ShinLeeLee:2017},
\begin{align}\label{problem: gradient flows}
  \frac{\zd u_h}{\zd t}=M_h\frac{\delta E}{\delta u_h}\quad\text{with}\quad
  E[v_h]:=\frac1{2}\myinnert{v_h,L_hv_h}+\myinnert{G_h(v_h),1}.
\end{align}
The dynamics approaching the steady-state solution $u_h^*$, that is $L_hu_h^*=g_h(u_h^*)$, of this dissipative system \eqref{problem: autonomous} satisfies the following \textit{original} energy dissipation law
\begin{align}\label{problem: energy dissipation law}
  \frac{\zd E}{\zd t}=\myinnerB{\frac{\delta E}{\delta u_h},\frac{\zd u_h}{\zd t}}=
  \myinnerB{M_h^{-1}\frac{\zd u_h}{\zd t},\frac{\zd u_h}{\zd t}}\le0.
\end{align}

Our recent interests in the Runge-Kutta (RK) methods are firstly to develop some high-order starting procedures for the backward differentiation formulas \cite{LiLiao:2022,LiaoZhang:2021}, which are proven to possess certain discrete gradient structures, cf. \cite[Section 5.6]{StuartHumphries:1998}, and then preserve the energy dissipation law of phase field models with the original energy added with a small positive term in their original forms. Thus, it would be expected that the RK methods also have some discrete gradient structures and preserve certain discrete energy dissipation laws for gradient flow problems.  
Combining the scalar auxiliary variable (SAV) method \cite{ShenXuYang:2018} and the Gaussian-type RK methods, arbitrarily high-order energy-stable schemes were constructed \cite{AkrivisLiLi:2019,GongZhaoWang:2020HEQRK}. However, these SAV schemes only satisfy a modified energy law, which may not necessarily ensure the energy stability in the original form \cite{JiangZhangZhao:2022,FuTangYang:2024}. Recently, Fu et al. \cite{FuYang:2022JCP,FuTangYang:2024} proved that some exponential time difference RK methods maintain the decaying of original energy. For some parameterized explicit exponential RK methods in \cite{HochbruckOstermann:2010ActaNu},  the discrete version of \eqref{problem: energy dissipation law} was established in \cite{LiaoWang:2024arxiv} for certain range of the method parameters. Since the functions of matrix exponential are always involved, these explicit RK methods would be also computationally intensive. Actually, the efficient algorithm to accurately compute the matrix exponential is still limited \cite{MolerVanLoan:2003,FasiGaudreaultLundSchweitzer:2024}. 

In this study, we will focus on whether and to what extent the IERK methods preserve
the original energy dissipation law \eqref{problem: energy dissipation law}.
Due to the relative ease of implementation, the diagonally implicit Runge-Kutta (DIRK) methods are possibly the most widely used implicit Runge-Kutta method in practical applications, see the recent review \cite{KennedyCarpenter:2016} and abundant references therein, especially involving stiff differential equations. The DIRK method has greatly reduced the computational complexity of the fully implicit Runge-Kutta method, but it still requires an iterative method to solve the nonlinear equations at each stage. In order to improve the computational efficiency, the IERK methods have attracted much attention and are widely used, see \cite{AscherRuuthSpiteri:1997,BoscarinoRusso:2009,BoscarinoPareschiRusso:2017,CardoneJackiewiczSanduZhang:2014MMA,CooperSayfy:1983,DimarcoPareschi:2013,IzzoJackiewicz:2017,GuiWangChen:2023,Jin:1995,KennedyCarpenter:2003,KanevskyCarpenterGottliebHesthaven:2007,LiuZou:2006,ShinLeeLee:2017,ShinLeeLee:2017CMA}.
Cooper and Sayfy \cite{CooperSayfy:1983} presented some IERK methods up to fourth-order accuracy with the implicit part is A-stable.
Kennedy and Carpenter \cite{KennedyCarpenter:2003} constructed high-order IERK methods from third- to fifth-order to simulate convection-diffusion-reaction equations. The widespread ARS-type IERK methods by Ascher, Ruuth and Spiteri  were developed in \cite{AscherRuuthSpiteri:1997} to solve the convection-diffusion problems. These ARS-type methods had better stability regions than implicit-explicit multi-step schemes over a wide parameter range.  Some ARS-type IERK methods from second- to fourth-order were also proposed in \cite{LiuZou:2006} for nonlinear differential equations with constraints, such as the Navier-Stokes equation. Cardone et.al. \cite{CardoneJackiewiczSanduZhang:2014MMA} proposed a class of ARS-type IERK methods up to fourth-order based on extrapolation of the stage values at the current step by stage values in the previous step. Izzo and Jackiewicz \cite{IzzoJackiewicz:2017} constructed some parameterized IERK methods up to fourth-order by choosing the method parameters to maximize the regions of absolute stability for the explicit part in the method with the assumption of A-stable implicit part.

We will consider another class of parameterized IERK methods and choose the method parameters to maintain the original energy dissipation law \eqref{problem: energy dissipation law} as much as possible.  
Let $u_h^k$ be the numerical approximation of $u_h(t_k)$ at the mesh point $t_k$ for $0\le k\le N$. For a $s$-stage Runge-Kutta method, let $U^{n,i}$ be the approximation of $u_h(t_{n-1}+c_{i}\tau)$ at the abscissas $c_1:=0$, $c_i>0$ for $2\le i\le s-1$, and $c_{s}:=1$.   
To integrate the semilinear parabolic problem \eqref{problem: autonomous} from the time $t_{n-1}$ ($n\ge1$) to the next grid point $t_n=t_{n-1}+\tau$ (typically, $\tau$ also represents a variable-step size), 
one can construct the following $s$-stage IERK methods with explicit first stage, see \cite{KennedyCarpenter:2003,LiuZou:2006,Zhong:1995}
\begin{subequations}\label{Scheme: general IERK}  
\begin{align}   
  &~U^{n,i}:=U^{n,1}+\tau\sum_{j=1}^{i}a_{i,j}M_hL_hU^{n,j}
    -\tau\sum_{j=1}^{i-1}\hat{a}_{i,j}M_hg_h(U^{n,j}),
    \quad\text{$2\le i\le s-1$,}\\
  &~U^{n,s}:=U^{n,1}+\tau\sum_{j=1}^{s}b_{j}M_hL_hU^{n,j}-\tau\sum_{j=1}^{s-1}\hat{b}_{j}M_hg_h(U^{n,j}),
\end{align}
\end{subequations}
where $U^{n,1}:=u_h^{n-1}$ and $u_h^{n}:=U^{n,s}$.
The implicit part of IERK methods \eqref{Scheme: general IERK} approximating the stiff linear term $L_h$ can be represented by the abscissa vector $\mathbf{c}$, the coefficient matrix $A$ and the vector of weights $\mathbf{b}$, 
\begin{equation*}
  \begin{array}{c|c}
    \mathbf{c} & A \\
    \hline\\[-8pt]  & \mathbf{b}^T
  \end{array}
  = \begin{array}{c|ccccc}
    c_{1}=0 & 0 &  &  &  &   \\
    c_{2} & a_{21} & a_{22} &  &  &   \\
    c_{3} & a_{31} & a_{32} & a_{33} &  &   \\
    \vdots & \vdots & \vdots & \ddots & \ddots &  \\[2pt]
    c_{s}=1 & a_{s,1} & a_{s,2} &  \cdots  & a_{s,s-1}   & a_{s,s} \\[2pt]
    \hline  & a_{s,1} & a_{s,2} &  \cdots  & a_{s,s-1}   & a_{s,s}
  \end{array},\end{equation*}
where the abscissa $c_i=\sum_{j=1}^{i}a_{i,j}$ for $1\le i\le s$, and we use the so-called first same as last methods \cite{OlssonSoderlind:2000,KennedyCarpenter:2016}, that is, the DIRK methods with the explicit first stage and the stiffly-accurate assumption: $b_{j}=a_{s,j}$ for $1\leq j\leq s$ such that 
$$c_s=\sum_{j=1}^{s}a_{s,j}=\sum_{j=1}^{s}b_{j}=1.$$
That is, we consider Lobatto-type DIRK methods, with the number of implicit stage $s_{\mathrm{I}}:=s-1,$ which differs from the widespread Radau-type DIRK methods with $a_{11}\neq0$ or $a_{i1}=0$ for $2\le i\le s$.

The explicit part of IERK methods \eqref{Scheme: general IERK}, namely explicit Runge-Kutta methods, for the nonlinear term $g_h$ can be represented by the abscissa vector $\hat{\mathbf{c}}$, the strictly lower triangular coefficient matrix $\widehat{A}$ and the vector of weights $\hat{\mathbf{b}}$, 
\begin{equation*}
\begin{array}{c|c}
  \hat{\mathbf{c}} & \widehat{A} \\
\hline\\[-8pt]  & \hat{\mathbf{b}}^T
\end{array}
=\begin{array}{c|ccccc}
    \hat{c}_{1}=0 & 0 &  &  &  &   \\
    \hat{c}_{2} & \hat{a}_{21} & 0 &  &  &   \\
    \hat{c}_{3} & \hat{a}_{31} & \hat{a}_{32} & 0 &  &   \\
    \vdots & \vdots & \vdots & \ddots & \ddots &  \\[2pt]
    \hat{c}_{s}=1 & \hat{a}_{s,1} & \hat{a}_{s,2} &  \cdots  & \hat{a}_{s,s-1}   & 0 \\[1pt]
    \hline\\[-8pt]   & \hat{a}_{s,1} & \hat{a}_{s,2} &  \cdots  & \hat{a}_{s,s-1}   & 0
\end{array},\end{equation*}
where the abscissa $\hat{c}_i=\sum_{j=1}^{i-1}\hat{a}_{i,j}$ for $1\le i\le s$, and introduce, to simplify our notations, 
\begin{align}\label{cond: Matrix setting}
  \hat{b}_j:=\hat{a}_{s,j},\quad 1\leq j\leq s-1\quad\text{such that}\quad \hat{c}_{s} = \sum_{j=1}^{s-1}\hat{a}_{s,j}=1.
\end{align} 
Always, we assume that $\hat{a}_{k+1,k}\neq0$ for any $1\le k\le s-1$. The IERK methods \eqref{Scheme: general IERK} become
\begin{align}\label{Scheme: general IERK2}        &~U^{n,i}=U^{n,1}+\tau\sum_{j=1}^{i}a_{i,j}M_hL_hU^{n,j}-\tau\sum_{j=1}^{i-1}\hat{a}_{i,j}M_hg_h(U^{n,j})\quad\text{for $2\le i\le s$.}
\end{align}
Since the  Lobatto-type DIRK methods are used in the implicit part, we call \eqref{Scheme: general IERK2} as Lobatto-type IERK methods. They are quite different from the Radau-type IERK method containing the Radau-type DIRK methods in the implicit part, which are also known as ARS-type IERK methods \cite{AscherRuuthSpiteri:1997,FuTangYang:2024,IzzoJackiewicz:2017,ShinLeeLee:2017,GuiWangChen:2023,LiuZou:2006}.

In general, the consistency of IERK methods  \eqref{Scheme: general IERK2} requires the canopy node condition 
\begin{align}\label{cond: canopy node setting}
  c_{i}=\hat{c}_{i}\quad\text{or}\quad \sum_{j=1}^{i}a_{i,j}=\sum_{j=1}^{i-1}\hat{a}_{i,j}\quad
  \text{for $1\le i\le s$},
\end{align} 
which makes the form \eqref{Scheme: general IERK2} invariant under the transformation of IERK method to the non-autonomous system. 
Under the canopy node condition $\hat{\mathbf{c}}={\mathbf{c}}$, Table \ref{table: order condition} 
lists the order conditions for the coefficient matrices and the weights vector to make the IERK method \eqref{Scheme: general IERK2}  accurate up to fourth-order. A detailed description of these order conditions can also be found in
\cite{AscherRuuthSpiteri:1997,CooperSayfy:1983,IzzoJackiewicz:2017,ShinLeeLee:2017}.

\begin{table}[htb!]\centering
  \begin{threeparttable}
    \centering 
    \renewcommand\arraystretch{1.3}
    \belowrulesep=0pt\aboverulesep=0pt
    \caption{Order conditions for IERK methods up to fourth-order accuracy.}
    \label{table: order condition}
    \vspace{2mm}
    \begin{tabular}{c|cc|c}
      \toprule 
      \multirow{2}*{Order} & \multicolumn{2}{c|}{Stand-alone conditions} & Coupling condition \\
      \cmidrule{2-4}
      & Implicit part & Explicit part & \\
      \midrule 
      1 & $\mathbf{b}^{T} \mathbf{1}=1$ & $\hat{\mathbf{b}}^{T} \mathbf{1} = 1$ & - \\[2pt]
      \hline 2 & $\mathbf{b}^{T} \mathbf{c} = \tfrac1{2}$ & $\hat{\mathbf{b}}^{T} \mathbf{c} = \tfrac1{2}$ & - \\[2pt]
      \hline 3 & $\mathbf{b}^{T} \mathbf{c}^{.2} = \tfrac1{3}$ & $\hat{\mathbf{b}}^{T} \mathbf{c}^{.2} = \tfrac1{3}$ & \\
      & $\mathbf{b}^{T} A  \mathbf{c} = \tfrac1{6}$ & $\hat{\mathbf{b}}^{T} \widehat{A}  \mathbf{c} =\tfrac1{6}$ & {$\mathbf{b}^{T} \widehat{A}  \mathbf{c} = \tfrac1{6}, \; \hat{\mathbf{b}}^{T} A  \mathbf{c} = \tfrac1{6}$} \\[2pt]
      \hline 
      & $\mathbf{b}^{T} \mathbf{c}^{.3} = \frac{1}{4}$ & $ \hat{\mathbf{b}}^{T} \mathbf{c}^{.3} = \frac{1}{4}$  &  \\[2pt]
       & { $\mathbf{b}^{T} [\mathbf{c} \odot (A  \mathbf{c})] = \frac{1}{8} $} &  {$ \hat{\mathbf{b}}^{T} [\mathbf{c} \odot (\widehat{A}  \mathbf{c})] = \frac{1}{8}$}  & {$\mathbf{b}^{T} [\mathbf{c} \odot (\widehat{A}  \mathbf{c})] = \frac{1}{8}, \;\hat{\mathbf{b}}^{T} [\mathbf{c} \odot (A  \mathbf{c})] = \frac{1}{8} $} \\[2pt]
      4 & {$\mathbf{b}^{T} A \mathbf{c}^{.2} = \frac{1}{12}$} & {$\hat{\mathbf{b}}^{T} \widehat{A} \mathbf{c}^{.2} = \frac{1}{12}$} & $\mathbf{b}^{T} \widehat{A} \mathbf{c}^{.2} = \frac{1}{12}, \; \hat{\mathbf{b}}^{T} A \mathbf{c}^{.2} = \frac{1}{12}$  \\[2pt]
      & $\mathbf{b}^{T} A^{2}  \mathbf{c} =\frac{1}{24}$ &  $\hat{\mathbf{b}}^{T} \widehat{A}^{2}  \mathbf{c} = \frac{1}{24}$  & { $ \mathbf{b}^{T} \widehat{A}^{2}  \mathbf{c} = \frac{1}{24}, \; \hat{\mathbf{b}}^{T} A \widehat{A}  \mathbf{c} = \frac{1}{24}, \; \hat{\mathbf{b}}^{T} \widehat{A} A  \mathbf{c} = \frac{1}{24}$} \\[2pt]
      & & & {$\mathbf{b}^{T} A \widehat{A}  \mathbf{c} = \frac{1}{24}, \; \mathbf{b}^{T} \widehat{A} A  \mathbf{c} = \frac{1}{24}, \; \hat{\mathbf{b}}^{T} A^{2} \mathbf{c} = \frac{1}{24}$} \\
      \bottomrule
    \end{tabular}
    \footnotesize
    \tnote{* For the vectors $\mathbf{x}=(x_1,x_2,\cdots,x_s)^T$ and $\mathbf{y}=(y_1,y_2,\cdots,y_s)^T$, $\mathbf{x}\odot\mathbf{y}:=(x_1y_1,x_2y_2,\cdots,x_sy_s)^T$ and $\mathbf{x}^{.m}:=\mathbf{x}\odot\mathbf{x}^{.(m-1)}$ for $m>1$.}
  \end{threeparttable}
\end{table}

In simulating the semilinear parabolic problems \eqref{problem: autonomous} and related gradient flow problems \eqref{problem: gradient flows}, IERK methods turned out to be very competitive. Shin et al. \cite{ShinLeeLee:2017CMA} observed that the Radau-type IERK methods combined with convex splitting technique (CSRK) exhibit the original energy stability in numerical experiments. 
Subsequently, they proved \cite{ShinLeeLee:2017} that the proposed CSRK schemes can unconditionally maintain the original energy stability by assuming $\hat{a}_{ij} = a_{i,j+1}$ in \eqref{Scheme: general IERK}. The proposed second-order CSRK scheme requires at least four stages, while the third-order CSRK scheme requires at least seven stages.
The proposed CSRK schemes in \cite{ShinLeeLee:2017} are nonlinear and would be computationally expensive due to the required inner iterations at each time stage. 
Very recently, Fu et al. \cite{FuTangYang:2024} derived some sufficient conditions of Radau-type IERK methods to maintain the decay of original energy for the gradient flow problems and presented some concrete schemes up to third-order accuracy in time. It is noted that, the Radau-type IERK methods are only focused in the mentioned works, while we will construct some Lobatto-type IERK methods and present a unified theoretical framework to examine the energy dissipation properties of IERK methods up to fourth-order accurate for gradient flow problems \eqref{problem: gradient flows}.

The unified theoretical framework, presented in the next section, is inspired by two aspects: one is the idea by treating the
stiffly-accurate (but not necessarily algebraically stable) DIRK method as the composite linear
multistep method; the other is the recent discrete energy technique for the stability and convergence of multistep backward differentiation formulas \cite{LiLiao:2022,LiaoZhang:2021}.
We will construct the differential forms and the associated differentiation matrices of IERK methods by using the difference coefficients of method and the so-called discrete orthogonal convolution kernels. It is proven that an IERK method can preserve the original energy dissipation law unconditionally if the associated differentiation  matrix is positive semi-definite, see Theorem \ref{thm: energy stability}. The recent indicator in \cite{LiaoWang:2024arxiv}, namely average energy dissipation rate, is also defined for the IERK methods to evaluate the overall energy dissipation rate when applied to the gradient flow problems \eqref{problem: gradient flows} such that one can choose proper parameters in some parameterized IERK methods, see Lemma \ref{lemma: average dissipation rate}. 
As described later, this framework will be fit for both Lobatto-type and Radau-type IERK methods, but is quite different from those in previous studies \cite{ShinLeeLee:2017,ShinLeeLee:2017CMA,FuTangYang:2024}.  


Section \ref{sec: IERK2} addresses the details in constructing three parameterized second-order IERK (IERK2) methods and establishes the original energy dissipation laws for the resulting IERK2 methods. Also, we choose the optimal parameter in these IERK2 methods by using the concept of average dissipation rates,  and present extensive tests to support our theoretical predictions, see Table \ref{table: comparison of 2nd-order methods}, in which the parameter choices for the energy stability
of six IERK2 methods are summarized. Third-order parameterized IERK (IERK3) methods are discussed and tested in Section \ref{section: third-order methods} and, at the end of section, Table \ref{table: comparison of 3rd-order methods} collects some parameter choices for the energy stability of
eight IERK3 methods. To show the existence of energy
stable IERK methods with fourth-order time accuracy, two approximately fourth-order IERK (IERK4-A) methods are presented and tested in Section \ref{sec: IERK4-A}. In the last section, we summarize this article and present some further issues to be studied.

\section{Stage energy laws of IERK methods}
\setcounter{equation}{0}


Requiring $U^{n,i}=u_h^n=u_h^*$ for all
$i$ and $n \geq 0$ immediately shows that the canopy node condition \eqref{cond: canopy node setting} makes the IERK method \eqref{Scheme: general IERK2} preserve naturally the equilibria $u_h^*$ of the gradient flow problem \eqref{problem: gradient flows}, that is, $L_hu_h^*=g_h(u_h^*)$. So one can reformulate the form \eqref{Scheme: general IERK2} into the following steady-state preserving form
\begin{align}\label{Scheme: IERK steady-state preserving}       U^{n,i}=&\,U^{n,1}+\tau\sum_{j=1}^{i}a_{i,j}M_h\brab{L_hU^{n,j}-L_hU^{n,1}}-\tau\sum_{j=1}^{i-1}\hat{a}_{i,j}M_h\kbrab{g_h(U^{n,j})-L_hU^{n,1}}\nonumber\\
=&\,U^{n,1}+\tau\sum_{j=2}^{i}a_{i,j}M_h\brab{L_hU^{n,j}-L_hU^{n,1}}-\tau\sum_{j=1}^{i-1}\hat{a}_{i,j}M_h\kbrab{g_h(U^{n,j})-L_hU^{n,1}}
\end{align}
for $2\le i\le s$, in which we drop the terms with the coefficients $a_{i1}$ for $2\le i\le s$.
In this sense, we define the lower triangular coefficient matrices for the implicit and explicit parts, respectively,
$$A_{\mathrm{I}}:=\brab{a_{i+1,j+1}}_{i,j=1}^{s_{\mathrm{I}}}\quad\text{and}\quad A_{\mathrm{E}}:=\brab{\hat{a}_{i+1,j}}_{i,j=1}^{s_{\mathrm{I}}}\,.$$
The two matrices are always required in our theory on the energy dissipation property, while the coefficient vector $\mathbf{a}_{1}:=(a_{21},a_{31},\cdots,a_{s1})^T$ would be not involved directly. Due to the addition of $s_{\mathrm{I}}$ coefficients (degree of freedom) in the vector $\mathbf{a}_{1}$, as shown later, it would be useful in designing some computationally effective Lobatto-type IERK methods.

\subsection{Our theoretical framework}

Motivated by the stabilization idea from Du et al. \cite{DuJuLiQiao:2021SIREV} and Fu  et al. \cite{FuYang:2022JCP,FuTangYang:2024}, 
we introduce the following stabilized operators with a parameter $\kappa\ge0$,
\begin{align}\label{def: stabilized parameter}
  L_{\kappa}:=L_h+\kappa I\quad\text{and}\quad g_{\kappa}(u):=g_h(u)+\kappa u,
\end{align} 
such that the problem \eqref{problem: autonomous} becomes the stabilized version
\begin{align}\label{problem: stabilized version}
  u_h'(t)=M_h\kbra{L_{\kappa}u_h(t)-g_{\kappa}(u_h)},\quad u_h(t_0)=u_h^0.
\end{align}
Thus, applying \eqref{Scheme: IERK steady-state preserving} to \eqref{problem: stabilized version}, 
we have the following  IERK method
\begin{align}\label{Scheme: general IERK stabilized2}     
  U^{n,i+1}=U^{n,1}+\tau\sum_{j=1}^{i}a_{i+1,j+1}\sum_{\ell=1}^{j}M_hL_{\kappa}\delta_{\tau}U^{n,\ell+1}
  -\tau\sum_{j=1}^{i}\hat{a}_{i+1,j}M_h\kbrab{g_{\kappa}(U^{n,j})-L_{\kappa}U^{n,1}}
\end{align}
for $1\le i\le s_{\mathrm{I}}$, where the (stage) time difference $\delta_{\tau}U^{n,\ell+1}:=U^{n,\ell+1}-U^{n,\ell}$ for $\ell\ge 1$.

To make our idea more concise, we assume further that the nonlinear function $g_h$ is Lipschitz continuous with a constant $\ell_{g}>0$, cf. \cite{StuartHumphries:1998} or the recent discussions in \cite{FuTangYang:2024}. 
In theoretical manner, the stabilization parameter $\kappa$ in \eqref{def: stabilized parameter} is chosen properly large (determining the minimum stabilized parameter is out of our current scope although it is also practically useful). In this sense,  if an IERK method is proven to maintain the original energy dissipation law \eqref{problem: energy dissipation law} unconditionally, we mean that this IERK method can be stabilized by setting a properly large parameter $\kappa$. 
To derive the energy dissipation law of the general IERK method \eqref{Scheme: general IERK stabilized2}, we need the following result. 

\begin{lemma}\label{lemma: origional energy derivation}\cite{LiaoWang:2024arxiv}
  If $g_h$ is Lipschitz-continuous with a constant $\ell_{g}>0$ and 
  $\kappa\ge2\ell_g$, then
  \begin{align*}
    \myinnerb{u-v,\tfrac12L_{\kappa} (u+v)-g_{\kappa}(v)} \geq E[u]-E[v]+\tfrac{1}2\brat{\kappa-2\ell_g}\mynormb{u-v}^2,
  \end{align*}
  where the energy $E$ is defined in \eqref{problem: gradient flows}.
\end{lemma}

    Our theoretical framework contains three main steps: 
    \begin{description}
      \item[(1)]  \underline{Compute difference coefficients}: we introduce a class of difference coefficients, for $i=2,\cdots,s$,
      \begin{align*}
        &\underline{a}_{i,i}:=a_{i,i}\quad\text{and}\quad
        \underline{a}_{i,j}:=a_{i,j}-a_{i-1,j}\quad\text{for $1\le j\le i-1$;}\\
        & \underline{\hat{a}}_{i,i-1}:=\hat{a}_{i,i-1}\quad\text{and}\quad
        \underline{\hat{a}}_{i,j-1}:=\hat{a}_{i,j-1}-\hat{a}_{i-1,j-1}\quad\text{for $2\le j\le i-1$.}
      \end{align*}
      It is not difficult to derive from \eqref{Scheme: general IERK stabilized2} that
      \begin{align}\label{Scheme: general IERK stabilized2II}     
        M_h^{-1}\delta_{\tau}U^{n,i+1}=&\,
        \tau\sum_{j=1}^{i}\underline{a}_{i+1,j+1}\sum_{\ell=1}^jL_{\kappa}\delta_{\tau}U^{n,\ell+1}
        +{\tau}\sum_{j=1}^{i}\underline{\hat{a}}_{i+1,j}\kbra{L_{\kappa}U^{n,1}-g_{\kappa}(U^{n,j})}
        \end{align}
       for $1\le i\le s_{\mathrm{I}}$. The associated Butcher difference tableaux read
      \begin{equation*}
        \begin{array}{c|ccccc}
          c_{1}=0 & 0 &  &  &  &   \\
          c_{2} & \underline{a}_{21} & \underline{a}_{22} &  &  &   \\
          c_{3} & \underline{a}_{31} & \underline{a}_{32} & \underline{a}_{33} &  &   \\
          \vdots & \vdots & \vdots & \ddots & \ddots &  \\[2pt]
          c_{s}=1 & \underline{a}_{s,1} & \underline{a}_{s,2} &  \cdots  & \underline{a}_{s,s-1}   & \underline{a}_{s,s} \\[2pt]
          \hline  & 0 & 0 & \cdots 
          & 0  & 0
          \end{array},\quad
        \begin{array}{c|ccccc}
          \hat{c}_{1}=0 & 0 &  &  &  &   \\
          \hat{c}_{2} & \underline{\hat{a}}_{21} & 0 &  &  &   \\
          \hat{c}_{3} & \underline{\hat{a}}_{31} & \underline{\hat{a}}_{32} & 0 &  &   \\
          \vdots & \vdots & \vdots & \ddots & \ddots &  \\[2pt]
          \hat{c}_{s}=1 & \underline{\hat{a}}_{s,1} & \underline{\hat{a}}_{s,2} &  \cdots  & \underline{\hat{a}}_{s,s-1}   & 0 \\[1pt]
          \hline\\[-8pt]   & 0 & 0 & \cdots 
          & 0  & 0
        \end{array}.
        \end{equation*}
    \item[(2)] \underline{Determine DOC kernels and differential form}: we introduce the so-called discrete orthogonal convolution (DOC) kernels 
    $\underline{\theta}_{k,j}(z)$ with respect to the explicit coefficient $\underline{\hat{a}}_{i+1,j}$, cf. \cite{LiLiao:2022,LiaoTangZhou:2024,LiaoWang:2024arxiv,LiaoZhang:2021}, 
    \begin{align}\label{eq: orthogonal procedureII}
      \underline{\theta}_{k,k}:=\frac1{\underline{\hat{a}}_{k+1,k}}\quad\text{and}\quad
      \underline{\theta}_{k,j}:=-\sum_{\ell=j+1}^{k}\underline{\theta}_{k,\ell}
      \frac{\underline{\hat{a}}_{\ell+1,j}}{\underline{\hat{a}}_{j+1,j}},
      \quad\text{$1\leq j\le k-1$}
    \end{align}
    for $k=1,2,\cdots,s_{\mathrm{I}}$.
    It is easy to check the following discrete orthogonal identity,
    \begin{align}\label{eq: orthogonal identityII}
      \sum_{\ell=j}^{m}\underline{\theta}_{m,\ell}\underline{\hat{a}}_{\ell+1,j}\equiv\delta_{m,j}\quad \text{for $1\leq j\leq m\leq s_{\mathrm{I}}$},
    \end{align}
    where $\delta_{m,j}$ is the Kronecker delta symbol with $\delta_{m,j}=0$ if $j\neq m$.
    Multiplying the third term of \eqref{Scheme: general IERK stabilized2II} by the DOC kernels  $\underline{\theta}_{k,i}$,  and summing $i$ from 1 to $k$, one can exchange the summation order and apply the discrete orthogonal identity \eqref{eq: orthogonal identityII} to find that
    \begin{align*}  
      \tau\sum_{i=1}^{k}\underline{\theta}_{k,i}
      \sum_{j=1}^{i}\underline{\hat{a}}_{i+1,j}&\,\kbra{L_{\kappa}U^{n,1}-g_{\kappa}(U^{n,j})}
      =\tau\sum\limits_{j=1}^{k}\kbra{L_{\kappa}U^{n,1}-g_{\kappa}(U^{n,j})}
      \sum_{i=j}^{k}\underline{\theta}_{k,i}\underline{\hat{a}}_{i+1,j}
      \\
      =&\,\tau L_{\kappa}U^{n,1}-\tau g_{\kappa}(U^{n,k})
      =\tau L_{\kappa}U^{n,k+1}-\tau g_{\kappa}(U^{n,k})
      -\tau\sum_{\ell=1}^{k}L_{\kappa}\delta_{\tau}U^{n,\ell+1}\\
      =&\,\tfrac{\tau}2L_{\kappa}(U^{n,k+1}+U^{n,k})-\tau g_{\kappa}(U^{n,k})
      -\tau\sum_{\ell=1}^{k}\brab{1-\tfrac{1}2\delta_{k,\ell}}L_{\kappa}\delta_{\tau}U^{n,\ell+1}
    \end{align*}
    for $1\le k\le s_{\mathrm{I}}$. Similarly, by 
    multiplying the second term of \eqref{Scheme: general IERK stabilized2II} by the DOC kernels  $\underline{\theta}_{k,i}$,  and summing $i$ from 1 to $k$, one has
    \begin{align*}  
      \tau\sum_{i=1}^{k}\underline{\theta}_{k,i}
      \sum_{j=1}^{i}\underline{a}_{i+1,j+1}\sum_{\ell=1}^jL_{\kappa}\delta_{\tau}U^{n,\ell+1}
      =\tau\sum\limits_{j=1}^{k}\sum_{\ell=1}^jL_{\kappa}\delta_{\tau}U^{n,\ell+1}  \sum_{i=j}^{k}\underline{\theta}_{k,i}\underline{a}_{i+1,j+1}        
    \end{align*}
    for $1\le k\le s_{\mathrm{I}}$. 
    With the help of the above two equalities, we have an equivalent form of 
    the IERK method \eqref{Scheme: general IERK stabilized2},
    \begin{align*}  
      \sum\limits_{\ell=1}^{k}\braB{\underline{\theta}_{k,\ell}-{\tau}M_hL_{\kappa}\sum_{j=\ell}^k
        \sum_{i=j}^{k}\underline{\theta}_{k,i}\underline{a}_{i+1,j+1}
        +{\tau}M_hL_{\kappa}-\tfrac{1}2{\tau}M_hL_{\kappa}\delta_{k,\ell}}\delta_{\tau}U^{n,\ell+1}\\
      =\tau M_h\kbra{\tfrac{1}2L_{\kappa}(U^{n,k+1}+U^{n,k})-g_{\kappa}(U^{n,k})}  \quad\text{for $1\le k\le s_{\mathrm{I}}$.}      
    \end{align*}
    Thus we obtain the following differential form of 
    the IERK method \eqref{Scheme: general IERK stabilized2}
    \begin{align}\label{Scheme: DOC IERK stabilized2II}     
      \sum_{\ell=1}^{k}d_{k,\ell}(\tau M_h L_{\kappa})\delta_{\tau}U^{n,\ell+1}    
      =\tau M_h\kbra{\tfrac{1}2L_{\kappa}(U^{n,k+1}+U^{n,k})-g_{\kappa}(U^{n,k})} 
    \end{align}
     for $1\le k\le s_{\mathrm{I}}$,  where the elements $ d_{k,\ell}$ are defined by $d_{k,\ell}(z):=0$ for $\ell> k$, and
    \begin{align}\label{Def: Differential Matrix DII}     
    d_{k,\ell}(z):=\underline{\theta}_{k,\ell}-z\sum_{j=\ell}^k
    \sum_{i=j}^{k}\underline{\theta}_{k,i}\underline{a}_{i+1,j+1}+z-\frac{z}2\delta_{k,\ell}\quad\text{for $1\le \ell\le k\le s_{\mathrm{I}}$}.
    \end{align}
    The associated lower triangular matrix $D:=(d_{k,\ell})_{s_{\mathrm{I}}\times s_{\mathrm{I}}}$  is called the differentiation  matrix.
    
    \hspace{5mm}Let $E_{s_{\mathrm{I}}}:=(1_{i\ge j})_{s_{\mathrm{I}}\times s_{\mathrm{I}}}$ be the lower triangular matrix full of element 1, and let $I_{s_{\mathrm{I}}}$ be the  identity matrix of the same size as $A_{\mathrm{I}}$. One has
    $$(\underline{\hat{a}}_{i+1,j})_{s_{\mathrm{I}}\times s_{\mathrm{I}}}=E_{s_{\mathrm{I}}}^{-1}A_{\mathrm{E}},\quad    
    (\underline{\theta}_{k,\ell})_{s_{\mathrm{I}}\times s_{\mathrm{I}}}=A_{\mathrm{E}}^{-1}E_{s_{\mathrm{I}}} \quad\text{and}\quad (\underline{a}_{i+1,j+1})_{s_{\mathrm{I}}\times s_{\mathrm{I}}}=E_{s_{\mathrm{I}}}^{-1}A_{\mathrm{I}}. $$ 
    Thus the above differentiation  matrix $D$ defined in 
    \eqref{Def: Differential Matrix DII} can be formulated as
    \begin{align}\label{Def: Differential Matrix D}     
      D(z)=D_{\mathrm{E}}-D_{\mathrm{EI}}z\quad\text{with}\quad D_{\mathrm{E}}:=A_{\mathrm{E}}^{-1}E_{s_{\mathrm{I}}},\quad D_{\mathrm{EI}}:=A_{\mathrm{E}}^{-1}A_{\mathrm{I}}E_{s_{\mathrm{I}}}-E_{s_{\mathrm{I}}}+\tfrac{1}2I_{s_{\mathrm{I}}}.
    \end{align}
    Always, if the symmetric part $\mathcal{S}(D;z):=\tfrac{1}{2}\kbrat{D(z)+D(z)^T}$ is positive (semi-)definite, we say that the matrix $D(z)$ is positive (semi-)definite.
      
    
    \item[(3)] \underline{Establish stage energy dissipation law}: this process is standard and we have the following result,
    which simulates the original energy dissipation law \eqref{problem: energy dissipation law} at all stages.    
    \begin{theorem}\label{thm: energy stability} If the differentiation  matrix $D(z)$ in \eqref{Def: Differential Matrix D} is positive (semi-)definite for $z\le0$, 
      the IERK methods \eqref{Scheme: general IERK stabilized2} with the stabilized parameter $\kappa\ge2\ell_g$ preserve 
      the original energy dissipation law \eqref{problem: energy dissipation law} 
      at all stages, 
      \begin{align}\label{thmResult: stage energy laws}     
        E[U^{n,j+1}]-E[U^{n,1}]\le&\,\frac1{\tau}\sum_{k=1}^{j}\myinnerB{M_h^{-1}\delta_{\tau}U^{n,k+1},
          \sum_{\ell=1}^{k}d_{k,\ell}({\tau}M_hL_{\kappa})\delta_{\tau}U^{n,\ell+1}}
      \end{align}
      for $1\le j\le s_{\mathrm{I}}$, and in particular, by taking $j:=s_{\mathrm{I}}$,
      \begin{align*}      
        E[u_h^{n}]-E[u_h^{n-1}]\le
        \frac1{\tau}\sum_{k=1}^{s_{\mathrm{I}}}\myinnerB{M_h^{-1}\delta_{\tau}U^{n,k+1},
          \sum_{\ell=1}^{k}d_{k,\ell}({\tau}M_hL_{\kappa})\delta_{\tau}U^{n,\ell+1}}
          \quad\text{for $n\ge1$.}
      \end{align*}
    \end{theorem}
    \begin{proof}Making the inner product of the equivalent form \eqref{Scheme: DOC IERK stabilized2II} 
      with $\frac1{\tau}M_h^{-1}\delta_{\tau}U^{n,k+1}$ and 
      summing $k$ from $k=1$ to $j$, one can find that
      \begin{align*}      
        \frac1{\tau}\sum_{k=1}^{j}
        &\,\myinnerB{M_h^{-1}\delta_{\tau}U^{n,k+1},
          \sum_{\ell=1}^{k}d_{k,\ell}({\tau}M_hL_{\kappa})\delta_{\tau}U^{n,\ell+1}}\\
        &\,=\sum_{k=1}^{j}\myinnerB{\delta_{\tau}U^{n,k+1},
          \tfrac{1}2L_{\kappa}(U^{n,k+1}+U^{n,k})-g_{\kappa}(U^{n,k})}
      \end{align*}
      for $1\le j\le s_{\mathrm{I}}$. Lemma \ref{lemma: origional energy derivation} yields
      the following energy dissipation law at each stage
      \begin{align*}      
        E[U^{n,j+1}]-E[U^{n,1}]
        -&\,\frac1{\tau}\sum_{k=1}^{j}\myinnerB{M_h^{-1}\delta_{\tau}U^{n,k+1},
          \sum_{\ell=1}^{k}d_{k,\ell}({\tau}M_hL_{\kappa})\delta_{\tau}U^{n,\ell+1}}\le0
      \end{align*}       
      for $1\le j\le s_{\mathrm{I}}$. It completes the proof.     
    \end{proof}
    
    For $1\le j\le s_{\mathrm{I}}$, let $D_{j}:=D[1:j,1:j]$ be the $j$-th sequential sub-matrix of the differentiation matrix $D(z)$.    
    The above stage energy dissipation law \eqref{thmResult: stage energy laws}
    can be formulated as 
    \begin{align}\label{thmResult: stage energy laws Matrix}      
      E[U^{n,j+1}]-E[U^{n,1}]\le\frac1{\tau}\myinnerB{M_h^{-1}\delta_{\tau}\vec{U}_{n,j+1},
        D_{j}({\tau}M_hL_{\kappa})\delta_{\tau}\vec{U}_{n,j+1}}
      \quad\text{for $1\le j\le s_{\mathrm{I}}$.}
    \end{align}    
   The involved vector  $\delta_{\tau}\vec{U}_{n,j+1}:=(\delta_{\tau}U^{n,2},\delta_{\tau}U^{n,3},\cdots,\delta_{\tau}U^{n,j+1})^T$ and the matrix $D(\tau M_h L_{\kappa})$ can be formulated as
    \begin{align}\label{Def: Matrix D Kronecker product}     
    	D(\tau M_h L_{\kappa}) =D_{\mathrm{E}} \otimes I+D_{\mathrm{EI}} \otimes (-\tau M_h L_{\kappa}),
    \end{align}
    where $I$ is the identity matrix of the same size as $L_{\kappa}$ and $\otimes$ represents the Kronecker product.
    Since the matrix $-\tau M_h L_{\kappa}$ is  symmetric and positive semi-definite, the properties of Kronecker product and the formula \eqref{Def: Matrix D Kronecker product} arrive at the following corollary. 

\begin{corollary}\label{corol: energy stability}
  If the two $s_{\mathrm{I}}\times s_{\mathrm{I}}$ matrices $D_{\mathrm{E}}$ and $D_{\mathrm{EI}}$ are positive (semi-)definite, the results of Theorem \ref{thm: energy stability} are valid.
\end{corollary}

    \end{description}
    
    \subsection{Average dissipation rate}
    Theorem \ref{thm: energy stability} shows that the IERK method \eqref{Scheme: general IERK stabilized2} is unconditionally energy stable if the differentiation  matrix $D(\tau M_h L_{\kappa})$ is positive semi-definite, that is, all eigenvalues of the symmetric part $\mathcal{S}(D;\tau M_h L_{\kappa})$
    are nonnegative. A necessary condition is that the average eigenvalue  is nonnegative, $$\overline{\lambda}\kbrab{\mathcal{S}(D;\tau M_h L_{\kappa})}=\frac1{s_{\mathrm{I}}m_{L_{\kappa}}}\mathrm{tr}\kbrab{D(\tau M_h L_{\kappa})}\ge0,$$    
    where $m_{L_{\kappa}}$ is the size of $L_{\kappa}$ and $\mathrm{tr}\bra{D}$ is the trace of $D$. By comparing the discrete energy dissipation law \eqref{thmResult: stage energy laws Matrix} with the continuous counterpart \eqref{problem: energy dissipation law}, the overall 
    energy dissipation rate of the  energy $E[u_h^{n}]$ could be roughly estimated by the average  
    eigenvalue. Following the idea in \cite{LiaoWang:2024arxiv}, we will use the \textit{average dissipation rate}, defined by
    \begin{align}\label{def: numerical rate}      
      \mathcal{R}:=\overline{\lambda}\kbrab{\mathcal{S}(D;\tau M_h L_{\kappa})}=\frac1{s_{\mathrm{I}}m_{L_{\kappa}}}\mathrm{tr}\kbrab{D(\tau M_h L_{\kappa})},
    \end{align} 
    to examine the energy dissipation behaviors of IERK methods. By using the definitions \eqref{Def: Differential Matrix DII} 
    and \eqref{eq: orthogonal procedureII}, one can compute the diagonal elements of the matrices $D_{\mathrm{E}}$ and $D_{\mathrm{EI}}$,
     \begin{align*}   
      (D_{\mathrm{E}})_{k,k}=&\,\underline{\theta}_{k,k}=\tfrac1{\underline{\hat{a}}_{k+1,k}}=\tfrac1{\hat{a}_{k+1,k}},\\
      (D_{\mathrm{EI}})_{k,k}=&\,\underline{\theta}_{k,k}\underline{a}_{k+1,k+1}-\tfrac{1}2
      =\tfrac{\underline{a}_{k+1,k+1}}{\underline{\hat{a}}_{k+1,k}}-\tfrac{1}2=\tfrac{a_{k+1,k+1}}{\hat{a}_{k+1,k}}-\tfrac12\quad\text{for $1\le k\le s_{\mathrm{I}}$.}
    \end{align*}
    Thus, by using the Kronecker product form \eqref{Def: Matrix D Kronecker product} and the following property
    $$\mathrm{tr}\kbrab{D_{\mathrm{E}}\otimes I+D_{\mathrm{EI}}\otimes (-\tau M_h L_{\kappa})} 
    = \mathrm{tr}(D_{\mathrm{E}})\mathrm{tr}(I)+\mathrm{tr}(D_{\mathrm{EI}}) \mathrm{tr}(-\tau M_h L_{\kappa}),$$
    it is easy to obtain the following result.
    \begin{lemma}\label{lemma: average dissipation rate}
      If the two matrices $D_{\mathrm{E}}$ and $D_{\mathrm{EI}}$ are positive (semi-)definite, then the average dissipation rate of the IERK method \eqref{Scheme: general IERK stabilized2} is nonnegative, that is,
      \begin{align}\label{def: numerical rate Formula}    
        \mathcal{R}=\frac1{s_{\mathrm{I}}}\mathrm{tr}(D_{\mathrm{E}})+\frac1{s_{\mathrm{I}}}\mathrm{tr}(D_{\mathrm{EI}})\tau\overline{\lambda}_{\mathrm{ML}}=\frac1{s_{\mathrm{I}}}\sum_{k=1}^{s_{\mathrm{I}}}\frac1{\hat{a}_{k+1,k}}
        +\frac{1}{s_{\mathrm{I}}}\sum_{k=1}^{s_{\mathrm{I}}}
        \braB{\frac{a_{k+1,k+1}}{\hat{a}_{k+1,k}}-\frac12}\tau\overline{\lambda}_{\mathrm{ML}}\ge0,  
      \end{align}   
      where the average eigenvalue $\overline{\lambda}_{\mathrm{ML}}=\overline{\lambda}_{\mathrm{ML}}(\kappa,h)>0$ of the symmetric, positive definite matrix $- M_h L_{\kappa}$ is determined by the stabilized parameter $\kappa$, the method of spatial approximation and the grid spacing $h$  as well.     
    \end{lemma}
    
  Since the average dissipation rate \eqref{def: numerical rate} is defined via  the similarity between the energy dissipation law \eqref{problem: energy dissipation law} and the discrete counterpart \eqref{thmResult: stage energy laws Matrix}, we say that an IERK method is a ``good" candidate to preserve the original energy dissipation law \eqref{problem: energy dissipation law} if the average dissipation rate $\mathcal{R}\ge0$ and is as close to 1 as possible within properly large range of $\tau\overline{\lambda}_{\mathrm{ML}}$. The specific expression of $\overline{\lambda}_{\mathrm{ML}}$ will depend on the phase field model involved; but always,  the larger the value of  $\kappa$, the greater the value of $\overline{\lambda}_{\mathrm{ML}}$, while the smaller the spatial step $h$, the greater the value of $\overline{\lambda}_{\mathrm{ML}}$.
  At the same time, the values of  $\frac1{s_{\mathrm{I}}}\mathrm{tr}(D_{\mathrm{E}})$ and $\frac1{s_{\mathrm{I}}}\mathrm{tr}(D_{\mathrm{EI}})$ are entirely determined by the IERK method itself, regardless of the time-step size $\tau$ and the spatial step $h$ adopted by potential users. 
 Throughout this paper, we aim to design some ``good" IERK methods by minimizing the values of $\frac1{s_{\mathrm{I}}}\mathrm{tr}(D_{\mathrm{EI}})$ and $\frac1{s_{\mathrm{I}}}\mathrm{tr}(D_{\mathrm{E}})$ as far as possible to make the associated average dissipation rate $\mathcal{R}$ as close to 1 as possible within a given number of implicit stages.

Obviously, the average dissipation rate $\mathcal{R}$ is only a rough indicator for describing discrete energy behavior of IERK method, that is, it would be a qualitative rather than quantitative indicator, since it ignores the computational accuracy of the algorithm.
As remarked in \cite{LiaoWang:2024arxiv}, one can
not expect that the long-time dynamics of the energy can be understood by such a scalar alone.
Nonetheless, Lemma \ref{lemma: average dissipation rate} provides us a very simple criterion to evaluate the overall energy dissipation rate of an IERK method so that one can choose proper parameters in some parameterized IERK methods or compare different IERK methods with the same time accuracy.

As an example, consider the first-order IERK (IERK1) method \cite[p. 383]{HundsdorferVerwer:2003book}, 
\begin{equation}\label{Scheme: IERK1}
 \begin{array}{c|c}
  \mathbf{c} & A \\
  \hline\\[-10pt]   & \mathbf{b}^T
 \end{array}
 =\begin{array}{c|cc}
  0 & 0 &  0  \\
  1 & 1-\theta & \theta \\
  \hline  & 1-\theta & \theta
\end{array}\;,\quad \begin{array}{c|c}
\hat{\mathbf{c}} & \widehat{A} \\
\hline\\[-10pt]   & \hat{\mathbf{b}}^T
\end{array}
=\begin{array}{c|cc}
  0 & 0 & 0   \\
  1 & 1 & 0 \\
  \hline  & 1 & 0 
  \end{array}\;.
\end{equation}
Obviously, $A_{\mathrm{I}}=(\theta)$ and $A_{\mathrm{E}}=(1)$ such that the differentiation matrix $D^{(1)}(z)=1-z\bra{\theta-\tfrac{1}2}$ is positive semi-definite for $z\le0$ provided the weighted parameter $\theta\ge\frac12$ according to Corollary \ref{corol: energy stability}.
Lemma \ref{lemma: average dissipation rate} says that the average dissipation rate
\begin{align}\label{AverRate: IERK1}
  \mathcal{R}^{(1)}(\theta)=1+\bra{\theta-\tfrac{1}2}\tau\overline{\lambda}_{\mathrm{ML}}\quad\text{for $\theta\ge\tfrac12$}.
\end{align}
To enlarge the choice of $\tau\overline{\lambda}_{\mathrm{ML}}$, one can choose $\theta=1/2$ such that the average dissipation rate, $\mathcal{R}^{(1)}(\tfrac12)=1$, is independent of  $\tau\overline{\lambda}_{\mathrm{ML}}$. We say that the stabilized Crank-Nicolson type scheme
\begin{align*}
  \delta_{\tau}u_h^{n}    
  =\tau M_h\kbra{\tfrac{1}2L_{\kappa}(u_h^{n}+u^{n-1})-g_{\kappa}(u^{n-1})}\quad \text{for $n\ge1$,}  
\end{align*}
is a ``good" candidate to preserve the original energy dissipation law \eqref{problem: energy dissipation law} unconditionally.
    
Here and hereafter, the superscript $(p)$ is always used to indicate the formal order of the method, that is to say,  $D^{(p)}$ and $\mathcal{R}^{(p)}$ denote the associated differential matrix and the average dissipation rate, respectively, of a $p$-th order IERK method.

\section{Discrete energy laws of IERK2 methods}\label{sec: IERK2}
\setcounter{equation}{0}

\subsection{Lobatto-type IERK2 methods}

Second-order methods require two implicit stages, $s_I=2$. Consider the 3-stage IERK methods that satisfy the canopy node condition and the two order conditions for the first-order accuracy,
\begin{equation*}
  \begin{array}{c|c}
    \mathbf{c} & A \\
    \hline\\[-10pt]   & \mathbf{b}^T
  \end{array}
  = \begin{array}{c|ccc}
    0 & 0 &  &      \\
    c_{2} & c_2-a_{22} & a_{22} &   \\
    1 & 1-a_{32}-a_{33} & a_{32} & a_{33}    \\
    \hline\\[-10pt]   & 1-a_{32}-a_{33} & a_{32} & a_{33} 
  \end{array}\;,\quad
\begin{array}{c|c}
  \hat{\mathbf{c}} & \widehat{A} \\
  \hline\\[-10pt]   & \hat{\mathbf{b}}^T
\end{array}
=\begin{array}{c|ccc}
  0 & 0 &  &     \\
  c_{2} & c_2 & 0 &   \\
  1 & 1-\hat{a}_{32} &\hat{a}_{32} & 0   \\
  \hline\\[-10pt]   & 1-\hat{a}_{32} &\hat{a}_{32} & 0
\end{array}\;.
\end{equation*}

We should determine five independent coefficients.  
By the two order conditions for second-order accuracy, see the second line in Table \ref{table: order condition}, 
one has three independent coefficients to be determined. In detail, 
from the stand-alone condition for explicit part, $\hat{\mathbf{b}}^T\mathbf{c}=\tfrac12$, one has  $\hat{a}_{32}=\tfrac1{2c_2}$. From the stand-alone condition for implicit part, $\mathbf{b}^T\mathbf{c}=\tfrac12$, one has
$a_{32}c_2+a_{33} =\tfrac1{2}$ such that 
$a_{32}=\frac{1-2a_{33}}{2c_2}$. According to Lemma \ref{lemma: average dissipation rate} and the formula \eqref{def: numerical rate Formula}, we always choose $$a_{22}=2c_2^2a_{33}>0$$ such that the trace $\mathrm{tr}\brat{D_{\mathrm{EI}}^{(2,1)}}$ can attain the minimum value, that is,
\begin{align}\label{IREK2: trace D1}
  \mathrm{tr}\kbrab{D_{\mathrm{EI}}^{(2,1)}(c_2,a_{33})}=\frac{a_{22}}{c_2}+2c_2a_{33}-1
  =4c_2a_{33}-1\quad\text{for $a_{33},c_2>0$.}
\end{align}
We obtain the following two-parameter class of second-order IERK (IERK2-1) methods with the associated Butcher tableaux, also see \cite{LiuZou:2006},
\begin{equation}\label{Scheme: IERK2-two parameters}
\begin{array}{c|ccc}
    0 & 0 &  &      \\
    c_{2} & c_2-2c_2^2a_{33} & 2c_2^2a_{33} &   \\[2pt]
    1 & 1-\frac{1}{2c_2}+\frac{a_{33}(1-c_2)}{c_2} & \frac{1-2a_{33}}{2c_2} & a_{33}    \\
    \hline\\[-10pt]   & 1-\frac{1}{2c_2}+\frac{a_{33}(1-c_2)}{c_2} & \frac{1-2a_{33}}{2c_2} & a_{33} 
  \end{array}\;,\quad
  \begin{array}{c|ccc}
    0 & 0 &  &     \\
    c_{2} & c_2 & 0 &   \\[2pt]
    1 & 1-\tfrac1{2c_2} &\tfrac1{2c_2} & 0   \\
    \hline\\[-10pt]   & 1-\tfrac1{2c_2} &\tfrac1{2c_2} & 0
  \end{array}\;.
\end{equation}

Now we are to determine the two independent coefficients $c_2$ and $a_{33}$ such that the IERK2-1 methods \eqref{Scheme: IERK2-two parameters} are ``good" candidates to preserve 
the original energy dissipation law \eqref{problem: energy dissipation law} unconditionally. By the definition \eqref{Def: Differential Matrix D}, one has
\begin{align*}
  D_{\mathrm{E}}^{(2,1)}(c_2):=A_{\mathrm{E}}^{-1}E_2=\begin{pmatrix}
    \frac{1}{c_2}& 0  \\[4pt] 
    2c_2+\frac1{c_{2}}-2& 2c_2
  \end{pmatrix}.
\end{align*}
The determinant of the symmetric part 
$$\mathrm{Det}\kbrab{\mathcal{S}(D_{\mathrm{E}}^{(2,1)};c_2)}
=\brab{c_2+\sqrt{2}-1+\tfrac1{2c_2}}\brab{-c_2+\sqrt{2}+1-\tfrac1{2c_2}}.$$
It is easy to check that 
\begin{align}\label{cond: IERK2-c2 condition}
  0.228788\approx\tfrac{1+\sqrt{2}-\sqrt{1+2 \sqrt{2}}}{2}\le c_2\le \tfrac{1+\sqrt{2}+\sqrt{1+2 \sqrt{2}}}{2}\approx 2.18543
\end{align}
is sufficient to ensure the positive semi-definiteness of  $D_{\mathrm{E}}^{(2,1)}(c_2)$.  
Now compute the lower triangular matrix $D_{\mathrm{EI}}^{(2,1)}$ defined by \eqref{Def: Differential Matrix D},
\begin{align*}
  D_{\mathrm{EI}}^{(2,1)}(c_2,a_{33}):=A_{\mathrm{E}}^{-1}A_{\mathrm{I}}E_2-E_2+\tfrac{1}2I_2=\begin{pmatrix}
    2c_2a_{33}-\frac12& 0  \\[4pt]  
    -2 a_{33}(2 c_2^2-2c_2+1) & 2c_2a_{33}-\frac12
  \end{pmatrix}.
\end{align*}
The positive semi-definiteness of  $D_{\mathrm{EI}}^{(2,1)}(c_2,a_{33})$ 
requires $a_{33}\ge\frac1{4c_2}$
and
\begin{align*}
  \mathrm{Det}\kbrab{\mathcal{S}(D_{\mathrm{EI}}^{(2,1)};c_2,a_{33})}
=&-a_{33}^2 \bra{2c_2^2+1}(2 c_2^2-4 c_2+1)-2 a_{33}c_2+\frac{1}{4}\\
=&\brat{2c_2^2+1}(4c_2-2 c_2^2-1)\braB{a_{33}-\tfrac1{ 2(2c_2^2+1)}}\braB{a_{33}-\tfrac1{2 (4c_2-2 c_2^2-1)}}\ge0,
\end{align*}
in which we set $4c_2-2 c_2^2-1>0$, or $1-\frac{\sqrt{2}}2<c_2<1+\frac{\sqrt{2}}2$. Otherwise, the fact $a_{33}\ge\frac1{4c_2}>\frac1{ 4 c_2^2+2}$ arrives at the negative determinant, $\mathrm{Det}\kbrab{\mathcal{S}(D_{\mathrm{EI}}^{(2,1)};c_2,a_{33})}<0$. Thus, by using the restriction \eqref{cond: IERK2-c2 condition} and the fact $\frac{1}{4c_2^2+2}\le \frac1{4c_2}\le \frac1{2 (4c_2-2 c_2^2-1)}$, one has the following conditions for the  coefficient $a_{33}$,
\begin{align}\label{cond: IERK2-c2 b3 condition}
  a_{33}\ge \frac{1}{2 (4c_2-2 c_2^2-1)}=\frac{1}{4(c_2-1+\frac{\sqrt{2}}{2})(1+\frac{\sqrt{2}}{2}-c_2)}\quad\text{for $1-\tfrac{\sqrt{2}}2<c_2<1+\tfrac{\sqrt{2}}2$.}
\end{align}
By using Corollary \ref{corol: energy stability} and Lemma \ref{lemma: average dissipation rate}, we have the following theorem.
\begin{theorem}\label{thm: IERK2-two parameters} 
  In simulating the gradient flow system \eqref{problem: stabilized version} with  $\kappa\ge2\ell_g$, the two-parameter IERK2-1 methods \eqref{Scheme: IERK2-two parameters} with the parameter setting  \eqref{cond: IERK2-c2 b3 condition} preserve
  the original energy dissipation law \eqref{problem: energy dissipation law} unconditionally
  at all stages in the sense that
  \begin{align*}    
    E[U^{n,j+1}]-E[U^{n,1}]\le\frac1{\tau}\myinnerB{M_h^{-1}\delta_{\tau}\vec{U}_{n,j+1},
      D_{j}^{(2,1)}(c_2,a_{33};{\tau}M_hL_{\kappa})\delta_{\tau}\vec{U}_{n,j+1}}
    \quad\text{for $1\le j\le 2$,}
  \end{align*}
  where the  differentiation  matrix $D^{(2,1)}$ is defined by  $$D^{(2,1)}(c_2,a_{33};{\tau}M_hL_{\kappa}):=D_{\mathrm{E}}^{(2,1)}(c_2)\otimes I+D_{\mathrm{EI}}^{(2,1)}(c_2,a_{33}) \otimes (-\tau M_h L_{\kappa}).$$
  The associated average dissipation rate 
  \begin{align}\label{AverRate: IERK2-two parameters}
    \mathcal{R}^{(2,1)}(c_2,a_{33})
    =\tfrac{1}{2}\mathrm{tr}\kbrab{D_{\mathrm{E}}^{(2,1)}(c_2)}+\tfrac{1}{2}\mathrm{tr}\kbrab{D_{\mathrm{EI}}^{(2,1)}(c_2,a_{33})}\tau\overline{\lambda}_{\mathrm{ML}}
    =c_2+\tfrac{1}{2c_2}+(2c_2a_{33}-\tfrac1{2})\tau\overline{\lambda}_{\mathrm{ML}}.
  \end{align}
\end{theorem}

The parameter setting  \eqref{cond: IERK2-c2 b3 condition} implies that one can minimize the trace $\mathrm{tr}\kbrab{D_{\mathrm{EI}}^{(2,1)}(c_2,a_{33})}$ by choosing  the lower bound of $a_{33}$, that is,
\begin{align*}
  \mathrm{tr}\kbrab{D_{\mathrm{EI}}^{(2,1)}(c_2,a_{33})}\ge&\, \frac{c_2}{4c_2-2 c_2^2-1}-\frac1{2}\quad\text{for $1-\tfrac{\sqrt{2}}2<c_2<1+\tfrac{\sqrt{2}}2$.}
\end{align*}
The minimum value is attained by setting $c_2=\frac{\sqrt{2}}{2}$ due to the fact $\frac{\zd}{\zd c_2}\brab{\frac{c_2}{4c_2-2 c_2^2-1}}=\frac{2 c_2^2-1}{(2 c_2^2-4 c_2+1)^2}$. 
In this case, the two-parameter IERK2-1 methods \eqref{Scheme: IERK2-two parameters} reduce into the one-parameter IERK2 (IERK2-2) methods with the parameter $a_{33}\ge\tfrac{1+\sqrt{2}}{4}$,
\begin{equation}\label{Scheme: IERK2-one parameter}
  \begin{array}{c|ccc}
    0 & 0 &  &      \\
    \frac{\sqrt{2}}2 & \frac{\sqrt{2}}2-a_{33} & a_{33}  &   \\[1pt]
    1 & \frac{\sqrt{2}-1+(2-\sqrt{2})a_{33}}{\sqrt{2}} & \frac{1-2a_{33}}{\sqrt{2}} & a_{33}     \\[1pt]
    \hline\\[-10pt]   & \frac{\sqrt{2}-1+(2-\sqrt{2})a_{33}}{\sqrt{2}} & \frac{1-2a_{33}}{\sqrt{2}} & a_{33}  
  \end{array}\;,\quad
  \begin{array}{c|ccc}
    0 & 0 &  &     \\
    \frac{\sqrt{2}}2 & \frac{\sqrt{2}}2 & 0 &   \\[1pt]
    1 & \frac{2-\sqrt{2}}2 &\frac{\sqrt{2}}2 & 0   \\[1pt]
    \hline\\[-10pt]   & \frac{2-\sqrt{2}}2 &\frac{\sqrt{2}}2 & 0
  \end{array}\;.
\end{equation}

As expected, from the perspective of average dissipation rate, the optimal IERK2-2 methods maintain a fixed diagonal ratio between the implicit and explicit parts, that is, $\frac{a_{k+1,k+1}}{\hat{a}_{k+1,k}}=\sqrt{2}a_{33}$ for $1\le k\le 2$. The associated average dissipation rate 
\begin{align}\label{AverRate: IERK2-one parameter}
	\mathcal{R}^{(2,2)}(\tfrac{\sqrt{2}}{2},a_{33})
	=\sqrt{2}+\sqrt{2}(a_{33}-\tfrac{\sqrt{2}}{4})\tau\overline{\lambda}_{\mathrm{ML}}\quad\text{for $a_{33}\ge\tfrac{1+\sqrt{2}}{4}$}.
\end{align}
Note that, the 3-stage Lobatto-type IERK2 method in \cite[Section 2]{CooperSayfy:1983} would be not suitable for gradient flow problems because the associated differentiation matrix is not positive definite.

\subsection{Radau-type IERK2 methods}
For completeness, we also briefly consider the Radau-type DIRK methods, also known as ARS-type IERK methods \cite{AscherRuuthSpiteri:1997,FuTangYang:2024,ShinLeeLee:2017}, for the implicit part by assuming $\mathbf{a}_1=\mathbf{0}$. In this case, one has three independent coefficients,  
see the associated Butcher tableaux as follows,
\begin{equation*}
  \begin{array}{c|c}
    \mathbf{c} & A \\
    \hline\\[-10pt]   & \mathbf{b}^T
  \end{array}
  = \begin{array}{c|ccc}
    0 & 0 &  &      \\
    c_{2} & 0 & c_{2} &   \\
    1 & 0 & 1-a_{33} & a_{33}    \\
    \hline\\[-10pt]   & 0 & 1-a_{33} & a_{33} 
  \end{array}\;,\quad
  \begin{array}{c|c}
    \hat{\mathbf{c}} & \widehat{A} \\
    \hline\\[-10pt]   & \hat{\mathbf{b}}^T
  \end{array}
  =\begin{array}{c|ccc}
    0 & 0 &  &     \\
    c_{2} & c_2 & 0 &   \\
    1 & 1-\hat{a}_{32} &\hat{a}_{32} & 0   \\
    \hline\\[-10pt]   & 1-\hat{a}_{32} &\hat{a}_{32} & 0
  \end{array}\;.
\end{equation*}
From $\hat{\mathbf{b}}^T\mathbf{c}=\tfrac12$, one has  $\hat{a}_{32}=\tfrac1{2c_2}$. The condition $\mathbf{b}^T\mathbf{c}=\tfrac12$ leads to
$a_{33}=\frac{1-2c_2}{2(1-c_2)}$. We obtain the following $c_2$-parameterized class of Radau-type IERK (IERK2-Radau) methods \cite{AscherRuuthSpiteri:1997}
\begin{equation}\label{Scheme: IERK2-Radau-one parameter}
  \begin{array}{c|ccc}
    0 & 0 &  &      \\
    c_{2} & 0 & c_{2} &   \\
    1 & 0 & \frac{1}{2(1-c_2)} & \frac{1-2c_2}{2(1-c_2)}    \\
    \hline\\[-10pt]   & 0 & \frac{1}{2(1-c_2)} & \frac{1-2c_2}{2(1-c_2)} 
  \end{array},\quad
  \begin{array}{c|ccc}
    0 & 0 &  &     \\
    c_{2} & c_2 & 0 &   \\
    1 & 1-\tfrac1{2c_2} &\tfrac1{2c_2} & 0   \\
    \hline\\[-10pt]   & 1-\tfrac1{2c_2} &\tfrac1{2c_2} & 0
  \end{array}.
\end{equation}
Obviously, the restriction \eqref{cond: IERK2-c2 condition} is sufficient to ensure the positive semi-definiteness of  $D_{\mathrm{rad},\mathrm{E}}^{(2)}$. We also apply
the definition \eqref{Def: Differential Matrix D} to find
\begin{align*}
  D_{\mathrm{rad},\mathrm{EI}}^{(2)}(c_2):=A_{\mathrm{E}}^{-1}A_{\mathrm{I}}E_2-E_2+\tfrac{1}2I_2
  =\begin{pmatrix}
    \frac12& 0  \\[4pt] 
    0& \frac{4c_2^2-3c_2+1}{2(c_2-1)}
  \end{pmatrix}.
\end{align*}
The positive semi-definiteness of  $D_{\mathrm{rad},\mathrm{EI}}^{(2)}(c_2)$ 
requires $c_2>1$. By using Corollary \ref{corol: energy stability} and Lemma \ref{lemma: average dissipation rate}, we  have the following theorem.

\begin{theorem}\label{thm: IERK2-Radau} 
  In simulating the gradient flow system \eqref{problem: stabilized version} with  $\kappa\ge2\ell_g$, the $c_2$-parameterized IERK2-Radau methods \eqref{Scheme: IERK2-Radau-one parameter} with the parameter $1<c_2\le \tfrac{1+\sqrt{2}+\sqrt{1+2 \sqrt{2}}}{2}$ preserve
  the original energy dissipation law \eqref{problem: energy dissipation law} unconditionally
  at all stages in the sense that
  \begin{align*}    
    E[U^{n,j+1}]-E[U^{n,1}]\le\frac1{\tau}\myinnerB{M_h^{-1}\delta_{\tau}\vec{U}_{n,j+1},
      D_{\mathrm{rad},j}^{(2)}(c_2;{\tau}M_hL_{\kappa})\delta_{\tau}\vec{U}_{n,j+1}}
    \quad\text{for $1\le j\le 2$,}
  \end{align*}
  where the associated differentiation  matrix $D_{\mathrm{rad}}^{(2)}$ is defined by  $$D_{\mathrm{rad}}^{(2)}(c_2;{\tau}M_hL_{\kappa}):=D_{\mathrm{rad},\mathrm{E}}^{(2)}(c_2)\otimes I+D_{\mathrm{rad},\mathrm{EI}}^{(2)}(c_2) \otimes (-\tau M_h L_{\kappa}).$$
  The associated average dissipation rate 
  \begin{align}\label{AverRate: IERK2-Radau-one parameters}
    \mathcal{R}_{\mathrm{rad}}^{(2)}(c_2)
    =\tfrac{1}{2c_2}+c_2+\brab{2c_2+1+\tfrac{1}{c_2-1}}\tau\overline{\lambda}_{\mathrm{ML}}
    \quad\text{for $1<c_2\le \tfrac{1+\sqrt{2}+\sqrt{1+2 \sqrt{2}}}{2}$}.
  \end{align}
\end{theorem}

To enlarge the possible choices of $\tau\overline{\lambda}_{\mathrm{ML}}$, we choose $c_2=1+\frac{\sqrt{2}}{2}$ such that $\mathrm{tr}\kbrab{D_{\mathrm{rad},\mathrm{EI}}^{(2)}(c_2)}$ attains the minimum value and $ \mathcal{R}_{\mathrm{rad}}^{(2)}(1+\tfrac{\sqrt{2}}{2})
=2+\brat{3+2\sqrt{2}}\tau\overline{\lambda}_{\mathrm{ML}}$.
However, it is always larger than the average dissipation rate $\mathcal{R}^{(2,2)}(\tfrac{\sqrt{2}}{2},a_{33})$ of the IERK2-2 methods \eqref{Scheme: IERK2-one parameter} for  $\tfrac{1+\sqrt{2}}{4}\le a_{33}\le\tfrac{8+7\sqrt{2}}{4}$.
Under the same stabilization strategy, the IERK2-2 methods \eqref{Scheme: IERK2-one parameter} would be more competitive than the IERK2-Radau methods \eqref{Scheme: IERK2-Radau-one parameter}, at least, from the perspective of average dissipation rate.

 To improve the average dissipation rate of Radau-type IERK2 methods, one can consider four-stage procedure such that there are nine independent coefficients to be determined. An example is the four-stage Radau-type IERK2 procedure in \cite[Section 4]{ShinLeeLee:2017},  which is constructed by imposing the simplifying setting $A_{\mathrm{I}}=A_{\mathrm{E}}$ and has a better average dissipation rate $R_{\mathrm{rad},S}^{(2)}=\tfrac{3}{2}+\tfrac{1}{2}\tau\overline{\lambda}_{\mathrm{ML}}$. It is to mention that,  the 3-stage Radau-type IERK method in \cite{AscherRuuthSpiteri:1997} would not be suitable for the gradient flow problems since it does not satisfy the assumptions of Corollary \ref{corol: energy stability}.

\subsection{Tests of IERK2 methods}

\begin{example}\label{ex: with source term}
  Consider the Cahn-Hilliard model, $\partial_tu=\partial_{xx}\brat{-\epsilon^{2}\partial_{xx} u-u+u^3}+f$, subject to the initial data $u_{0} = \sin x$ on $\Omega=(0,2\pi)$ with the interface parameter $\epsilon=0.2$. The source term $f$ is set by choosing the exact solution $u=e^{-t}\sin x$. Always, the spatial operators are approximated by the Fourier pseudo-spectral approximation with 256 grid points. 
\end{example}

%

\begin{figure}[htb!]
  \centering
  \subfigure[IERK2-1 \eqref{Scheme: IERK2-two parameters} with $c_2=1$]{
  \includegraphics[width=2in]{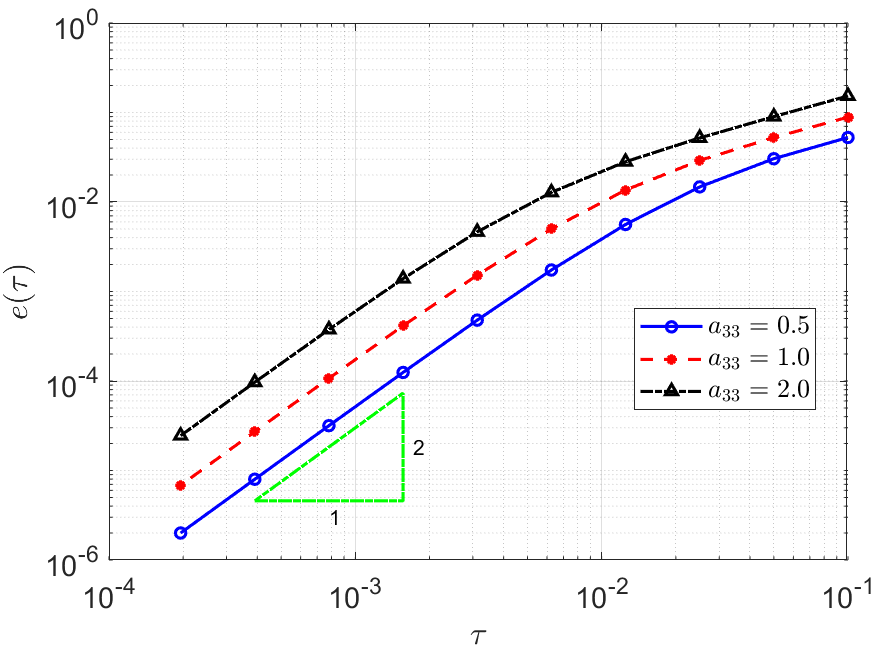}}
  \subfigure[IERK2-2 \eqref{Scheme: IERK2-one parameter}]{
  \includegraphics[width=2in]{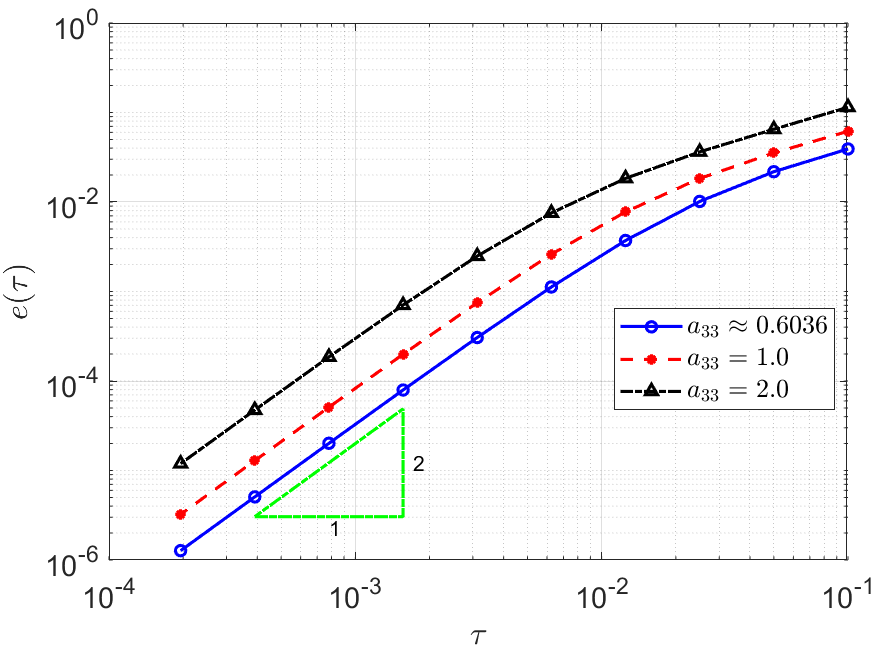}}
  \subfigure[IERK2-Radau \eqref{Scheme: IERK2-Radau-one parameter}]{
  \includegraphics[width=2in]{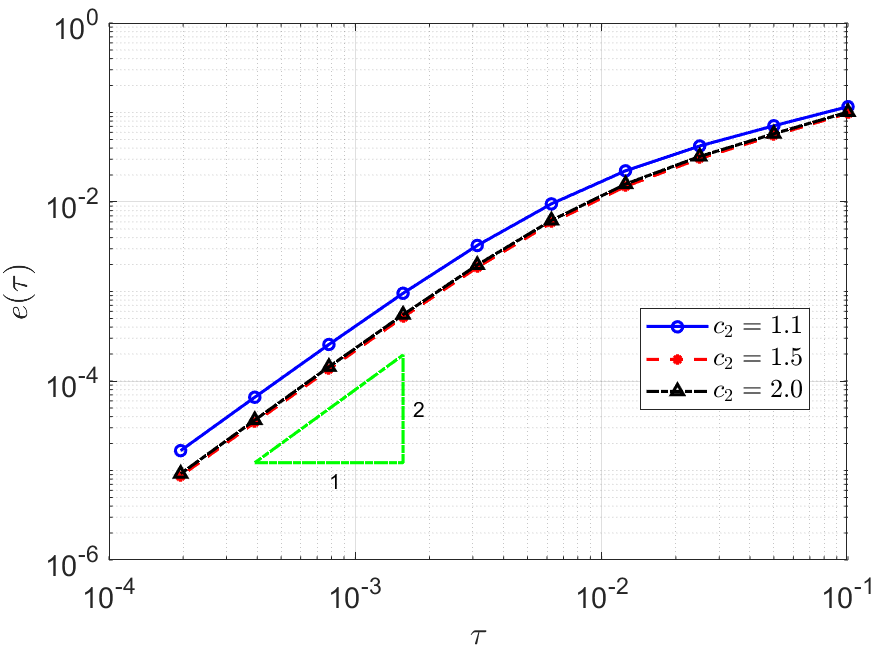}}
  \caption{Solution errors of IERK2 methods with different method parameters.}
  \label{fig: IERK2 convergence}
\end{figure}

At first, we test the convergence of IERK2-1 \eqref{Scheme: IERK2-two parameters}, IERK2-2 \eqref{Scheme: IERK2-one parameter} and IERK2-Radau \eqref{Scheme: IERK2-Radau-one parameter} methods by choosing the final time $T=1$ and the stabilized parameter $\kappa=4$. Figure \ref{fig: IERK2 convergence} lists the $L^\infty$ norm error $e(\tau):=\max_{1\le n\le N}\mynorm{u_h^{n} - u(t_n)}_\infty$ for the three classes of IERK2 methods on halving time steps $\tau=2^{-k}/10$ for $0\le k\le9$. As expected, the IERK2-1 \eqref{Scheme: IERK2-two parameters}, IERK2-2 \eqref{Scheme: IERK2-one parameter} and IERK2-Radau \eqref{Scheme: IERK2-Radau-one parameter} methods are second-order accurate in time. It seems that the different parameters for the IERK2-1 \eqref{Scheme: IERK2-two parameters} and IERK2-2 \eqref{Scheme: IERK2-one parameter} methods would arrive at different numerical precision; while the IERK2-Radau  methods \eqref{Scheme: IERK2-Radau-one parameter} with two different parameters $c_2=1.5$ and $c_2=2$ generates almost the same solution. Also, it can be observed that the error of IERK2-2 method for $a_{33}=\tfrac{1+\sqrt{2}}{4}\approx 0.6036$ is always the smallest for the same time-step size $\tau$.

\begin{example}\label{example: energy_2}
  Consider the Cahn-Hilliard model, $\partial_tu=\partial_{xx}\brat{-\epsilon^{2}\partial_{xx} u-u+u^3}$ on $\Omega=(-\pi,\pi)$ with the interface parameter $\epsilon=0.1$, subject to the initial data 
  \begin{align*}
    u_{0} =&\, \tfrac{1}{3} \tanh(2\sin x) - \tfrac{1}{10} e^{-23.5(\left |x \right |-1)^{2}} + e^{-27(\left |x \right |-4.2)^{2}} + e^{-38(\left |x \right |-5.4)^2}.
  \end{align*}
  The reference solution is generated with a small time-step size $\tau=10^{-3}$ by the IERK2-1 method \eqref{Scheme: IERK2-two parameters}  for the parameters $c_2=1$ and $a_{33}=1$.
\end{example}

\begin{figure}[htb!]
  \centering
  \subfigure[solution $u_h^N$ for $\tau=0.01,\kappa=2$]{
    \includegraphics[width=2.05in]{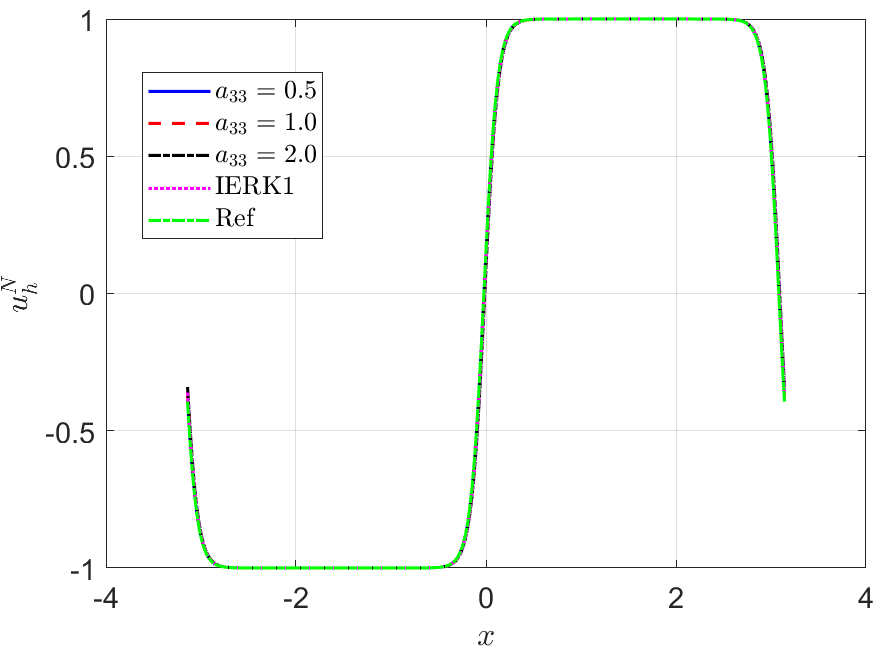}}\hspace{10mm}
  \subfigure[energy for $\tau=0.01,\kappa=2$]{
    \includegraphics[width=2.05in]{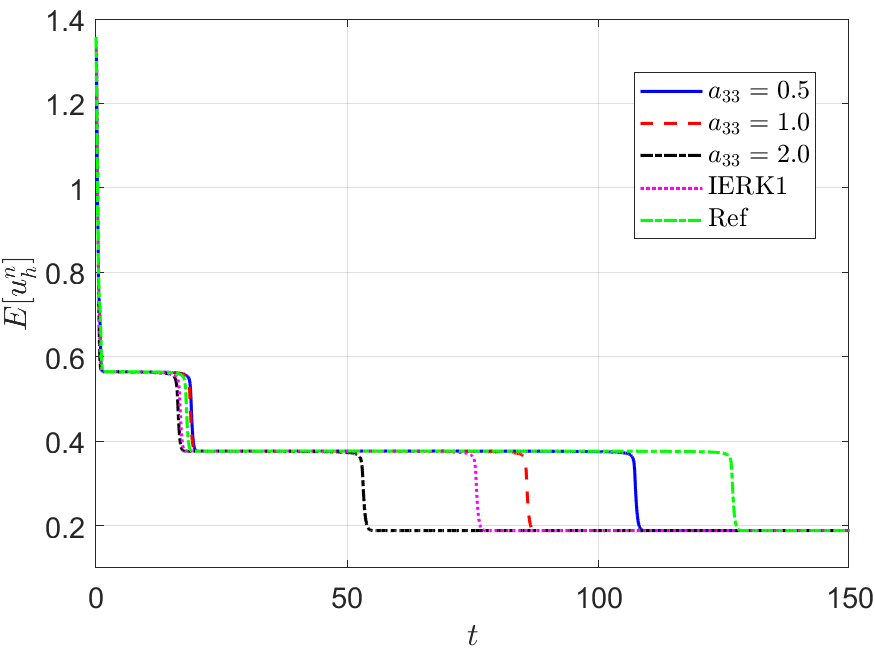}}
  \subfigure[energy for $\tau=0.05,\kappa=2$]{
    \includegraphics[width=2.05in]{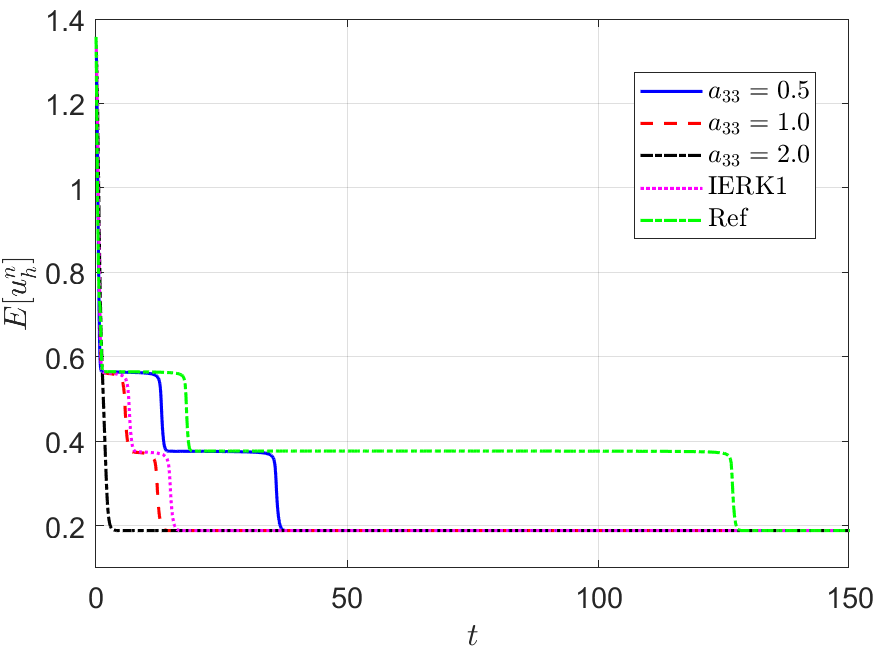}}\hspace{10mm} 
    \subfigure[energy for $\tau=0.05,\kappa=3$]{
      \includegraphics[width=2.05in]{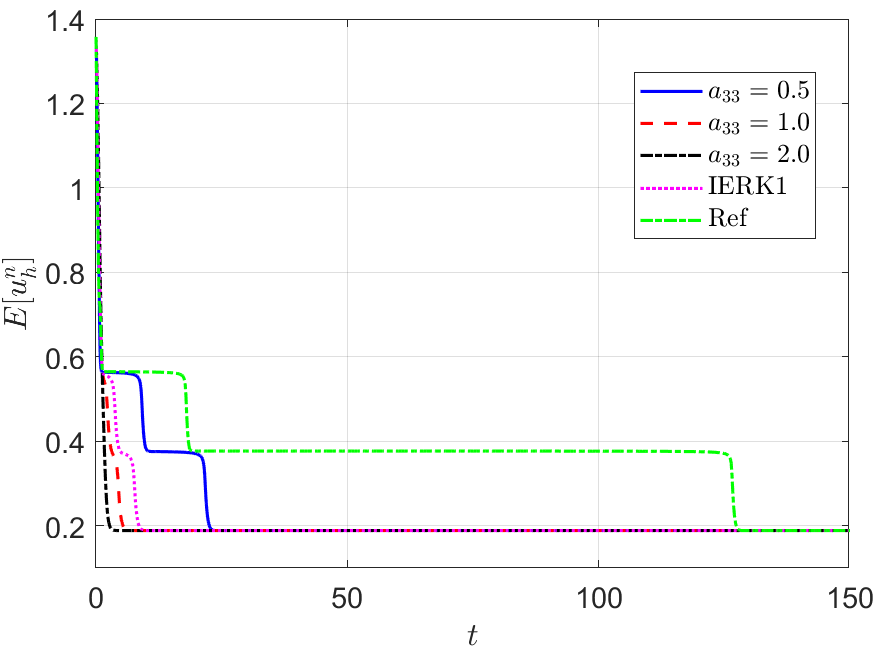}}
  \caption{Final solution and energy behaviors of IERK2-1 methods \eqref{Scheme: IERK2-two parameters} with $c_2=1$.}
  \label{fig: IERK2-I decay a33}
\end{figure}

Figure \ref{fig: IERK2-I decay a33} depicts the final solution $u_h^N$ at $T=150$ for $\tau=0.01,\kappa=2$, and the discrete energy $E[u_h^n]$ generated by the IERK2-1 methods \eqref{Scheme: IERK2-two parameters} for three different scenes: (i) $\tau=0.01,\kappa=2$; (ii) $\tau=0.05,\kappa=2$ and (iii) $\tau=0.05,\kappa=3$. 
As predicted by Theorem \ref{thm: IERK2-two parameters}, the original energy curves in all cases always monotonically decreasing. Also, the numerical dissipation rates of original energy seem quite different for three different choices of the parameter $a_{33}$. For three different discrete scenes in Figure \ref{fig: IERK2-I decay a33} (b)-(d), the energy curves for $a_{33}=0.5$ are always closest to the reference energy, while the energies for $a_{33}=2$ are the farthest away. More interestingly, these phenomena can be predicted by \eqref{AverRate: IERK2-two parameters}, or $\mathcal{R}^{(2,1)}(1,a_{33})=\tfrac{3}{2}+(2a_{33}-\tfrac1{2})\tau\overline{\lambda}_{\mathrm{ML}}$, which suggests that the discrete energy gradually deviates from the continuous energy as the method parameter $a_{33}$ increases.

It seems that, the above differences between discrete energy curves and continuous energy shown in Figure \ref{fig: IERK2-I decay a33} (b)-(d) can be explained by the different precision of numerical solutions. Actually, Figure \ref{fig: IERK2 convergence} (a) shows that the solution for the case $a_{33}=0.5$ is a bit more accurate than that for $a_{33}=2$ although both of them are second-order accurate. Maybe, this would not be the whole story. In Figure \ref{fig: IERK2-I decay a33} (b)-(d), we also include the discrete energy generated by the first-order IERK1 method \eqref{Scheme: IERK1} with $\theta=1/2$, which has the average dissipation rate $\mathcal{R}^{(1)}(\tfrac12)=1$. More surprisingly, the energy curve generated by the first-order IERK1 method with $\theta=\tfrac12$ is even closer to the reference one than some IERK2-1 schemes especially when the time-step size is properly large. Although we are unsure of the complete mechanism behind it, they suggest that the selection of method parameters in IERK methods is at least as important as the selection of different IERK methods, if not more important. Obviously, the choice $a_{33}=\frac{1}{2}$ is a good choice for IERK2-1 methods \eqref{Scheme: IERK2-two parameters}, at least, for this example.
 
\begin{figure}[htb!]
  \centering
  \subfigure[$\tau=0.01,\kappa=2$]{
    \includegraphics[width=2.05in]{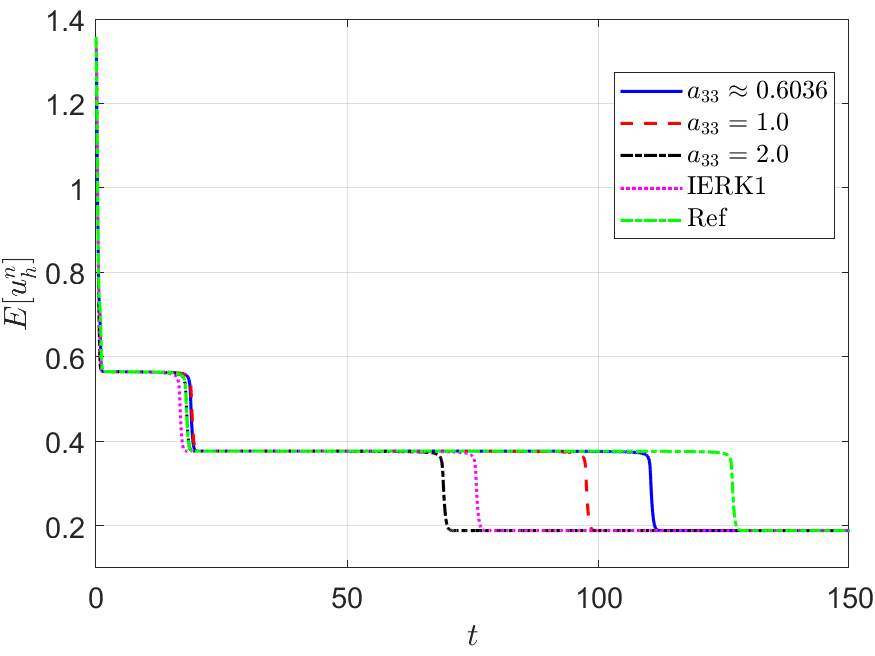}}  
  \subfigure[ $\tau=0.05,\kappa=2$]{
    \includegraphics[width=2.05in]{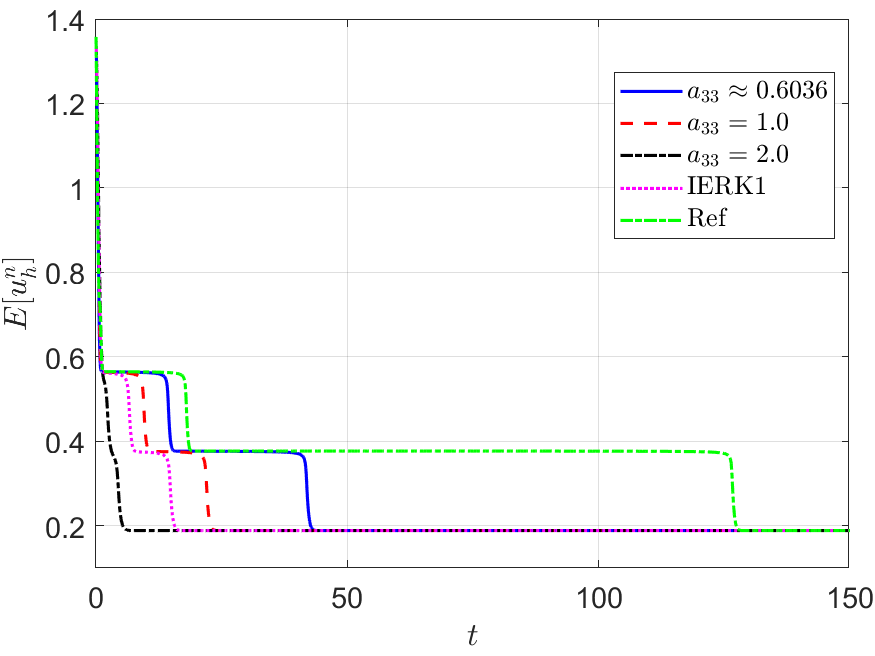}}
    \subfigure[$\tau=0.05,\kappa=3$]{
      \includegraphics[width=2.05in]{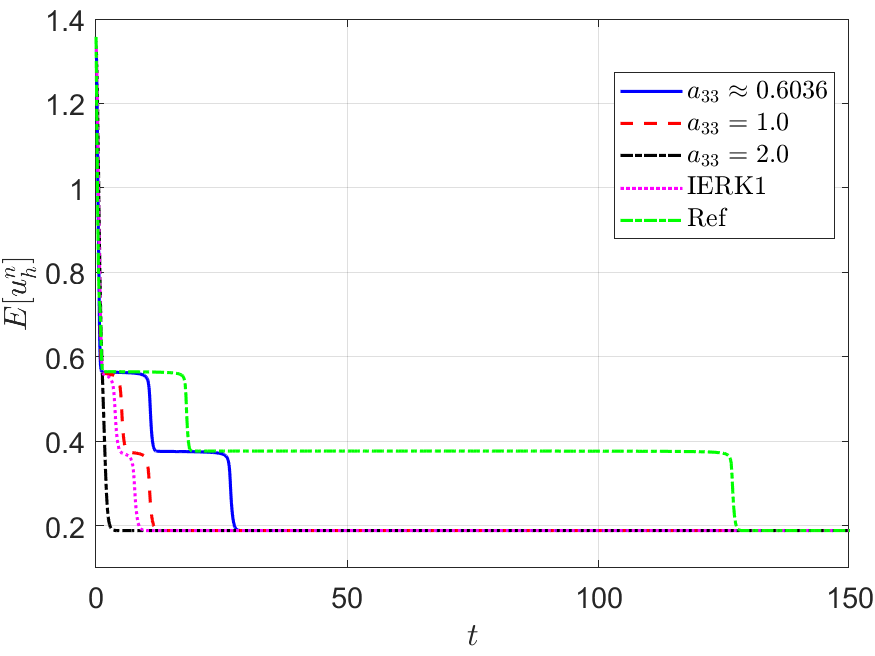}}
  \caption{Energy behaviors of $a_{33}$-parameterized IERK2-2 methods \eqref{Scheme: IERK2-one parameter}.}
  \label{fig: IERK2-2 decay a33_tau_kappa}
\end{figure}

\begin{figure}[htb!]
  \centering
  \subfigure[$\tau=0.01,\kappa=2$]{
    \includegraphics[width=2.05in]{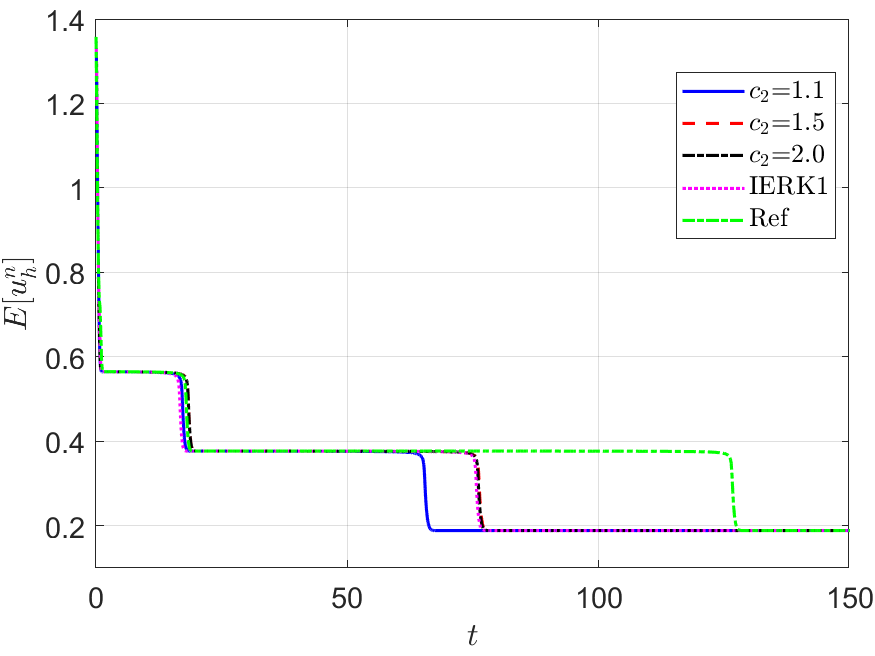}}
    \subfigure[$\tau=0.05,\kappa=2$]{
    \includegraphics[width=2.05in]{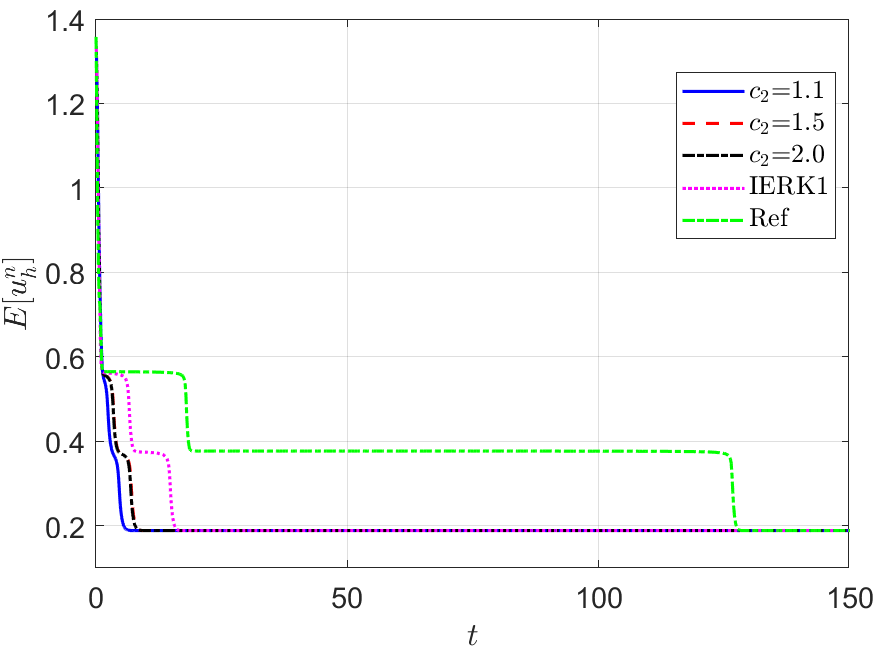}}
    \subfigure[$\tau=0.05,\kappa=3$]{
      \includegraphics[width=2.05in]{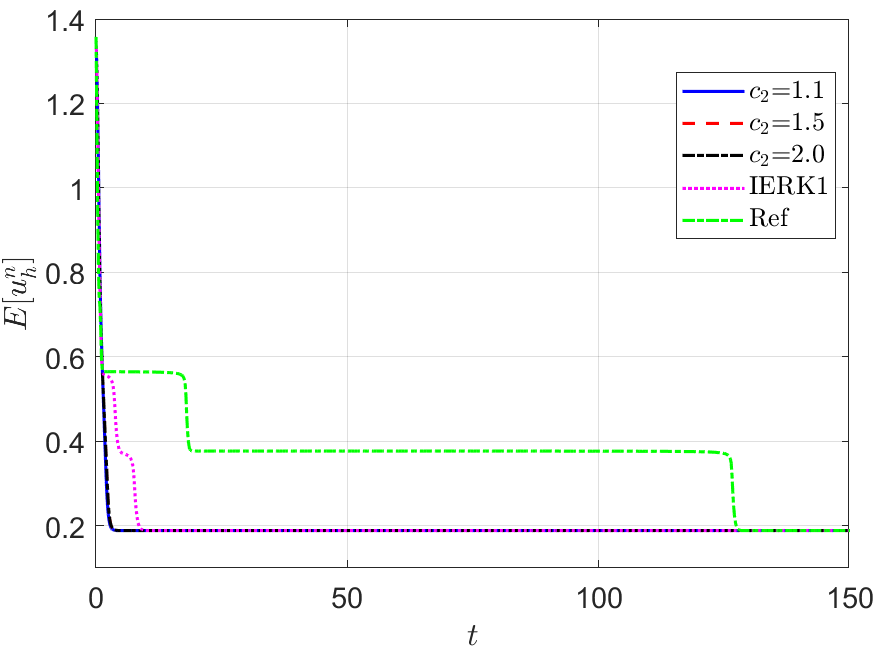}}
  \caption{Energy behaviors of IERK2-Radau methods \eqref{Scheme: IERK2-Radau-one parameter}.}
  \label{fig: IERK2-3 decay c2_tau_kappa}
\end{figure}

\begin{figure}[htb!]
  \centering
  \subfigure[$\tau=0.01,\kappa=2$]{
    \includegraphics[width=2in]{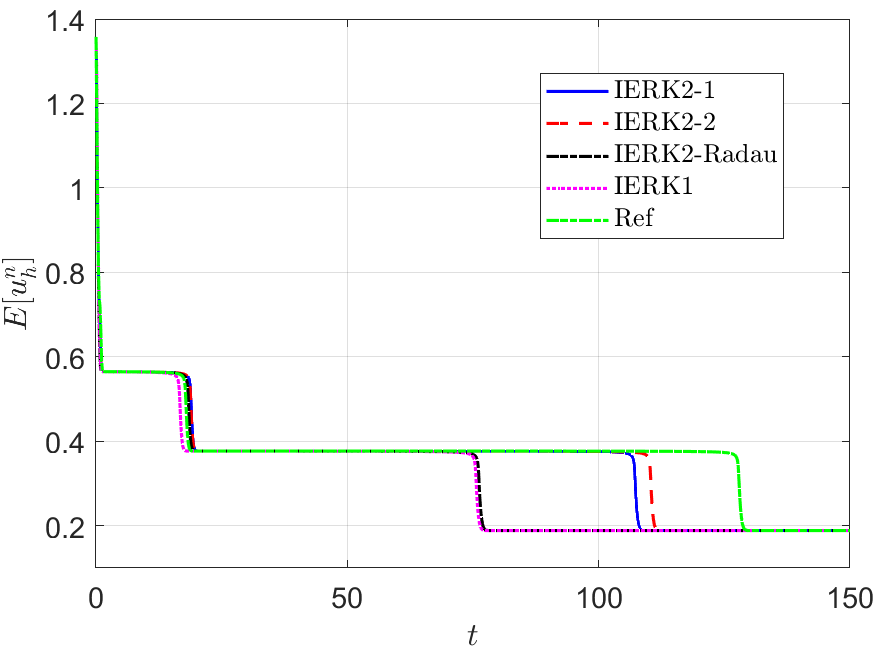}} 
  \subfigure[$\tau=0.05,\kappa=2$]{
    \includegraphics[width=2in]{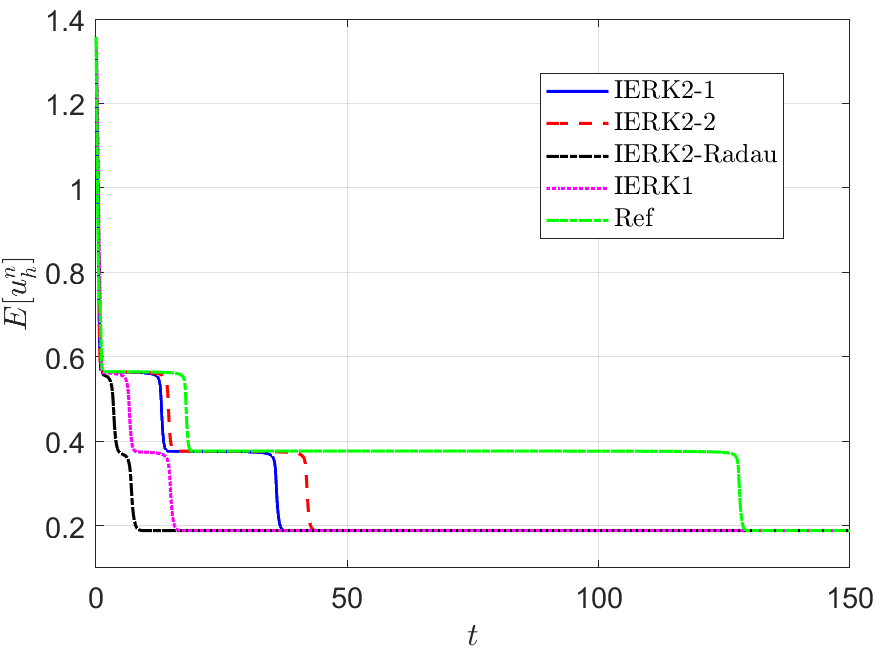}}
    \subfigure[$\tau=0.05,\kappa=3$]{
      \includegraphics[width=2in]{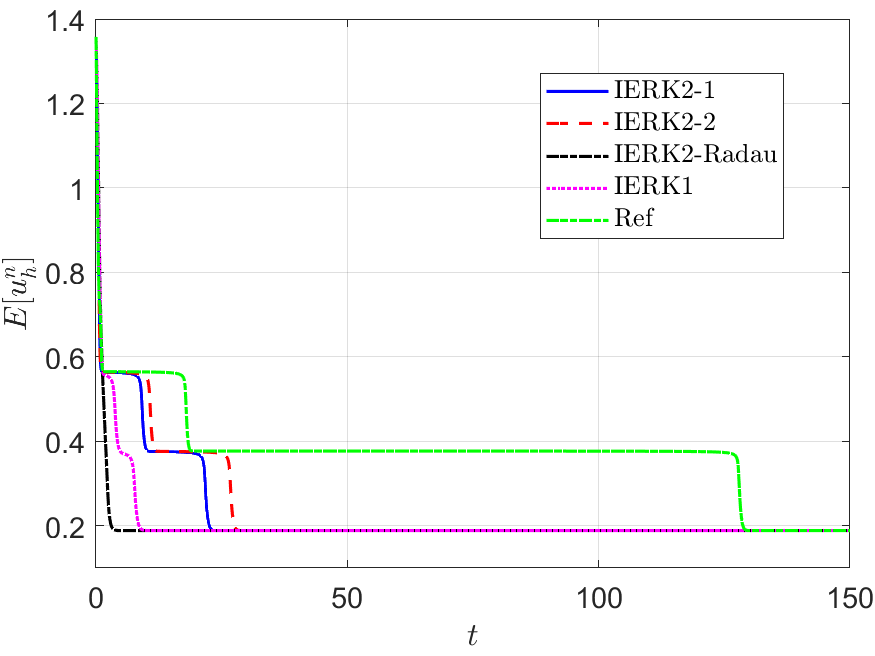}}
  \caption{Energy behaviors of IERK2-1 method for $c_2=1$ and $a_{33}=\frac1{2}$, IERK2-2 method for $a_{33}=\tfrac{1+\sqrt{2}}{4}\approx 0.6036$, and IERK2-Radau method for $c_{2}=\frac3{2}$.}
  \label{fig: IERK2 decay tau_kappa}
\end{figure}

We also run the one-parameterized IERK2-2 \eqref{Scheme: IERK2-one parameter} and IERK2-Radau \eqref{Scheme: IERK2-Radau-one parameter} methods up to the time $T=150$, and present the discrete energy curves in Figures \ref{fig: IERK2-2 decay a33_tau_kappa} and \ref{fig: IERK2-3 decay c2_tau_kappa}, respectively, for three different scenes: (i) $\tau=0.01,\kappa=2$; (ii) $\tau=0.05,\kappa=2$ and (iii) $\tau=0.05,\kappa=3$. As predicted by our theory, the discrete energies in all cases always monotonically decrease, and the discrete energies generated by the same method gradually deviate from the continuous energy as the associated average dissipation rate deviates from 1 (reminding that the average dissipation rate is an increasing function with respect to the step size $\tau$ and the parameter $\kappa$), cf. the formulas \eqref{AverRate: IERK2-one parameter} and \eqref{AverRate: IERK2-Radau-one parameters}. For this example, one can find that the case $a_{33}=\tfrac{1+\sqrt{2}}{4}\approx 0.6036$ of IERK2-2 methods \eqref{Scheme: IERK2-one parameter}, and the case $c_2=\frac{3}{2}$ of IERK2-Radau methods \eqref{Scheme: IERK2-Radau-one parameter}  seem better than other parameters for the same time-step size and the stabilized parameter. Again, from the difference between discrete energy and continuous energy (for the same time-step size), the IERK1 scheme \eqref{Scheme: IERK1} with $\theta=\frac1{2}$ is also comparable to some second-order IERK2 schemes, although it only has first-order time accuracy.

Also, in Figure \ref{fig: IERK2 decay tau_kappa}, we compare the energy behaviors of three ``good" schemes, namely, the IERK2-1 method for $c_2=1$ and $a_{33}=\frac1{2}$, the IERK2-2 method for $a_{33}=\tfrac{1+\sqrt{2}}{4}\approx 0.6036$, and the IERK2-Radau method for $c_{2}=\frac3{2}$, in previous tests shown in Figures \ref{fig: IERK2-I decay a33}-\ref{fig: IERK2-3 decay c2_tau_kappa}. We see that, for this example, the discrete energy generated by the IERK2-2 scheme is the closest to the continuous energy; while the discrete energy generated by the IERK2-Radau method is the farthest from the continuous one.
Experimentally, they confirm our theoretical predictions in previous subsection.

\begin{table}[htb!]\centering
  \begin{threeparttable}
    \centering \small
    \renewcommand\arraystretch{1.3}
    \belowrulesep=0pt\aboverulesep=0pt
    \caption{Average dissipation rate of second-order IERK methods.}
    \label{table: comparison of 2nd-order methods}
    \vspace{2mm}
    \begin{tabular}{c|c|c|c}
      \toprule 
      Type & Method & Average dissipation rate & Best  choice \\
      \midrule 
      \multirow{3}*{Lobatto-type} & 3-stage IERK2-1 \eqref{Scheme: IERK2-two parameters} & $c_2+\tfrac{1}{2c_2}+(2c_2a_{33}-\tfrac1{2})\tau\overline{\lambda}_{\mathrm{ML}}$ & $a_{33}=\frac{1}{2 (4c_2-2 c_2^2-1)}$ \\[2pt]
      \cline{2-4}
      & 3-stage IERK2-2 \eqref{Scheme: IERK2-one parameter} & $\sqrt{2}+\sqrt{2}(a_{33}-\tfrac{\sqrt{2}}{4})\tau\overline{\lambda}_{\mathrm{ML}}$ & $a_{33}=\tfrac{1+\sqrt{2}}{4}$ \\[2pt]
      \cline{2-4} & 3-stage IERK2 \cite{CooperSayfy:1983} & NPD$^*$ \\[2pt]
      \midrule
      \multirow{3}*{Radau-type} & 3-stage IERK2-Radau \eqref{Scheme: IERK2-Radau-one parameter} & $\tfrac{1}{2c_2}+c_2+\brab{2c_2+1+\tfrac{1}{c_2-1}}\tau\overline{\lambda}_{\mathrm{ML}}$ & $c_2=1+\frac{\sqrt{2}}{2}$ \\[2pt]
      \cline{2-4} & 4-stage IERK2 \cite{ShinLeeLee:2017} & $\tfrac{3}{2}+\tfrac{1}{2}\tau\overline{\lambda}_{\mathrm{ML}}$ \\[2pt]
      \cline{2-4} & 3-stage IERK2 \cite{AscherRuuthSpiteri:1997} & NPD \\[2pt]
      \bottomrule
    \end{tabular}
    \footnotesize
    \tnote{* NPD means there exists a $z_0<0$ such that the associated differential matrix $D(z_0)$ is not positive (semi-)definite.}
  \end{threeparttable}
\end{table}

To end this section, we present the average dissipation rates and practical parameter choices in Table \ref{table: comparison of 2nd-order methods}, which summarizes the above numerical results, the theoretical results in Theorems \ref{thm: IERK2-two parameters}-\ref{thm: IERK2-Radau}, and some results of existing IERK algorithms in the literature.

\section{Discrete energy laws of IERK3 methods}
\label{section: third-order methods}
\setcounter{equation}{0}

\subsection{IERK3 methods with four stages}
Third-order methods require three implicit stages, $s_I=3$, at least. We should determine eleven independent coefficients in the following four-stage IERK methods
that satisfy the canopy node condition and the two order conditions for first-order accuracy,
\begin{equation*}
  \begin{array}{c|c}
    \mathbf{c} & A \\
    \hline\\[-10pt]   & \mathbf{b}^T
  \end{array}
  = \begin{array}{c|cccc}
    0 & 0 &  &    &  \\
    c_{2} & c_2-a_{22} & a_{22} &   \\
    c_3 & c_3-a_{32}-a_{33} & a_{32} & a_{33}    \\
    1 & 1-a_{42}-a_{43}-a_{44} & a_{42} & a_{43}  & a_{44} \\
    \hline\\[-10pt]   & 1-a_{42}-a_{43}-a_{44} & a_{42} & a_{43} & a_{44}
  \end{array},
  \begin{array}{c|c}
    \hat{\mathbf{c}} & \widehat{A} \\
    \hline\\[-10pt]   & \hat{\mathbf{b}}^T
  \end{array}
  =\begin{array}{c|cccc}
    0 & 0 &  &   &  \\
    c_{2} & c_2 & 0 & &  \\
    c_3 & c_3-\hat{a}_{32} &\hat{a}_{32} & 0&   \\
      1 & 1-\hat{a}_{42}-\hat{a}_{43}  &\hat{a}_{42} & \hat{a}_{43} & 0 \\
    \hline\\[-10pt]   & 1-\hat{a}_{42}-\hat{a}_{43} &\hat{a}_{42} & \hat{a}_{43}&0
  \end{array}\;.
\end{equation*}

From Table \ref{table: order condition}, eight order conditions for second- and third-order accuracy are required such that there are three independent coefficients.
Considering the stand-alone conditions for explicit part, $\hat{\mathbf{b}}^{T} \widehat{A}  \mathbf{c} =\tfrac1{6}$, $\hat{\mathbf{b}}^T\mathbf{c}^{.2}=\tfrac13$ and $\hat{\mathbf{b}}^T\mathbf{c}=\tfrac12$, 
we find that $\hat{a}_{42}=\frac{3 c_3-2}{6 c_2 (c_3-c_2)}$, $\hat{a}_{43}=\frac{2-3c_2}{6 c_3(c_3-c_2)}$ 
and $\hat{a}_{32}=\frac{c_3 (c_3-c_2)}{c_2 (2-3 c_2)}$
with the independent variables $c_2$ and $c_3$. Here we have excluded the case $c_2 (c_3-c_2)=0$ since the case $c_2=c_3$ can not arrive at our aim.  From the stand-alone conditions for implicit part, ${\mathbf{b}}^T\mathbf{c}^{.2}=\tfrac13$ and ${\mathbf{b}}^T\mathbf{c}=\tfrac12$, one has
${a}_{42} = \frac{6 {a}_{44} c_3-6 {a}_{44}-3 c_3+2}{6 c_2 (c_2 - c_3)}$
and ${a}_{43} = \frac{6 {a}_{44} c_2-6 {a}_{44}-3 c_2+2}{6 c_3(c_3 - c_2)}.$
Thus the coupling condition $\mathbf{b}^{T} \widehat{A} \mathbf{c} = \tfrac{1}{6}$ arrives at $\tfrac{{a}_{44} c_2}{3 c_2-2}=0$ and then $$a_{44} = 0,$$ that is, the procedure at the fourth stage is purely explicit. Then it follows that ${a}_{42} = \hat{a}_{42}$ and ${a}_{43} = \hat{a}_{43}$. It means that the coupling condition $\hat{\mathbf{b}}^{T} A \mathbf{c} = \tfrac{1}{6}$ is equivalent to the stand-alone condition for implicit part $\mathbf{b}^{T} A \mathbf{c} = \tfrac{1}{6}$. Therefore we have four independent coefficients. 

By choosing $c_2 = \frac{1}{3}$, $c_3 = \frac{2}{3}$ and $a_{33} =\frac{1}{3}$, we obtain the following $a_{22}$-parameterized IERK3 methods with the associated Butcher tableaux
\begin{align}\label{Scheme: IERK3-4stage-a22}
\begin{array}{c|c}
	\mathbf{c} & A \\
	\hline\\[-10pt]   & \mathbf{b}^T
\end{array}
= \begin{array}{c|cccc}
    0 & 0 &  &    &  \\
   \tfrac{1}{3} &\tfrac{1}{3}-a_{22} & a_{22} \\[2pt]
  \tfrac{2}{3} & \tfrac{1}{3} & 0 & \tfrac{1}{3} \\[2pt]
  1 & \tfrac{1}{4} & 0 & \tfrac{3}{4} & 0 \\[1pt]
  \hline\\[-10pt]   & \tfrac{1}{4} & 0 & \tfrac{3}{4} & 0
\end{array}\,,\quad
\begin{array}{c|c}
	\hat{\mathbf{c}} & \widehat{A} \\
	\hline\\[-10pt]   & \hat{\mathbf{b}}^T
\end{array}
=\begin{array}{c|cccc}
   0 & 0  \\
  \tfrac{1}{3} &  \tfrac{1}{3} & 0  \\[2pt]
  \tfrac{2}{3} & 0 & \tfrac{2}{3} & 0 \\[2pt]
  1 & \tfrac{1}{4} & 0 & \tfrac{3}{4} & 0 \\[1pt]
   \hline\\[-10pt]   & \tfrac{1}{4} & 0 & \tfrac{3}{4} & 0
\end{array}\;.
\end{align}

By the definition \eqref{Def: Differential Matrix DII}, one can get
\begin{align*}
 D_{\mathrm{E}}^{(3)}:=\begin{pmatrix}
 3 & 0 & 0 \\[3pt]
 \tfrac{3}{2} & \tfrac{3}{2} & 0 \\[3pt]
 \tfrac{1}{3} & \tfrac{4}{3} & \tfrac{4}{3} 
  \end{pmatrix},\quad  D_{\mathrm{EI}}^{(3)}:=\begin{pmatrix}
 3 a_{22}-\tfrac{1}{2} & 0 & 0 \\[3pt]
 -\tfrac{1}{2} & 0 & 0 \\[3pt]
 -a_{22} & 0 & -\tfrac{1}{2}
  \end{pmatrix}.
\end{align*}
Simple computing yields that the matrix $D_{\mathrm{E}}^{(3)}$ is positive definite.  Also, it is easy to find that $\mathrm{Det}\kbrab{\mathcal{S}(D_{\mathrm{EI},1}^{(3)};a_{22})}
=3 a_{22}-\frac{1}{2}$, $\mathrm{Det}\kbrab{\mathcal{S}(D_{\mathrm{EI},2}^{(3)};a_{22})}
= -\tfrac{1}{16}$ and $\mathrm{Det}\kbrab{\mathcal{S}(D_{\mathrm{EI}}^{(3)};a_{22})} = \tfrac{1}{32}$. That is,  the associated differentiation  matrix $D^{(3)}(z)=D_{\mathrm{E}}^{(3)}-D_{\mathrm{EI}}^{(3)}z$ can not be positive semi-definite for all $z < 0$, and  the sufficient condition of Theorem \ref{thm: energy stability} is not fulfilled. Although the third-order methods \eqref{Scheme: IERK3-4stage-a22} may maintain the discrete energy dissipation law  \eqref{problem: energy dissipation law} under certain time-step conditions,  they may not be stabilized to preserve the energy dissipation law unconditionally  no matter how large the stabilization parameter $\kappa$ we set in the current stabilization strategy \eqref{problem: stabilized version}.

\begin{figure}[htb!]
  \centering
  \subfigure[$\tau=0.01,\kappa=4$]{
    \includegraphics[width=2.05in]{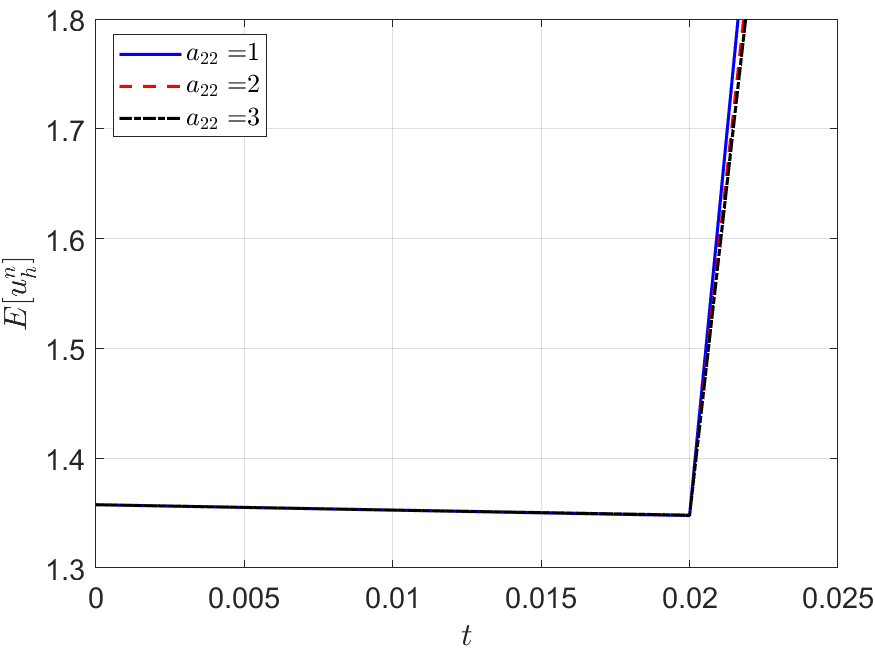}}  
  \subfigure[ $\tau=0.005,\kappa=4$]{
    \includegraphics[width=2.05in]{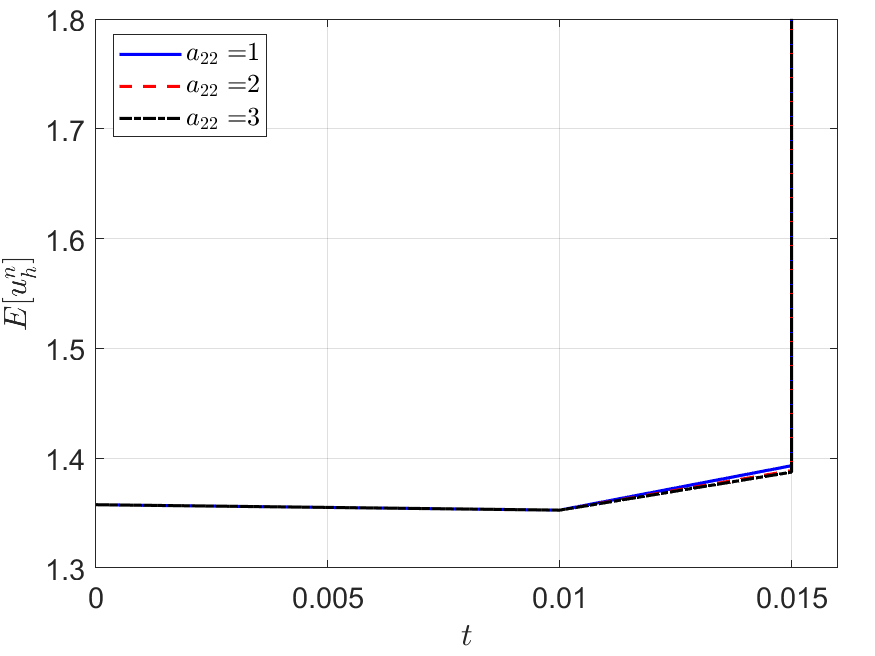}}
  \subfigure[$\tau=0.001,\kappa=4$]{
    \includegraphics[width=2.05in]{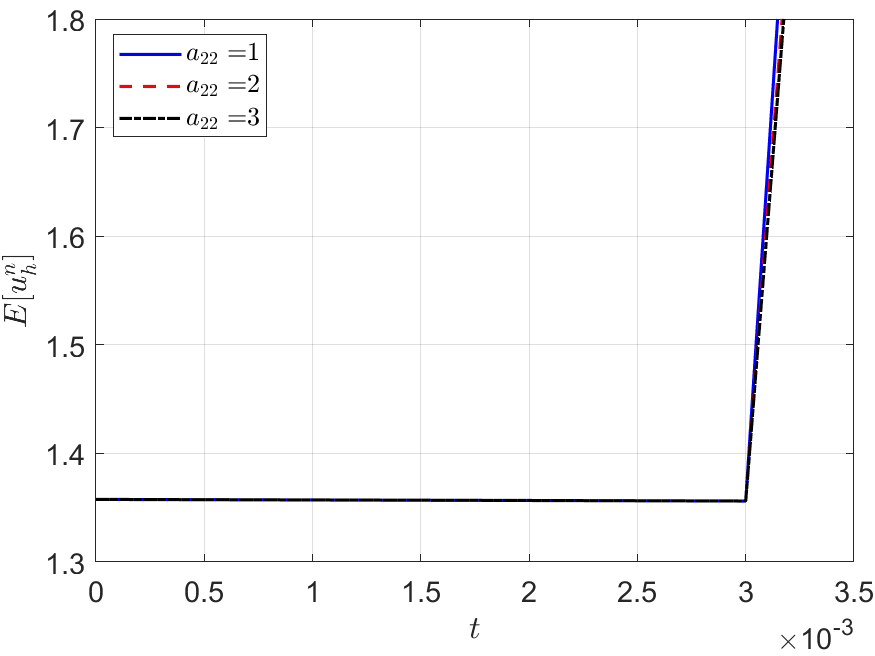}}
  \caption{Energy behaviors of IERK3 methods \eqref{Scheme: IERK3-4stage-a22}.}
  \label{fig: IERK3-1 decay a22_tau_kappa}
\end{figure}

As the numerical evidence,  we run the IERK3 methods \eqref{Scheme: IERK3-4stage-a22} with  $a_{22}=1$, 2 and 3 for Example \ref{example: energy_2} with a fixed stabilized factor $\kappa=4$. 
The associated discrete energy curves, depicted in Figure \ref{fig: IERK3-1 decay a22_tau_kappa}, for three different time-steps $\tau=0.01$, $0.005$ and $0.001$ can not maintain the decaying property. 
It is to mention that, the differentiation  matrix $D(z)=D_{\mathrm{E}}-D_{\mathrm{EI}}z$ can not be positive semi-definite for all $z < 0$ if there are one zero diagonal element in the coefficient matrix $A$ of the implicit part. For example, if $a_{k_0+1, k_0+1} = 0$, the formula \eqref{Def: Differential Matrix DII} gives the diagonal element $d_{k_0,k_0} = \tfrac{1}{\hat{a}_{k_0+1, k_0}} + \tfrac{z}{2}$, which cannot be always non-negative for all $z<0$. Then the condition of Theorem \ref{thm: energy stability} can not be fulfilled because
the diagonal elements of a positive semi-definite matrix must be non-negative. 
As mentioned before, four-stage IERK3 methods always have a zero diagonal element, $a_{44} \equiv 0$. In this sense, we say that there are no four-stage IERK3 methods fulfill the condition of Theorem \ref{thm: energy stability}. Note that, similar situations can be found in \cite[Section 2]{CooperSayfy:1983} and \cite[Sections 8 and 9]{LiuZou:2006}, where many IERK methods up to fourth-order accuracy have zero diagonal elements.

\subsection{IERK3 methods with five stages}
We consider the five-stage third-order IERK methods with the following Butcher tableaux 
\begin{equation*}
  \begin{array}{c|ccccc}
    0 & 0 &  &    &  &\\
    c_{2} & c_2-a_{55} & a_{55} &   &\\
    c_3 & c_3-a_{32}-a_{55} & a_{32} & a_{55} &   \\
    c_4 & c_4-a_{42}-a_{43}-a_{55} & a_{42} & a_{43} & a_{55}  \\
    1 & 1-a_{52}-a_{53}-a_{54}-a_{55} & a_{52} & a_{53}  &   a_{54} & a_{55}\\
    \hline\\[-10pt]   & 1-a_{52}-a_{53}-a_{54}-a_{55} & a_{52} & a_{53}  &   a_{54} & a_{55}
  \end{array},\quad
  \begin{array}{c|ccccc}
    0 & 0 &  &   & & \\
    c_{2} & c_2 & 0 & & & \\
    c_3 & c_3-c_2 &c_2 & 0&  & \\
    c_4 & c_4-\hat{a}_{42}-c_2 &\hat{a}_{42} & c_2&0  & \\
    1 & 1-\hat{a}_{52}-\hat{a}_{53}-c_2  &\hat{a}_{52} & \hat{a}_{53} & c_2&0\\
    \hline\\[-10pt]   & 1-\hat{a}_{52}-\hat{a}_{53}-c_2 &\hat{a}_{52} & \hat{a}_{53}&c_2&0
  \end{array},
\end{equation*}
where we assume $a_{k+1,k+1}=a_{55}$ and $\hat{a}_{k+1,k}=c_2$ for $1\le k\le 4$ to reduce the degree of freedom so that there are five independent variables. Such simplifying assumptions are set to maintain a fixed diagonal ratio $\frac{a_{k+1,k+1}}{\hat{a}_{k+1,k}}$ between the implicit and explicit parts so that one can minimize the value of $\frac1{s_{\mathrm{I}}}\mathrm{tr}(D_{\mathrm{EI}})$ as far as possible according to the formula \eqref{def: numerical rate Formula} of average dissipation rate.

After lots of trial and error processes, we keep one free variable $a_{55}$ by choosing the abscissas $c_2 = \tfrac45$, $c_3 = \tfrac75$, $c_4 = \tfrac65$ and the coefficient $a_{43} = -\tfrac35$,
and obtain the five-stage IERK3 (IERK3-1) methods with the following Butcher tableaux 
\begin{align}\label{Scheme: IERK3-5stage-a55}
\begin{array}{c|ccccc}
    0 & 0 &  &    &  \\[3pt]
  \tfrac{4}{5} &  \tfrac{4}{5}-a_{55} & a_{55} \\[3pt]
  \tfrac{7}{5} &  \tfrac{3}{5}-\tfrac{5 a_{55}}{16} & \tfrac{4}{5}-\tfrac{11 a_{55}}{16} & a_{55} \\[3pt]
  \tfrac{6}{5} &  \tfrac{977 a_{55}}{4032}-\tfrac{473}{10080} & \tfrac{18617}{10080}-\tfrac{5009 a_{55}}{4032} & -\tfrac{3}{5} & a_{55} \\[3pt]
  1 & a_{51} & a_{52} & a_{53} &  a_{54} & a_{55} \\[3pt]
  \hline\\[-10pt]   & a_{51} & a_{52} & a_{53} &  a_{54} & a_{55}
  \end{array}\,,\quad\begin{array}{c|ccccc}
  0 & 0 &  &   &  \\[3pt]
  \tfrac{4}{5} & \tfrac{4}{5} & 0  \\[3pt]
  \tfrac{7}{5} & \tfrac{3}{5} & \tfrac{4}{5} & 0  \\[3pt]
  \tfrac{6}{5} & \tfrac{10111}{10080} & -\tfrac{6079}{10080} & \tfrac{4}{5} & 0  \\[3pt]
  1 & \tfrac{313}{840} & \tfrac{131}{360} & -\tfrac{169}{315} & \tfrac{4}{5} & 0 \\[3pt]
  \hline\\[-10pt] & \tfrac{313}{840} & \tfrac{131}{360} & -\tfrac{169}{315} & \tfrac{4}{5} & 0
\end{array}\,,
\end{align}
where $a_{51}=\tfrac{313}{840}-\tfrac{191 a_{55}}{9590}$, $a_{52}=\tfrac{131}{360}-\tfrac{797 a_{55}}{4110}$, $a_{53}=\tfrac{7087 a_{55}}{14385}-\tfrac{169}{315}$ 
and $a_{54}=\tfrac{4}{5}-\tfrac{876 a_{55}}{685}$.

By the definitions in \eqref{Def: Differential Matrix D}, one can write out the matrices $D_{\mathrm{E}}^{(3,1)}$ and $D_{\mathrm{EI}}^{(3,1)}$ as follows,
\begin{small}\begin{align*}
 D_{\mathrm{E}}^{(3,1)}:=\begin{pmatrix}
   \tfrac{5}{6} & 0 & 0 & 0 \\[3pt]
   \tfrac{49}{30} & \tfrac{48}{35} & 0 & 0 \\[3pt]
   \tfrac{23}{9} & \tfrac{62}{21} & \tfrac{5}{2} & 0 \\[3pt]
   \tfrac{583}{1134} & \tfrac{851}{1323} & \tfrac{125}{126} & \tfrac{5}{3} 
  \end{pmatrix},  D_{\mathrm{EI}}^{(3,1)}:=\begin{pmatrix}
 \tfrac{5 a_{55}-2}{4} & 0 & 0 & 0 \\[3pt]
 -\tfrac{35 a_{55}}{64} & \tfrac{5 a_{55}-2}{4} & 0 & 0 \\[3pt]
 \tfrac{21}{16}-\tfrac{130885 a_{55}}{57344} & \tfrac{70715 a_{55}}{32256}-\tfrac{7}{4} & \tfrac{5 a_{55}-2}{4} & 0 \\[3pt]
 \tfrac{169}{192}-\tfrac{1213756279 a_{55}}{659914752} & \tfrac{1301752603 a_{55}}{1113606144}-\tfrac{169}{144} & \tfrac{67633 a_{55}}{138096} & \tfrac{5 a_{55}-2}{4}
  \end{pmatrix}.
\end{align*}\end{small}
Simple calculations can verify the positive semi-definiteness of $D_{\mathrm{E}}^{(3,1)}$.
For the matrix $D_{\mathrm{EI}}^{(3,1)}(a_{55})$,
We have $\mathrm{Det}\kbrab{\mathcal{S}(D_{\mathrm{EI},1}^{(3,1)};a_{55})}
=\tfrac{5 a_{55}-2}{4}$, $\mathrm{Det}\kbrab{\mathcal{S}(D_{\mathrm{EI},2}^{(3,1)};a_{55})}
= \tfrac{24375 a_{55}^2}{16384}-\tfrac{5 a_{55}}{4}+\tfrac{1}{4}$, and
\begin{align*}
  \mathrm{Det}\kbrab{\mathcal{S}(D_{\mathrm{EI},3}^{(3,1)};a_{55})}
  =&\,-\tfrac{156124878125 a_{55}^3}{266355081216}+\tfrac{4850385555625 a_{55}^2}{2130840649728}-\tfrac{9210215 a_{55}}{4718592}+\tfrac{969}{2048},\\
  \mathrm{Det}\kbrab{\mathcal{S}(D_{\mathrm{EI},4}^{(3,1)};a_{55})}
  =&\,-\tfrac{7024871482286543654375 a_{55}^4}{5079525965637831622656}+\tfrac{11933410165242198623465 a_{55}^3}{2539762982818915811328} \\
  &\, -\tfrac{26459898695651552391913 a_{55}^2}{5079525965637831622656}+\tfrac{192673313809999 a_{55}}{82103953784832}-\tfrac{1969849}{5308416}.
\end{align*}
One can check that the following condition 
\begin{align}\label{cond: IERK3-5stage a55 condition}
  0.717374\leq a_{55}\leq 1.74727
\end{align}
is sufficient for the positive semi-definiteness of $D_{\mathrm{EI}}^{(3,1)}$.
By using Corollary \ref{corol: energy stability} and Lemma \ref{lemma: average dissipation rate}, we have the following theorem.
\begin{theorem}\label{thm: IERK3-5stage-a55s} 
  In simulating the gradient flow system \eqref{problem: stabilized version} with  $\kappa\ge2\ell_g$, the one-parameter IERK3-1 methods \eqref{Scheme: IERK3-5stage-a55} with the parameter setting  \eqref{cond: IERK3-5stage a55 condition} preserve
  the original energy dissipation law \eqref{problem: energy dissipation law} unconditionally
  at all stages in the sense that
  \begin{align*}    
    E[U^{n,j+1}]-E[U^{n,1}]\le\frac1{\tau}\myinnerB{M_h^{-1}\delta_{\tau}\vec{U}_{n,j+1},
      D_{j}^{(3,1)}(a_{55};{\tau}M_hL_{\kappa})\delta_{\tau}\vec{U}_{n,j+1}}
    \quad\text{for $1\le j\le 4$,}
  \end{align*}
  where the  differentiation  matrix $D^{(3,1)}$ is defined by  $$D^{(3,1)}(a_{55};{\tau}M_hL_{\kappa}):=D_{\mathrm{E}}^{(3,1)}\otimes I+D_{\mathrm{EI}}^{(3,1)}(a_{55}) \otimes (-\tau M_h L_{\kappa}).$$
  The associated average dissipation rate reads
  \begin{align}\label{AverRate: IERK3-5stage-a55}
    \mathcal{R}^{(3,1)}(a_{55}):=\tfrac{5}{4}+\tfrac{5a_{55}-2}{4}\tau\overline{\lambda}_{\mathrm{ML}}
    \quad\text{for $0.717374\leq a_{55}\leq 1.74727$}.
  \end{align}
\end{theorem}

Also, one can take $a_{43}$ as the free variable by choosing the abscissas $c_2 = \tfrac45$, $c_3 = \tfrac75$, $c_4 = \tfrac65$ and the diagonal element $a_{55} =\tfrac{18}{25}$,
and obtain the following one-parameter IERK3 (IERK3-2) methods with the associated Butcher tableaux 
\begin{small}\begin{equation}\label{Scheme: IERK3-5stage-a43}
  \begin{array}{c|ccccc}
    0 & 0  \\[3pt]
    \tfrac{4}{5} & \tfrac{2}{25} & \tfrac{18}{25} \\[3pt]
   \tfrac{7}{5} & \tfrac{3}{8} & \tfrac{61}{200} & \tfrac{18}{25} \\[3pt]
   \tfrac{6}{5} & \tfrac{3 a_{43}}{4}+\tfrac{7277}{12600} & -\tfrac{7 a_{43}}{4}-\tfrac{1229}{12600} & a_{43} & \tfrac{18}{25} \\[3pt]
  1 & \tfrac{1030769}{2877000} & \tfrac{276523}{1233000} & -\tfrac{196127}{1078875} & -\tfrac{2068}{17125} & \tfrac{18}{25} \\[3pt]
    \hline\\[-10pt]   & \tfrac{1030769}{2877000} & \tfrac{276523}{1233000} & -\tfrac{196127}{1078875} & -\tfrac{2068}{17125} & \tfrac{18}{25}
\end{array},\;\;\;\begin{array}{c|ccccc}
  0 & 0 &  &   &  \\[3pt]
  \tfrac{4}{5} & \tfrac{4}{5} & 0  \\[3pt]
  \tfrac{7}{5} & \tfrac{3}{5} & \tfrac{4}{5} & 0  \\[3pt]
  \tfrac{6}{5} & \tfrac{10111}{10080} & -\tfrac{6079}{10080} & \tfrac{4}{5} & 0  \\[3pt]
  1 & \tfrac{313}{840} & \tfrac{131}{360} & -\tfrac{169}{315} & \tfrac{4}{5} & 0 \\[3pt]
  \hline\\[-10pt] & \tfrac{313}{840} & \tfrac{131}{360} & -\tfrac{169}{315} & \tfrac{4}{5} & 0
  \end{array}\,.
\end{equation}
\end{small}
Note that, the IERK3-2 methods \eqref{Scheme: IERK3-5stage-a43} is quite different from the IERK3-1 methods \eqref{Scheme: IERK3-5stage-a55} in the sense that the average dissipation rate $\mathcal{R}^{(3,2)}$ is independent of the parameter $a_{43}$, while the average dissipation rate $\mathcal{R}^{(3,1)}(a_{55})$ of the latter is always dependent on the choice of $a_{55}$.

By the definitions in \eqref{Def: Differential Matrix D}, one has  $D_{\mathrm{E}}^{(3,2)}:=D_{\mathrm{E}}^{(3,1)}$ and
\begin{align*}
  D_{\mathrm{EI}}^{(3,2)}(a_{43}):=\begin{pmatrix}
   \tfrac{2}{5} & 0 & 0 & 0 \\[3pt]
   -\tfrac{63}{160} & \tfrac{2}{5} & 0 & 0 \\[3pt]
   -\tfrac{15 a_{43}}{16}-\tfrac{128073}{143360} & \tfrac{5 a_{43}}{4}+\tfrac{5183}{8960} & \tfrac{2}{5} & 0 \\[3pt]
   -\tfrac{845 a_{43}}{1344}-\tfrac{6774788911}{8248934400} & \tfrac{845 a_{43}}{1008}+\tfrac{264498203}{1546675200} & \tfrac{67633}{191800} & \tfrac{2}{5}
  \end{pmatrix}\,.
\end{align*}
Simple calculations gives $\mathrm{Det}\kbrab{\mathcal{S}(D_{\mathrm{EI},1}^{(3,2)};a_{43})} = \tfrac{2}{5}$, $\mathrm{Det}\kbrab{\mathcal{S}(D_{\mathrm{EI},2}^{(3,2)};a_{43})} = \tfrac{2483}{20480}$, 
\begin{align*}
  \mathrm{Det}\kbrab{\mathcal{S}(D_{\mathrm{EI},3}^{(3,2)};a_{43})}
  =&-\tfrac{1055 a_{43}^2}{8192}-\tfrac{2184569 a_{43}}{14680064}-\tfrac{91452759}{6576668672},\\
  \mathrm{Det}\kbrab{\mathcal{S}(D_{\mathrm{EI}}^{(3,2)};a_{43})}
  =&-\tfrac{14878761752461399 a_{43}^2}{499921851934310400}-\tfrac{523061384625360827 a_{43}}{17497264817700864000}-\tfrac{342702875992479487397}{48992341489562419200000}.
\end{align*}
It is easy to check that the following condition 
\begin{align}\label{cond: IERK3-5stage-a43 condition}
	-0.633312 \leq a_{43} \leq -0.371114
\end{align}
is sufficient for the positive semi-definiteness of $D_{\mathrm{EI}}^{(3,2)}$. 
Thus we have the following result.

\begin{theorem}\label{thm: IERK3-5stage-a43} 
	In simulating the gradient flow system \eqref{problem: stabilized version} with  $\kappa\ge2\ell_g$, the one-parameter IERK3-2 methods \eqref{Scheme: IERK3-5stage-a43} with the parameter setting  \eqref{cond: IERK3-5stage-a43 condition} preserve
	the original energy dissipation law \eqref{problem: energy dissipation law} unconditionally
	at all stages in the sense that
	\begin{align*}    
		E[U^{n,j+1}]-E[U^{n,1}]\le\frac1{\tau}\myinnerB{M_h^{-1}\delta_{\tau}\vec{U}_{n,j+1},
			D_{j}^{(3,2)}(a_{43};{\tau}M_hL_{\kappa})\delta_{\tau}\vec{U}_{n,j+1}}
		\quad\text{for $1\le j\le 4$,}
	\end{align*}
	where the  differentiation  matrix $D^{(3,2)}$ is defined by  $$D^{(3,2)}(a_{43};{\tau}M_hL_{\kappa}):=D_{\mathrm{E}}^{(3,2)}\otimes I+D_{\mathrm{EI}}^{(3,2)}(a_{43}) \otimes (-\tau M_h L_{\kappa}).$$
	The associated average dissipation rate is parameter-independent, that is,
	\begin{align}\label{AverRate: IERK3-5stage-a43}
		\mathcal{R}^{(3,2)}
		:=\tfrac{5}{4} + \tfrac{2}{5}\tau\overline{\lambda}_{\mathrm{ML}}\quad\text{for $-0.633312 \leq a_{43} \leq -0.371114$}.
	\end{align}
\end{theorem}

\subsection{Radau-type IERK3 methods with five stages}
We also briefly discuss the five-stage Radau-type IERK3 methods by assuming $\mathbf{a}_1=\mathbf{0}$. That is, consider the following Butcher tableaux that satisfy the canopy node condition and the two order conditions for first-order accuracy,
\begin{equation*}
  \begin{array}{c|ccccc}
    0 & 0 &  &    &  &\\
    c_{2} & 0 & c_{2} &   &\\
    c_3 & 0 & c_3-c_2 & c_{2} &   \\
    c_4 & 0 & c_4-a_{43}-c_{2} & a_{43} & c_{2}  \\
    1 & 0 & 1-a_{53}-a_{54}-c_{2} & a_{53}  &   a_{54} & c_{2}\\
    \hline\\[-10pt]   & 0 & 1-a_{53}-a_{54}-c_{2} & a_{53}  &   a_{54} & c_{2}
  \end{array},\quad
  \begin{array}{c|ccccc}
    0 & 0 &  &   & & \\
    c_{2} & c_2 & 0 & & & \\
    c_3 & c_3-\hat{a}_{32} &\hat{a}_{32} & 0&  & \\
    c_4 & c_4-\hat{a}_{42}-\hat{a}_{43} &\hat{a}_{42} & \hat{a}_{43}&0  & \\
    1 & 1-\hat{a}_{52}-\hat{a}_{53}-\hat{a}_{54}  &\hat{a}_{52} & \hat{a}_{53} & \hat{a}_{54}&0\\
    \hline\\[-10pt]  & 1-\hat{a}_{52}-\hat{a}_{53}-\hat{a}_{54}  &\hat{a}_{52} & \hat{a}_{53} & \hat{a}_{54}&0
  \end{array},
\end{equation*}
where we set $a_{k+1,k+1}=c_2$ ($1\le k\le 4$) to reduce the degree of freedom.  According to Table \ref{table: order condition}, eight order conditions for second- and third-order accuracy are required so that there are four independent coefficients. This simplifying setting is also useful for practical applications since it provides the same iteration matrix for the systems of linear equations at all stages.

After lots of trial and error processes, we choose  $c_{2}=\tfrac{4}{5}$, $c_{3}=\tfrac{93}{200}$ and $c_4=\frac{171}{200}$ and obtain the $\hat{a}_{43}$-parameterized Radau-type IERK3 (IERK3-Radau) methods with the Butcher tableaux
\begin{align}\label{Scheme: IERK3-5stage-Radau}
  &\begin{array}{c|ccccc}
     0 & 0 &  \\[3pt]
    \tfrac{4}{5} & 0 & \tfrac{4}{5} \\[3pt]
    \tfrac{93}{200} & 0 & -\tfrac{67}{200} & \tfrac{4}{5} \\[3pt]
    \tfrac{171}{200} & 0 & -\tfrac{9361649}{5132200} & \tfrac{241098}{128305} & \tfrac{4}{5} \\[3pt]
    1 & 0 & -\tfrac{5309}{11055} & \tfrac{9998}{7839} & -\tfrac{766}{1287} & \tfrac{4}{5} \\[3pt]
    \hline\\[-10pt]   & 0 & \tfrac{25}{1488} & \tfrac{4024}{1395} & -\tfrac{1397}{144} & \tfrac{39}{5}
  \end{array},\quad
  \begin{array}{c|ccccc}
    0 & 0 &  \\[3pt]
  \tfrac{4}{5} &  \tfrac{4}{5} & 0 \\[3pt]
  \tfrac{93}{200} &  \tfrac{10391}{32000} & \tfrac{4489}{32000} & 0 \\[3pt]
  \tfrac{171}{200} & \hat{a}_{41}  & \hat{a}_{42} & \hat{a}_{43} & 0 \\[3pt]
  1 & \tfrac{2053}{11066} & \tfrac{3785983}{24466926} & \tfrac{20893310}{43373187} & \tfrac{1267730}{7120971} & 0 \\[3pt]
    \hline\\[-10pt]   & \tfrac{2053}{11066} & \tfrac{3785983}{24466926} & \tfrac{20893310}{43373187} & \tfrac{1267730}{7120971} & 0
  \end{array},
\end{align}
where $\hat{a}_{42} = \tfrac{9690263}{12256000}-\tfrac{93 \hat{a}_{43}}{160}$ and $\hat{a}_{41}=\tfrac{171}{200}-\hat{a}_{42}-\hat{a}_{43}$.

By the definition \eqref{Def: Differential Matrix D}, one can get the associated matrices
\begin{align*}
[D_{\mathrm{rad},\mathrm{E}}^{(3)}(\hat{a}_{43})]^T:=\begin{pmatrix}
 \tfrac{5}{4} & \tfrac{1135}{268} & \tfrac{200}{67}-\tfrac{4986267}{2052880 \hat{a}_{43}} & \tfrac{114597584391147}{17436733668080 \hat{a}_{43}}-\tfrac{2528991857}{339751640} \\[3pt]
 0 & \tfrac{32000}{4489} & \tfrac{18600}{4489}-\tfrac{7970976}{1719287 \hat{a}_{43}} & \tfrac{183194079827616}{14603264447017 \hat{a}_{43}}-\tfrac{67097015181}{5690839970} \\[3pt]
 0 & 0 & \tfrac{1}{\hat{a}_{43}} & \tfrac{7120971}{1267730}-\tfrac{22982641}{8493791 \hat{a}_{43}} \\[3pt]
 0 & 0 & 0 & \tfrac{7120971}{1267730}
\end{pmatrix},
\end{align*}
and
\begin{align*}
D_{\mathrm{rad},\mathrm{EI}}^{(3)}(\hat{a}_{43}):=\begin{pmatrix}       
 \tfrac{1}{2} & 0 & 0 & 0 \\[3pt]
 0 & \tfrac{46711}{8978} & 0 & 0 \\[3pt]
 0 & \tfrac{10391}{4489}-\tfrac{15730338}{8596435 \hat{a}_{43}} & \tfrac{4}{5 \hat{a}_{43}}-\tfrac{1}{2} & 0 \\[3pt]
 0 & \tfrac{361524711062658}{73016322235085 \hat{a}_{43}}-\tfrac{94060069247}{14227099925} & \tfrac{476922}{3169325}-\tfrac{91930564}{42468955 \hat{a}_{43}} & \tfrac{25314559}{6338650}
  \end{pmatrix}.
\end{align*}
As done in the above subsection, one can check that, if $0.598442 \leq \hat{a}_{43} \leq 1.05134$, $D_{\mathrm{rad},\mathrm{E}}^{(3)}(\hat{a}_{43})$ and $D_{\mathrm{rad},\mathrm{EI}}^{(3)}(\hat{a}_{43})$ are positive semi-definite.
Corollary \ref{corol: energy stability} and Lemma \ref{lemma: average dissipation rate} yield the following theorem.
\begin{theorem}\label{thm: IERK3-5stage-Radau} 
  In simulating the gradient flow system \eqref{problem: stabilized version} with  $\kappa\ge2\ell_g$, the one-parameter IERK3-Radau methods \eqref{Scheme: IERK3-5stage-Radau} for $0.598442 \leq \hat{a}_{43} \leq 1.05134$ preserve the original energy dissipation law \eqref{problem: energy dissipation law} unconditionally
  at all stages in the sense that
  \begin{align*}    
    E[U^{n,j+1}]-E[U^{n,1}]\le\frac1{\tau}\myinnerB{M_h^{-1}\delta_{\tau}\vec{U}_{n,j+1},
      D_{\mathrm{rad},j}^{(3)}(\hat{a}_{43};{\tau}M_hL_{\kappa})\delta_{\tau}\vec{U}_{n,j+1}}
    \quad\text{for $1\le j\le 3$,}
  \end{align*}
  where the  differentiation  matrix $D_{\mathrm{rad}}^{(3)}$ is defined by  $$D_{\mathrm{rad}}^{(3)}(\hat{a}_{43};{\tau}M_hL_{\kappa}):=D_{{\mathrm{rad}},\mathrm{E}}^{(3)}(\hat{a}_{43})\otimes I+D_{{\mathrm{rad}},\mathrm{EI}}^{(3)}(\hat{a}_{43}) \otimes (-\tau M_h L_{\kappa}).$$
  The associated average dissipation rate is 
  \begin{small}\begin{align}\label{AverRate: IERK3-5stage-Radau}
    \mathcal{R}_{\mathrm{rad}}^{(3)}(\hat{a}_{43})
    :=\tfrac{1}{4}\brab{\tfrac{1}{\hat{a}_{43}}+\tfrac{159293897563}{11381679940}}+\brab{\tfrac{1}{5 \hat{a}_{43}}+\tfrac{130839697713}{56908399700}}\tau\overline{\lambda}_{\mathrm{ML}}\quad\text{for $0.598442 \leq \hat{a}_{43} \leq 1.05134$.}
  \end{align}
  \end{small}
\end{theorem}

To enlarge the possible choices of $\tau\overline{\lambda}_{\mathrm{ML}}$, one can choose the parameter $\hat{a}_{43}=1$ such that $$ \mathcal{R}_{\mathrm{rad}}^{(3)}(1)=\tfrac{170675577503}{45526719760}+\tfrac{142221377653}{56908399700}\tau\overline{\lambda}_{\mathrm{ML}}\approx 3.74891+2.49913\tau\overline{\lambda}_{\mathrm{ML}}.$$
It is a little larger than the average dissipation rate $\mathcal{R}^{(3,2)}=\tfrac{5}{4} + \tfrac{2}{5}\tau\overline{\lambda}_{\mathrm{ML}}$ of the IERK3-2 methods \eqref{Scheme: IERK3-5stage-a43}.
From the perspective of average dissipation rate, the IERK3-Radau methods \eqref{Scheme: IERK3-5stage-Radau} are more competitive than the existing Radau-type IERK3 method in \cite[Example 5]{FuTangYang:2024} with the associated average dissipation rate $\mathcal{R}_{\mathrm{rad},F}^{(3)}\approx3.24727+16.3779 \tau\overline{\lambda}_{\mathrm{ML}}$, 
although they would be inferior to the one-parameter IERK3-1 \eqref{Scheme: IERK3-5stage-a55} and IERK3-2 \eqref{Scheme: IERK3-5stage-a43} methods. 

 To search a ``better'' IERK3 method, that is to say, to obtain a IERK method with the average dissipation rate closer to 1, one can consider higher-stage procedures. An example is the 7-stage IERK3 method in \cite{ShinLeeLee:2017}, which achieves a better average dissipation rate, $\mathcal{R}_{\mathrm{rad},S}^{(3)} = 2 + \tfrac{1}{2}\tau\overline{\lambda}_{\mathrm{ML}}$, at the cost of two additional implicit stages. 
 We also mentioned that, the 5-stage  Radau-type IERK method in \cite{AscherRuuthSpiteri:1997} does not satisfy the sufficient conditions of Corollary \ref{corol: energy stability} so that it is not a proper candidate for the gradient flow system \eqref{problem: gradient flows}. 

\subsection{Tests for IERK3 methods}

\begin{figure}[htb!]
	\centering
	\subfigure[IERK3-1 \eqref{Scheme: IERK3-5stage-a55}]{
		\includegraphics[width=2in]{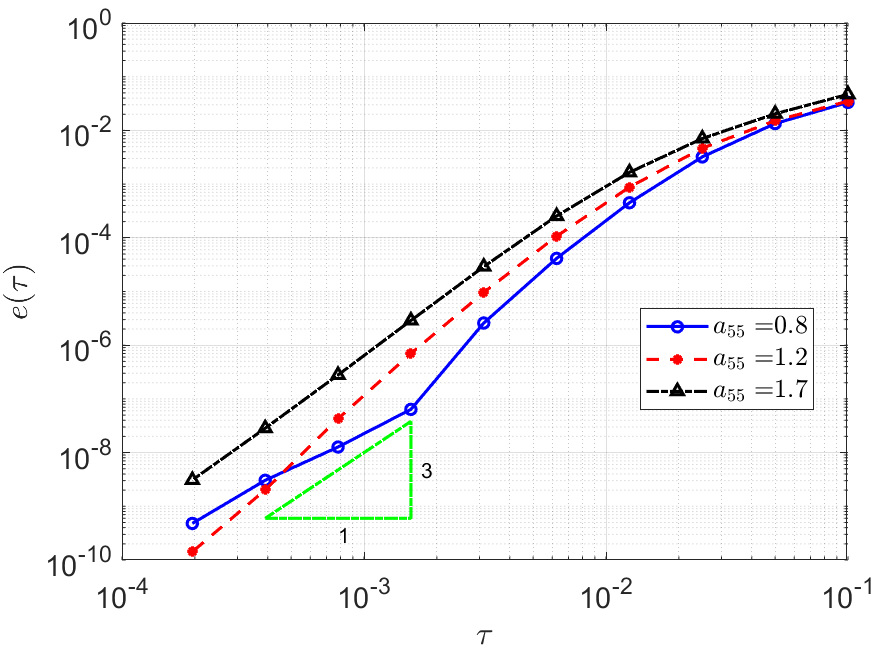}}
	\subfigure[IERK3-2 \eqref{Scheme: IERK3-5stage-a43}]{
		\includegraphics[width=2in]{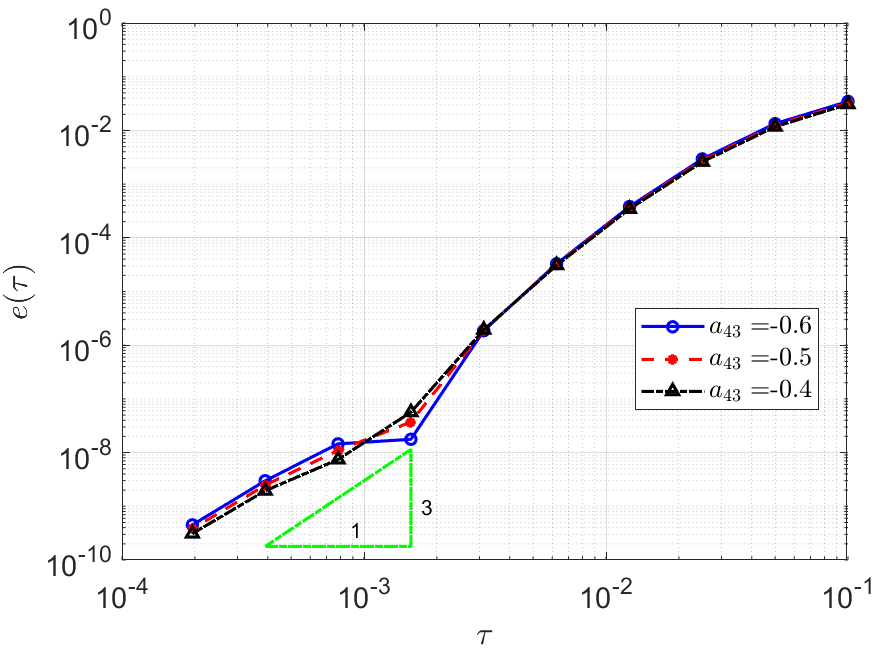}}
	\subfigure[IERK3-Radau \eqref{Scheme: IERK3-5stage-Radau}]{
		\includegraphics[width=2in]{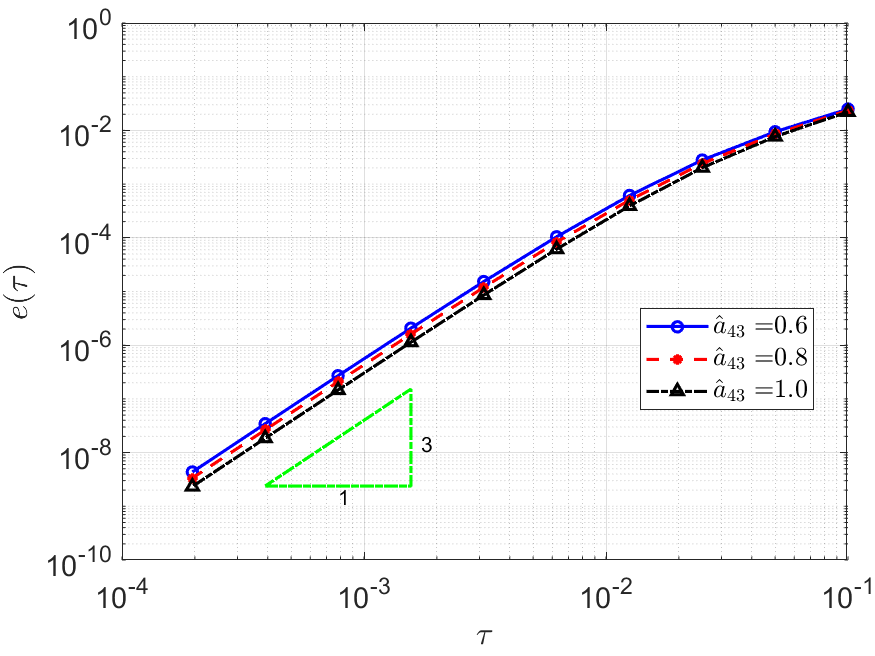}}
	\caption{Solution errors of IERK3 methods with different method parameters.}
	\label{fig: IERK3 convergence}
\end{figure}

We use Example \ref{ex: with source term} to test the convergence of IERK3-1 \eqref{Scheme: IERK3-5stage-a55}, IERK3-2 \eqref{Scheme: IERK3-5stage-a43} and IERK3-Radau \eqref{Scheme: IERK3-5stage-Radau} methods by choosing the final time $T=1$ and the stabilized parameter $\kappa=4$. Figure \ref{fig: IERK3 convergence} lists the $L^\infty$ norm error $e(\tau):=\max_{1\le n\le N}\mynorm{u_h^{n} - u(t_n)}_\infty$ for the three classes of IERK3 methods on halving time steps $\tau=2^{-k}/10$ for $0\le k\le9$. As expected, the IERK3-1 \eqref{Scheme: IERK3-5stage-a55}, IERK3-2 \eqref{Scheme: IERK3-5stage-a43} and IERK3-Radau \eqref{Scheme: IERK3-5stage-Radau} methods are third-order accurate in time. It is observed that the different parameters for the IERK3-1 \eqref{Scheme: IERK3-5stage-a55} and IERK3-2 \eqref{Scheme: IERK3-5stage-a43} methods would arrive at different numerical precision; while the IERK3-Radau  methods \eqref{Scheme: IERK3-5stage-Radau} with three different parameters $\hat{a}_{43}=0.6, 0.8$ and $1$ generates almost the same solution.

\begin{figure}[htb!]
	\centering
	\subfigure[$\tau=0.01,\kappa=2$]{
		\includegraphics[width=2.05in]{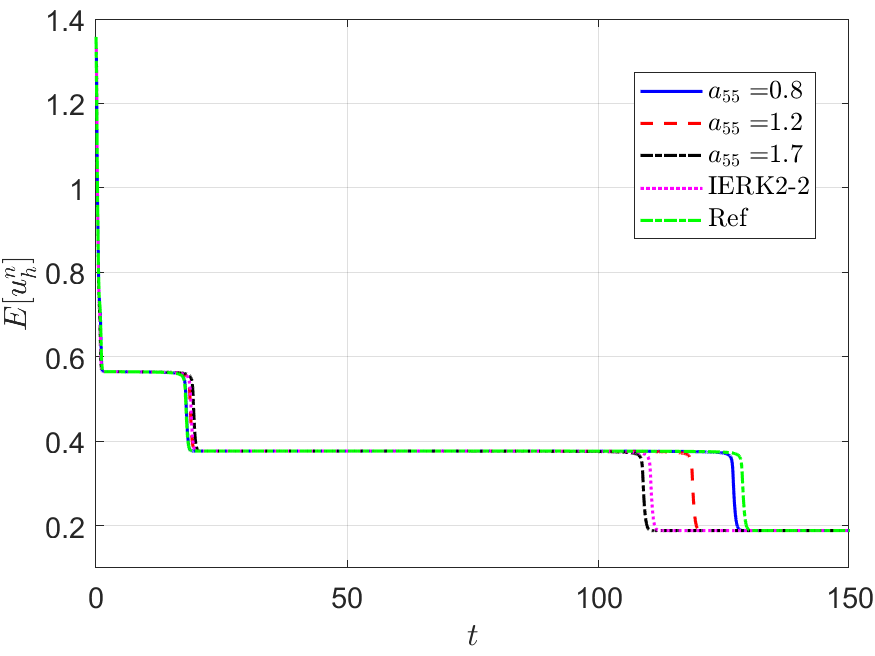}}
	\subfigure[ $\tau=0.05,\kappa=2$]{
		\includegraphics[width=2.05in]{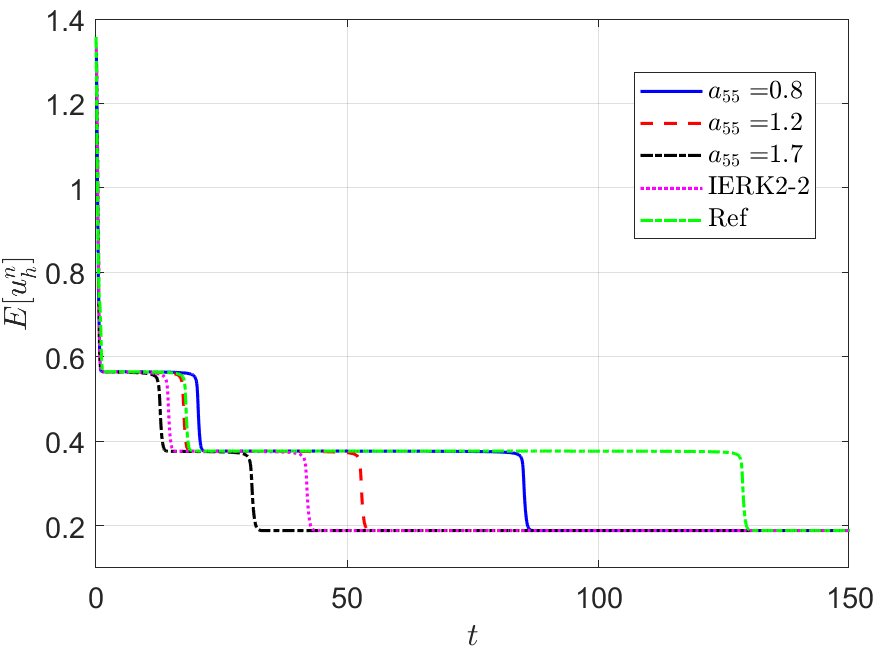}}
	\subfigure[$\tau=0.05,\kappa=3$]{
		\includegraphics[width=2.05in]{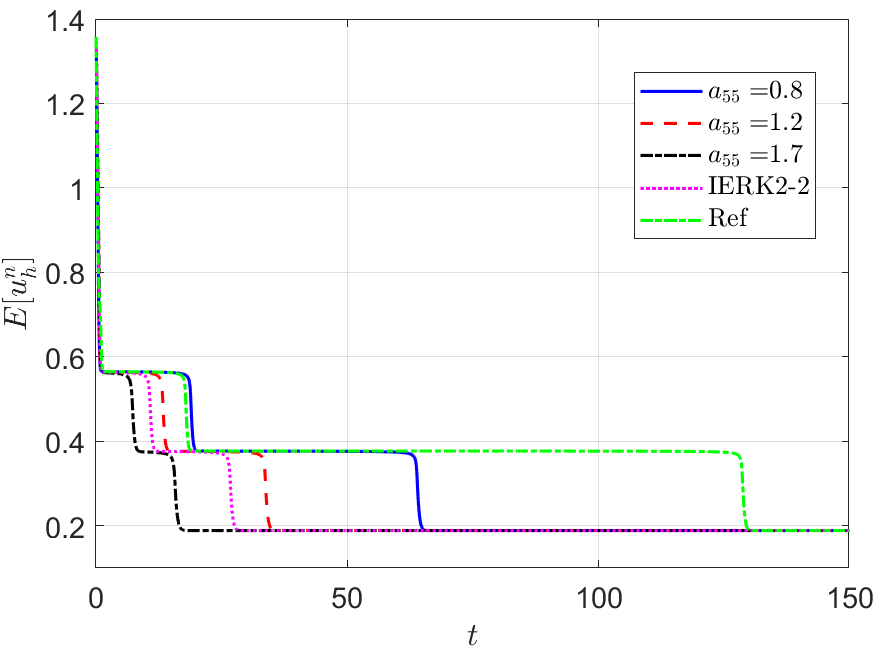}}
	\caption{Energy behaviors of IERK3-1 methods \eqref{Scheme: IERK3-5stage-a55}.}
	\label{fig: IERK3-2 decay a55_tau_kappa}
\end{figure}

\begin{figure}[htb!]
	\centering
	\subfigure[$\tau=0.01,\kappa=2$]{
		\includegraphics[width=2.05in]{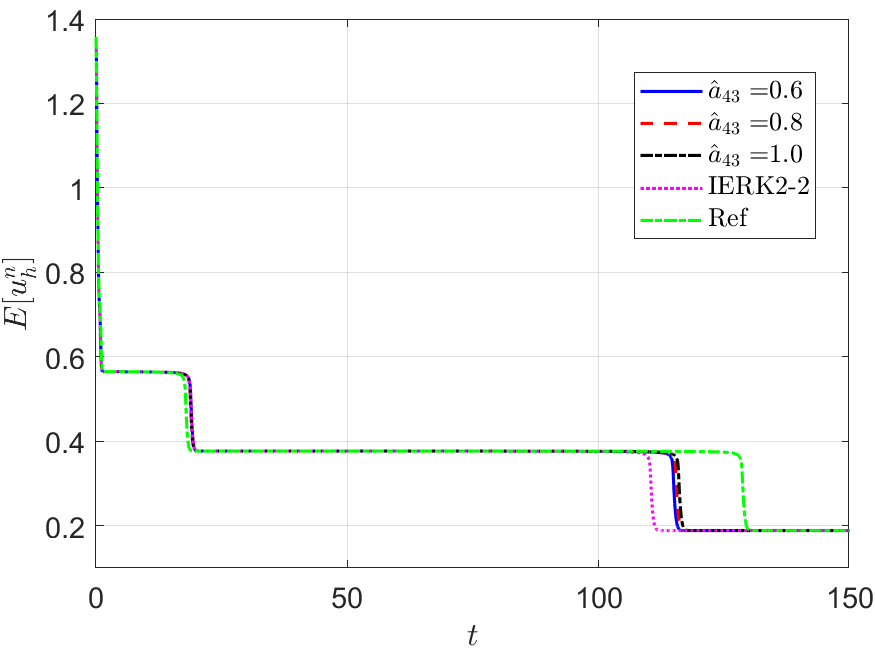}}
	\subfigure[ $\tau=0.05,\kappa=2$]{
		\includegraphics[width=2.05in]{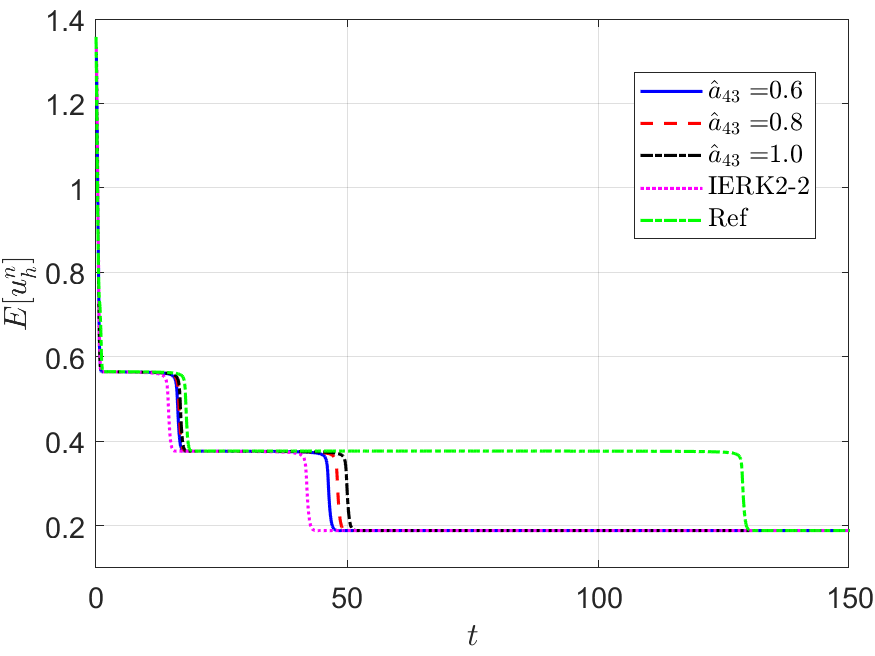}}
	\subfigure[$\tau=0.05,\kappa=3$]{
		\includegraphics[width=2.05in]{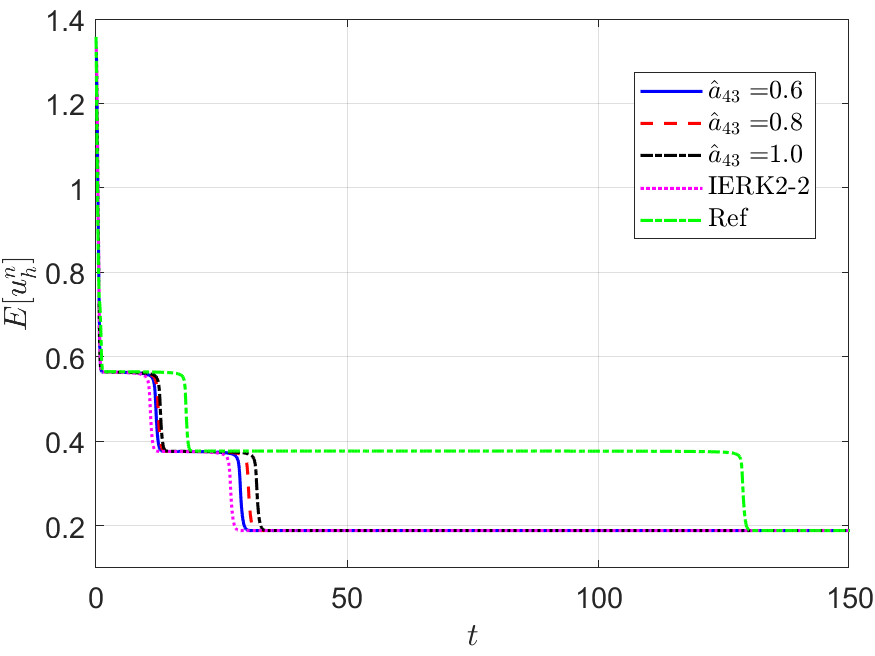}}
	\caption{Energy behaviors of IERK3-Radau methods \eqref{Scheme: IERK3-5stage-Radau}.}
	\label{fig: IERK3-4 decay ah43_tau_kappa}
\end{figure}

Next, we address the discrete energy behaviors of three IERK3 methods for Example \ref{example: energy_2}. Always, we use a small time-step size $\tau=10^{-3}$ to compute the reference solution and energy by the IERK3-2 method \eqref{Scheme: IERK3-5stage-a43} for $a_{43}=-0.4$. Figures \ref{fig: IERK3-2 decay a55_tau_kappa} and \ref{fig: IERK3-4 decay ah43_tau_kappa} depict the discrete energies generated by the IERK3-1 \eqref{Scheme: IERK3-5stage-a55} and IERK3-Radau \eqref{Scheme: IERK3-5stage-Radau} methods, respectively, for three different scenes: (i) $\tau=0.01,\kappa=2$; (ii) $\tau=0.05,\kappa=2$ and (iii) $\tau=0.05,\kappa=3$. As predicted by our theory, the discrete energies in all cases always monotonically decrease, and the discrete energies generated by the same method gradually deviate from the continuous energy as the associated average dissipation rates deviate from 1, cf. the formulas \eqref{AverRate: IERK3-5stage-a55} and \eqref{AverRate: IERK3-5stage-Radau}. It is seen that there are always significant differences between the reference curve and the three discrete energy curves (corresponds to three different parameters) generated by  the IERK3-1 methods \eqref{Scheme: IERK3-5stage-a55}; but such differences are not obvious for those of the IERK3-Radau  methods \eqref{Scheme: IERK3-5stage-Radau}. 

The above differences between discrete energy curves and continuous energy shown in Figures \ref{fig: IERK3-2 decay a55_tau_kappa} and \ref{fig: IERK3-4 decay ah43_tau_kappa} can be partly explained by the different precision of numerical solutions. Actually,  Figure \ref{fig: IERK3 convergence} (a) shows that the solution for the case $a_{55}=0.8$ of the IERK3-1 methods \eqref{Scheme: IERK3-5stage-a55} is a bit more accurate than that for $a_{55}=1.7$ although both of them are third-order accurate. Figure \ref{fig: IERK3 convergence} (c) shows that the solutions for the three cases $\hat{a}_{43}=0.6,$ 0.8 and $1$ of the IERK3-Radau  methods \eqref{Scheme: IERK3-5stage-Radau} are very close. At the same time, this would not be the whole story. In Figures \ref{fig: IERK3-2 decay a55_tau_kappa} and \ref{fig: IERK3-4 decay ah43_tau_kappa}, we also include the discrete energy generated by the IERK2-2  method \eqref{Scheme: IERK2-one parameter} for $a_{33}=\tfrac{1+\sqrt{2}}{4}\approx 0.6036$, which seems to be the ``best'' scheme of second-order accuracy. One can observe that
the discrete energy curve generated by the IERK2-2 method is even closer to the reference energy curve than the IERK3-1 method \eqref{Scheme: IERK3-5stage-a55} with $a_{55}=1.7$. They suggest that the selection of method parameters in IERK3 methods is at least as important as the selection of different IERK3 methods. For this example, the parameters $a_{55}=0.8$ and $\hat{a}_{43}=1$ are ``good" choices for the IERK3-1 \eqref{Scheme: IERK3-5stage-a55} and IERK3-Radau \eqref{Scheme: IERK3-5stage-Radau} methods, respectively.  


\begin{figure}[htb!]
	\centering
	\subfigure[$\tau=0.01,\kappa=2$]{
		\includegraphics[width=2.05in]{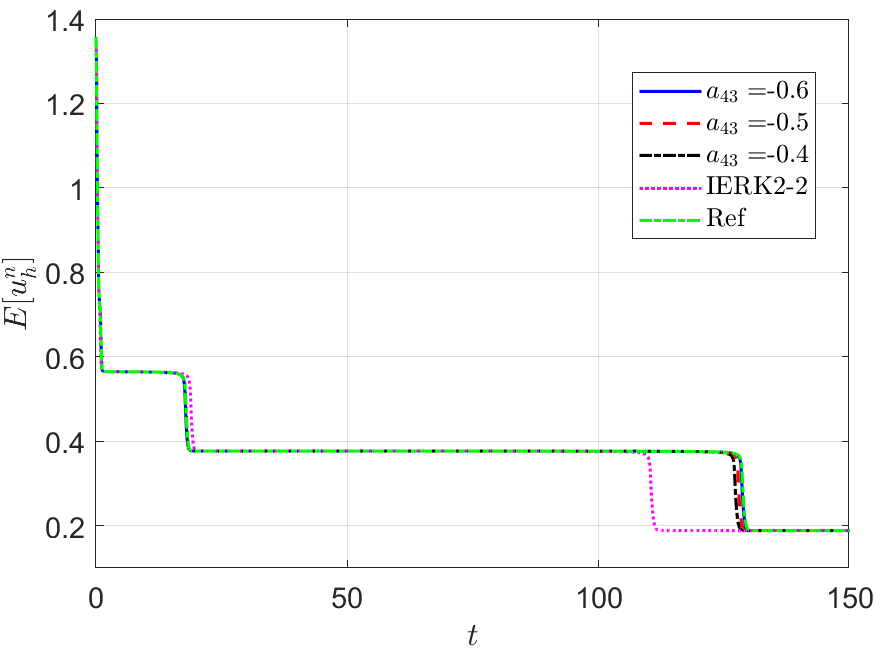}}
	\subfigure[ $\tau=0.05,\kappa=2$]{
		\includegraphics[width=2.05in]{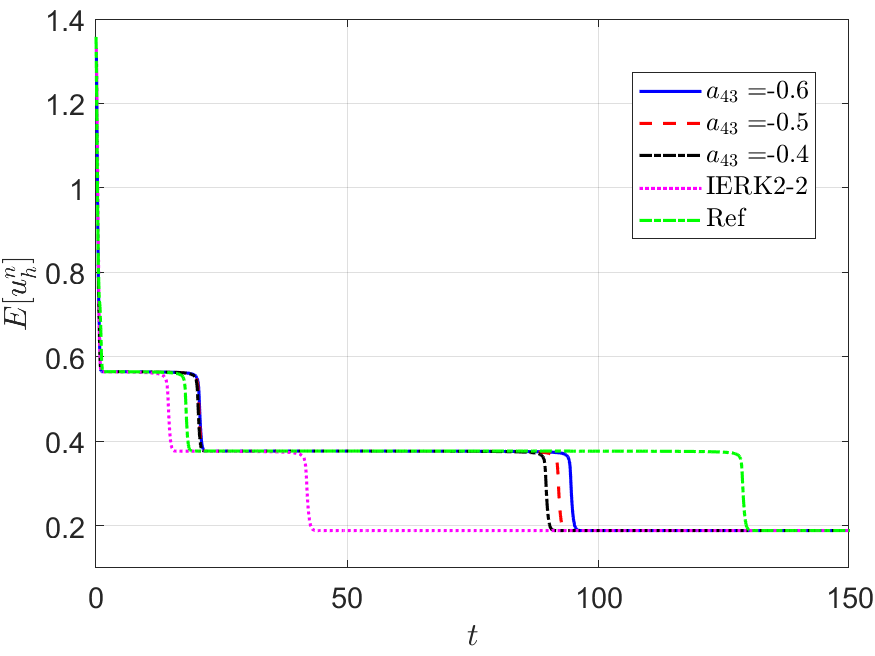}}
	\subfigure[$\tau=0.05,\kappa=3$]{
		\includegraphics[width=2.05in]{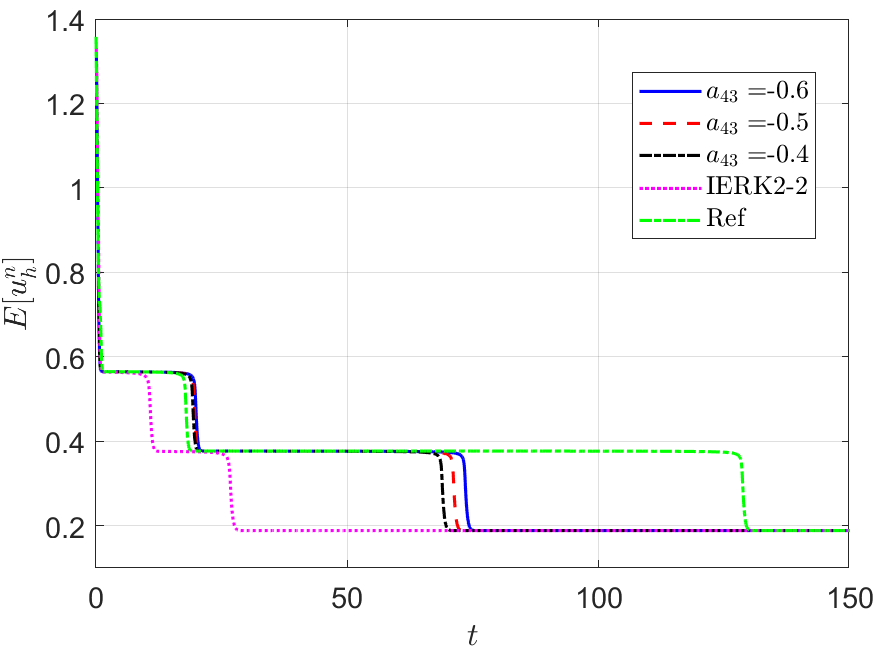}}
	\caption{Energy behaviors of IERK3-2 methods \eqref{Scheme: IERK3-5stage-a43}.}
	\label{fig: IERK3-3 decay a43_tau_kappa}
\end{figure}

\begin{figure}[htb!]
  \centering
  \subfigure[$\tau=0.01,\kappa=2$]{
    \includegraphics[width=2in]{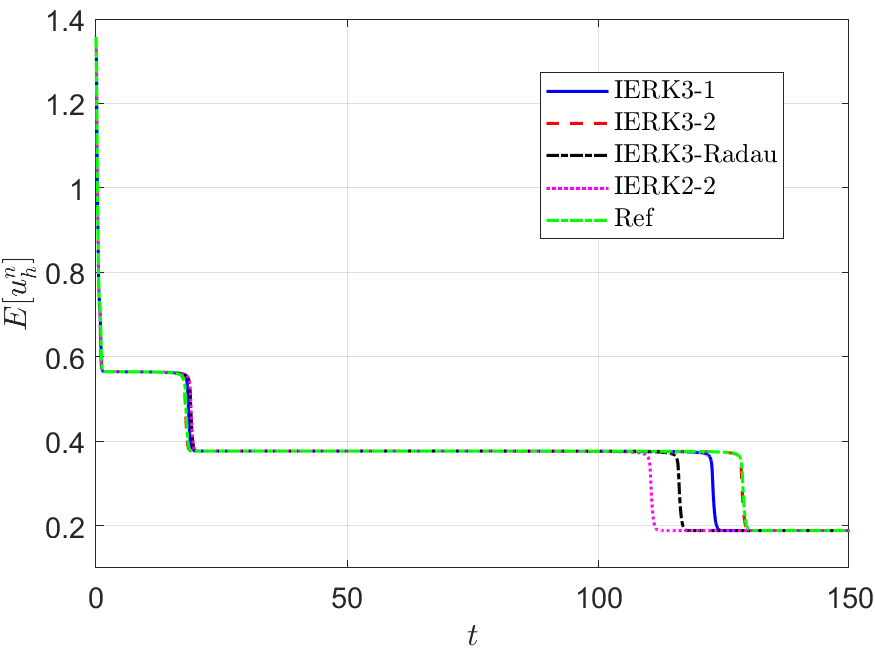}}
  \subfigure[$\tau=0.05,\kappa=2$]{
    \includegraphics[width=2in]{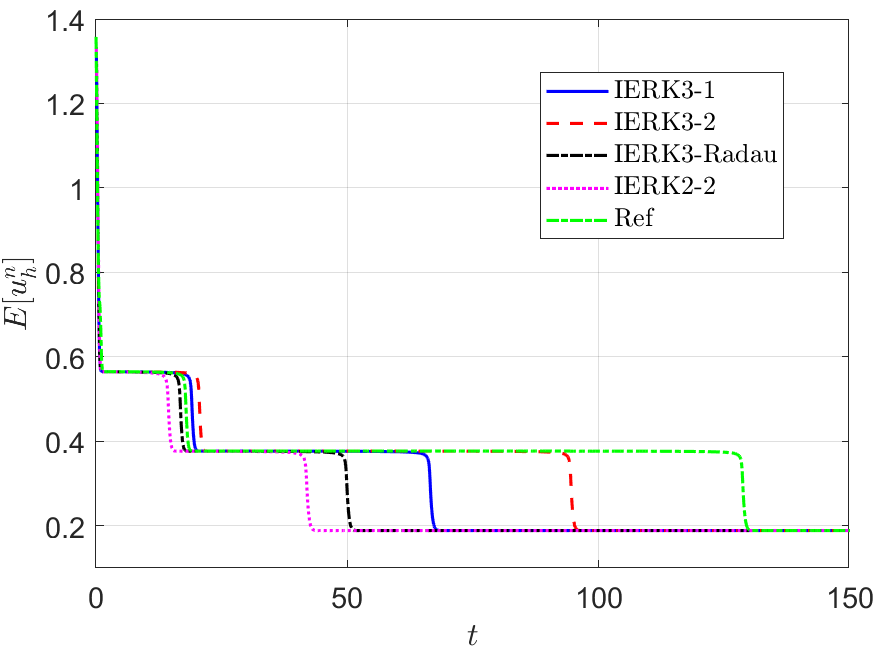}}
      \subfigure[$\tau=0.05,\kappa=3$]{
      \includegraphics[width=2in]{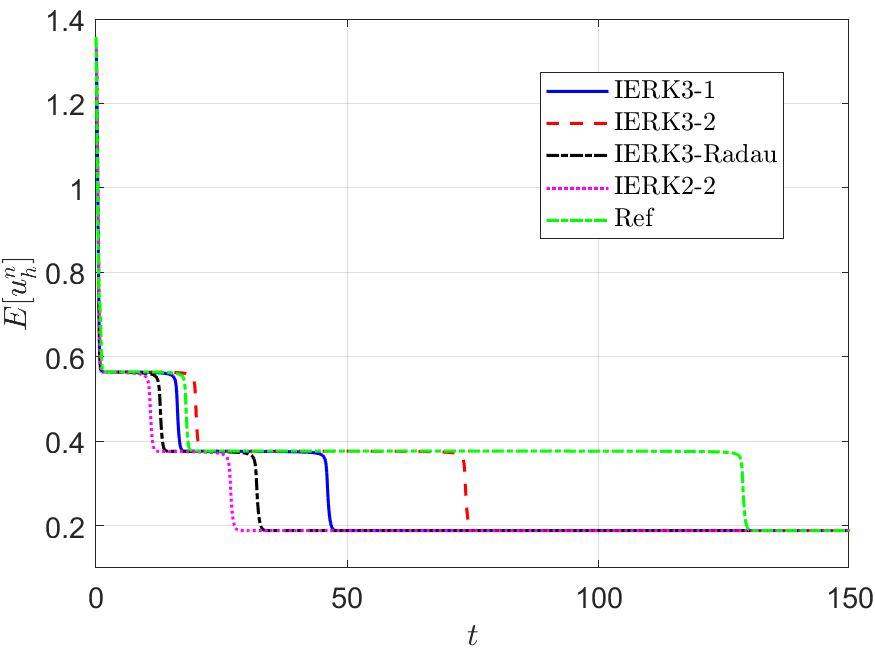}}
  \caption{Energy behaviors of IERK3-1 method for $a_{55}=0.8$, IERK3-2 method for $a_{43}=-0.6$, and IERK3-Radau method for $\hat{a}_{43}=1$.}
  \label{fig: IERK3 decay tau_kappa}
\end{figure}

We run the IERK3-2  methods \eqref{Scheme: IERK3-5stage-a43} with three parameters $a_{43}=-0.6, -0.5$ and $-0.4$ up to the time $T=150$, and exhibit the discrete energy curves in Figure \ref{fig: IERK3-3 decay a43_tau_kappa} for three scenes: (i) $\tau=0.01,\kappa=2$; (ii) $\tau=0.05,\kappa=2$ and (iii) $\tau=0.05,\kappa=3$. As predicted by Theorem \ref{Scheme: IERK3-5stage-a43}, the discrete energies in all cases always monotonically decreasing; however, the dissipation rates of discrete energies are quite different from those in Figures \ref{fig: IERK3-2 decay a55_tau_kappa} and \ref{fig: IERK3-4 decay ah43_tau_kappa}. We see that the discrete energies for three parameters $a_{43}=-0.6, -0.5$ and $-0.4$ approach the reference curve in a staggered manner (fast and slow) at different time periods. For example, the energy curves for $a_{43}=-0.6$ are always closest to the reference energy in the third fast varying period, while they are not the closest one to the reference energy in the second fast varying period. This phenomenon cannot be directly understood by the parameter-independent dissipation rate, $\mathcal{R}^{(3,2)}=\tfrac{5}{4} + \tfrac{2}{5}\tau\overline{\lambda}_{\mathrm{ML}}$, and remains mysterious to us. As expected, the IERK3-2  methods \eqref{Scheme: IERK3-5stage-a43} perform better than the ``best" second-order method \eqref{Scheme: IERK2-one parameter}.

Furthermore, the IERK3-2  methods \eqref{Scheme: IERK3-5stage-a43} would be the ``best" candidates among the three classes of IERK3 methods. In Figure \ref{fig: IERK3 decay tau_kappa}, we compare the energy behaviors of three ``good" schemes, namely, the IERK3-1 method for $a_{55}=0.8$, the IERK3-2 method for $a_{43}=-0.6$, and the IERK3-Radau method for $\hat{a}_{43}=1$, in previous tests shown in Figures \ref{fig: IERK3-2 decay a55_tau_kappa}-\ref{fig: IERK3-3 decay a43_tau_kappa}. The discrete energy curve generated by the IERK3-1 method is always closer to the reference energy curve than that of the IERK3-Radau method. Thanks to the minimum dissipation rate $\mathcal{R}^{(3,2)}=\tfrac{5}{4} + \tfrac{2}{5}\tau\overline{\lambda}_{\mathrm{ML}}$ among the three methods, the energy curve of the IERK3-2 method for $a_{43}=-0.6$ is always the closest one to the reference curve. By contrast, the discrete energy generated by the IERK2-2 method stays farthest away from the reference energy. Experimentally, they support our theoretical results.

To end this section, we summarize the theoretical results in Theorems \ref{thm: IERK3-5stage-a55s}-\ref{thm: IERK3-5stage-Radau}, the numerical results in this section and results of some existing IERK algorithms in the literature  in Table \ref{table: comparison of 3rd-order methods}, which presents the average dissipation rates and some practical choices.

\begin{table}[htb!]
	\centering
  \begin{threeparttable}
    \centering \small
    \renewcommand\arraystretch{1.3}
    \belowrulesep=0pt\aboverulesep=0pt
    \caption{Average dissipation rate of third-order IERK methods.}
    \label{table: comparison of 3rd-order methods}
    \vspace{2mm}
    \begin{tabular}{c|c|c|c}
      \toprule 
      Type & Method & Average dissipation rate & Practical choice \\
      \midrule 
      \multirow{4}*{Lobatto-type} & 5-stage IERK3-1 \eqref{Scheme: IERK3-5stage-a55} & $\tfrac{5}{4}+\tfrac{5a_{55}-2}{4}\tau\overline{\lambda}_{\mathrm{ML}}$ & $a_{55}=0.8$ \\[2pt]
      \cline{2-4}
      & 5-stage IERK3-2 \eqref{Scheme: IERK3-5stage-a43} & $\tfrac{5}{4} + \tfrac{2}{5}\tau\overline{\lambda}_{\mathrm{ML}}$ & $a_{43}=-0.6$ \\[2pt]
      \cline{2-4} & 4-stage IERK3 \cite{CooperSayfy:1983} & NPD \\[2pt]
      \cline{2-4} & 5-stage IERK3 \cite{LiuZou:2006} & NPD \\[2pt]
      \midrule
      \multirow{4}*{Radau-type} & 5-stage IERK3-Radau \eqref{Scheme: IERK3-5stage-Radau} & {\tiny$\tfrac{1}{4\hat{a}_{43}}+3.49891+\brab{\tfrac{1}{5 \hat{a}_{43}}+2.29913}\tau\overline{\lambda}_{\mathrm{ML}}$} & $\hat{a}_{43} =1.0$ \\[2pt]
      \cline{2-4} & 5-stage IERK3 \cite{FuTangYang:2024} & $3.24727+16.3779 \tau\overline{\lambda}_{\mathrm{ML}}$ \\[2pt]
      \cline{2-4} & 7-stage IERK3 \cite{ShinLeeLee:2017} & $2 + \tfrac{1}{2}\tau\overline{\lambda}_{\mathrm{ML}}$ \\[2pt]
      \cline{2-4} & 5-stage IERK3 \cite{AscherRuuthSpiteri:1997} & NPD \\[2pt]
      \bottomrule
    \end{tabular}
    \footnotesize
  \end{threeparttable}
\end{table}

\section{Discrete energy laws of IERK4-A methods}\label{sec: IERK4-A}
\setcounter{equation}{0}

\subsection{IERK4-A methods}

For six-stage IERK methods that satisfy the canopy node condition and the two order conditions for first-order accuracy, there are 29 coefficients to be determined by 26 order conditions up to fourth-order accuracy in Table \ref{table: order condition}. That is to say, there are three independent coefficients; however, we are not able to find any fourth-order IERK methods within six stages to fulfill the assumptions of Corollary \ref{corol: energy stability}, that is, the positive semi-definiteness of the two matrices $D_{\mathrm{E}}$ and $D_{\mathrm{EI}}$.

Unfortunately, we also fail to find energy stable fourth-order IERK methods with seven stages that fulfill all order conditions in Table \ref{table: order condition}. Alternatively, one can present some \emph{approximately fourth-order} IERK (IERK4-A) methods,
where the method coefficients are computed by \emph{the least squares approximation} of 18 order conditions for fourth-order accuracy with a given tolerance $\varepsilon_{\text{tol}}$, such as $\varepsilon_{\text{tol}} = 10^{-6}$. It is to note that, for the IERK4-A methods satisfying the canopy node condition and 10 order conditions up to third-order accuracy, the truncation errors are about of $O(\varepsilon_{\text{tol}} \tau^{3} + \tau^{4})$ if the solution is sufficiently regular. That is to say, a lower bound $\tau\geq O(\varepsilon_{\text{tol}})$ will be required to attain the fourth-order accuracy of IERK4-A methods although such time-step restriction is rather weak in practical applications.

The first one of fourth-order accuracy with  the tolerance $\varepsilon_{\text{tol}} = 10^{-6}$, called IERK4-A1 method, has the following Butcher tableaux
\begin{align}\label{Scheme: IERK4-A1}
  &\begin{array}{c|c}
    \mathbf{c} & A \\
    \hline\\[-10pt]   & \mathbf{b}^T
  \end{array}
  =\begin{array}{c|ccccccc}
   0 & 0  \\[3pt]
  \frac{95341769}{200000000} & a_{21} & \tfrac{2400249}{2000000}  \\[3pt]
  \frac{292103}{800000} & a_{31} & -\tfrac{173504613}{50000000} & \tfrac{2486}{625}  \\[3pt]
  \frac{59556813}{200000000} & a_{41} & \tfrac{13944041}{20000000} & a_{43} & \tfrac{1326491}{1000000}  \\[3pt]
  \frac{2580667}{5000000} & a_{51} & -\tfrac{92214113}{200000000} & -\tfrac{942329}{1000000} & \tfrac{11063869}{10000000} & \tfrac{1185669331}{2000000000}  \\[3pt]
  \frac{150085929}{200000000} & a_{61} & -\tfrac{2010707}{2500000} & -\tfrac{151161011}{200000000} & a_{64} & -\tfrac{23680671}{40000000} & \tfrac{15240463}{20000000}  \\[3pt]
  1 & a_{71} & \tfrac{188918701}{250000000} & -\tfrac{205244563}{500000000} & -\tfrac{463474901}{250000000} & a_{75} & a_{76} & \tfrac{1339351}{2000000} \\[3pt]
  \hline\\[-10pt] & a_{71} & \tfrac{188918701}{250000000} & -\tfrac{205244563}{500000000} & -\tfrac{463474901}{250000000} & a_{75} & a_{76} & \tfrac{1339351}{2000000}
  \end{array}\notag\\
  &\begin{array}{c|c}
    \hat{\mathbf{c}} & \widehat{A} \\
    \hline\\[-10pt]   & \hat{\mathbf{b}}^T
  \end{array}
  = \begin{array}{c|ccccccc}
    0 & 0  \\[3pt]
  \frac{95341769}{200000000} & \hat{a}_{21} & 0  \\[3pt]
  \frac{292103}{800000} & \hat{a}_{31} & \tfrac{558887}{2000000} & 0 \\[3pt]
  \frac{59556813}{200000000} & \hat{a}_{41} & -\tfrac{13454127}{250000000} & \tfrac{121}{1000} & 0  \\[3pt]
  \frac{2580667}{5000000} & \hat{a}_{51} & \hat{a}_{52} & \tfrac{32018089}{200000000} & \tfrac{593979}{2000000} & 0  \\[3pt]
  \frac{150085929}{200000000} & \hat{a}_{61} & \tfrac{18031311}{200000000} & \hat{a}_{63} & \tfrac{1177953}{10000000} & \tfrac{313441}{1000000} & 0  \\[3pt]
  1 & \hat{a}_{71} & -\tfrac{23838287}{200000000} & \tfrac{31003}{160000} & \hat{a}_{74} & \hat{a}_{75} & \tfrac{70189993}{100000000} & 0 \\[3pt]
  \hline\\[-10pt] & \hat{a}_{71} & -\tfrac{23838287}{200000000} & \tfrac{31003}{160000} & \hat{a}_{74} & \hat{a}_{75} & \tfrac{70189993}{100000000} & 0
  \end{array}
\end{align}
where the coefficients $a_{k,1}=c_k-\sum_{j=2}^{k} a_{k,j}, \hat{a}_{k,1}=c_{k}-\sum_{j=2}^{k-1}\hat{a}_{k,j}$ for $2\leq k \leq 7$ and
\begin{align*}
  a_{43}=&\,-\tfrac{1585409690050693626522959}{10^{24}},&\quad& a_{64} = \tfrac{2344693633530028154195338}{10^{24}}, 
  \\
  a_{75} =&\, \tfrac{12541381954770160599120029}{3627843526172490000000000}, &\quad& a_{76} = -\tfrac{5666098495504076766418243}{2637342719402745375000000},\\[2pt]
  \hat{a}_{52}=&\,\tfrac{33045877519515636367149}{10^{25}}, &\quad&  \hat{a}_{63} = \tfrac{983273174147729070884395}{10^{25}}, \\
  \hat{a}_{74}=&\,\tfrac{952796123534851512831817}{1560502861592322600000000}, &\quad& \hat{a}_{75} = -\tfrac{145732093978331037774401}{338092153983867000000000}.
\end{align*}
The associated matrices $D_{\mathrm{E}}^{(4,1)}$ and $D_{\mathrm{EI}}^{(4,1)}$ defined by \eqref{Def: Differential Matrix D} are positive definite with the eigenvalue vectors, respectively,
\begin{align*}
  &\brab{12.4798, 4.72381, 3.05289, 1.11175, 0.530827, 0.0238641}^T\\
  \text{and}\quad & \brab{14.4357, 10.3404, 2.47970, 1.62679, 1.21228, 0.000779}^T.
\end{align*}
The corresponding average dissipation rate is
$$\mathcal{R}^{(4,1)} 
\approx 3.65382+ 5.01594 \tau\overline{\lambda}_{\mathrm{ML}}.$$

The second one of fourth-order accuracy  with the tolerance $\varepsilon_{\text{tol}} = 10^{-6}$, called IERK4-A2 method, has the following Butcher tableaux
\begin{align}\label{Scheme: IERK4-A2}
  &\begin{array}{c|c}
    {\mathbf{c}} & {A} \\
    \hline\\[-10pt]   & {\mathbf{b}}^T
  \end{array}
  =  \begin{array}{c|ccccccc}
  0 & 0  \\[3pt]
  \tfrac{429533}{1000000} & a_{21} & \tfrac{63137}{200000}  \\[3pt]
  \tfrac{4785663}{10000000} & a_{31} & -\tfrac{917757}{1000000} & \tfrac{100379}{100000}  \\[3pt]
  \tfrac{1182276}{1000000} & a_{41} & -\tfrac{1929}{1250} & a_{43} & \tfrac{15281}{20000}  \\[3pt]
  \tfrac{915703}{1000000} & a_{51} & \tfrac{98637}{1000000} & \tfrac{196933}{1000000} & a_{54} & \tfrac{531}{625}  \\[3pt]
  \tfrac{7336053}{10000000} & a_{61} & \tfrac{302663}{1000000} & \tfrac{5967}{125000} & \tfrac{150781}{1000000} & -\tfrac{156231}{125000} & \tfrac{142387}{100000}  \\[3pt]
  1 & a_{71} & \tfrac{1843487}{10000000} & a_{73} & -\tfrac{129513}{1000000} & -\tfrac{820173}{2000000} & a_{76} & \tfrac{88673}{125000} \\[3pt]
  \hline\\[-10pt] & a_{71} & \tfrac{1843487}{10000000} & a_{73} & -\tfrac{129513}{1000000} & -\tfrac{820173}{2000000} & a_{76} & \tfrac{88673}{125000}
\end{array}\,,\notag\\
&\begin{array}{c|c}
    \hat{\mathbf{c}} & \widehat{A} \\
    \hline\\[-10pt]   & \hat{\mathbf{b}}^T
  \end{array}
  =  \begin{array}{c|ccccccc}
   0 & 0  \\[3pt]
  \tfrac{429533}{1000000} & \hat{a}_{21} & 0 \\[3pt]
  \tfrac{4785663}{10000000} & \hat{a}_{31} & \hat{a}_{32} & 0  \\[3pt]
  \tfrac{1182276}{1000000} & \hat{a}_{41} & \tfrac{4507}{6250} & \tfrac{1025153}{2000000} & 0 \\[3pt]
  \tfrac{915703}{1000000} & \hat{a}_{51} & \tfrac{587731}{2000000} & \hat{a}_{53} & \tfrac{371251}{2000000} & 0  \\[3pt]
  \tfrac{7336053}{10000000} & \hat{a}_{61} & \tfrac{82583}{200000} & -\tfrac{1415767}{10000000} & -\tfrac{33393}{200000} & \tfrac{119457}{250000} & 0  \\[3pt]
  1 & \hat{a}_{71} & \tfrac{28277}{100000} & \hat{a}_{73} & \hat{a}_{74} & -\tfrac{7683}{100000} & \tfrac{99459}{250000} & 0 \\[3pt]
  \hline\\[-10pt] & \hat{a}_{71} & \tfrac{28277}{100000} & \hat{a}_{73} & \hat{a}_{74} & -\tfrac{7683}{100000} & \tfrac{99459}{250000} & 0
\end{array}\,,
\end{align}
where the first coefficients $a_{k,1}=c_k-\sum_{j=2}^{k}a_{k,j}$ and $\hat{a}_{k,1}=c_k-\sum_{j=2}^{k-1}\hat{a}_{k,j}$ for $2\le k\le 7$, and
\begin{align*}
  a_{43} = &\, \tfrac{219830108841087453607347}{2\times10^{23}}, &\quad& a_{54} = -\tfrac{4498694297454501655541833}{10^{25}}, \\  
  a_{73} = &\, \tfrac{2298242610563399947}{4576990146963750000}, &\quad& a_{76} = -\tfrac{556251214988653}{1754043394750312500}, \\[2pt]
  \hat{a}_{32} =&\, \tfrac{1017529648895428446045183}{25\times10^{23}}, &\quad& \hat{a}_{53} = \tfrac{3161854834370097094143699}{1\times10^{25}}, \\
  \hat{a}_{73} =&\, \tfrac{4063870960730480933}{25257881055233250000}, &\quad& \hat{a}_{74} = \tfrac{422441222477275261}{6239843169579000000}.
\end{align*}
The associated matrices $D_{\mathrm{E}}^{(4,2)}$ and $D_{\mathrm{EI}}^{(4,2)}$ defined by \eqref{Def: Differential Matrix D} are positive definite with the eigenvalue vectors, respectively,
\begin{align*}
  &\brab{8.31573, 4.30483, 1.87790, 1.55745, 0.673643, 0.0000139}^T\\
  \text{and}\quad & \brab{5.12769, 3.33053, 1.57402, 0.952105, 0.0473756, 0.0000135}^T.
\end{align*}
The corresponding average dissipation rate is
$\mathcal{R}^{(4,2)} 
\approx 2.78826+ 1.83862 \tau\overline{\lambda}_{\mathrm{ML}}.$

By using Corollary \ref{corol: energy stability} and Lemma \ref{lemma: average dissipation rate}, we have the following theorem.
\begin{theorem}\label{thm: IERK4-A} 
  In simulating the gradient flow system \eqref{problem: stabilized version} with  $\kappa\ge2\ell_g$, the seven-stage IERK4-1 \eqref{Scheme: IERK4-A1} and IERK4-2 \eqref{Scheme: IERK4-A2} methods preserve the original energy dissipation law \eqref{problem: energy dissipation law},
  \begin{align*}    
    E[U^{n,j+1}]-E[U^{n,1}]\le\frac1{\tau}\myinnerB{M_h^{-1}\delta_{\tau}\vec{U}_{n,j+1},
      D_{j}^{(4,\mu)}({\tau}M_hL_{\kappa})\delta_{\tau}\vec{U}_{n,j+1}}
    \quad\text{for $1\le j\le 6$,}
  \end{align*}
  where  $\mu=1,2$ and the  differentiation  matrices $D^{(4,\mu)}$ are defined by  $$D^{(4,\mu)}({\tau}M_hL_{\kappa}):=D_{\mathrm{E}}^{(4,\mu)}\otimes I+D_{\mathrm{EI}}^{(4,\mu)}\otimes (-\tau M_h L_{\kappa}).$$
  The associated average dissipation rates are
  \begin{align}\label{AverRate: seven-stage IERK4}
    \mathcal{R}^{(4,1)}
    \approx 3.65382+ 5.01594 \tau\overline{\lambda}_{\mathrm{ML}} \quad\text{and}\quad 
    \mathcal{R}^{(4,2)} \approx   2.78826+ 1.83862 \tau\overline{\lambda}_{\mathrm{ML}}.
  \end{align}
\end{theorem}

The seven-stage IERK4-1 \eqref{Scheme: IERK4-A1} and IERK4-2 \eqref{Scheme: IERK4-A2} methods are considered only to show the existence of the energy
stable methods with fourth-order time accuracy.
As mentioned early, Liu and Zou \cite{LiuZou:2006} proposed several Lobatto-type fourth-order  IERK methods; but these methods have zero diagonal entries and do not satisfy our assumptions in Corollary \ref{corol: energy stability}. Actually, in the literature, we do not find any  Lobatto-type or Radau-type fourth-order IERK methods that the associated differentiation matrices are positive (semi-)definite.  

\subsection{Tests for IERK4-A methods}

\begin{figure}[htb!]
	\centering
	\subfigure[$\kappa=2$]{
		\includegraphics[width=2.1in]{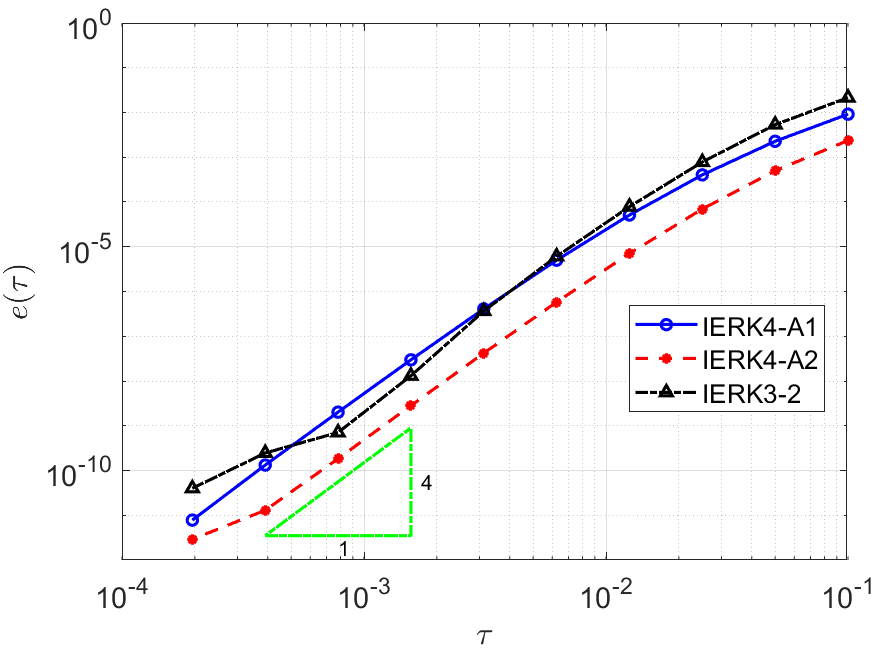}} 
	\subfigure[$\kappa=3$]{
		\includegraphics[width=2.1in]{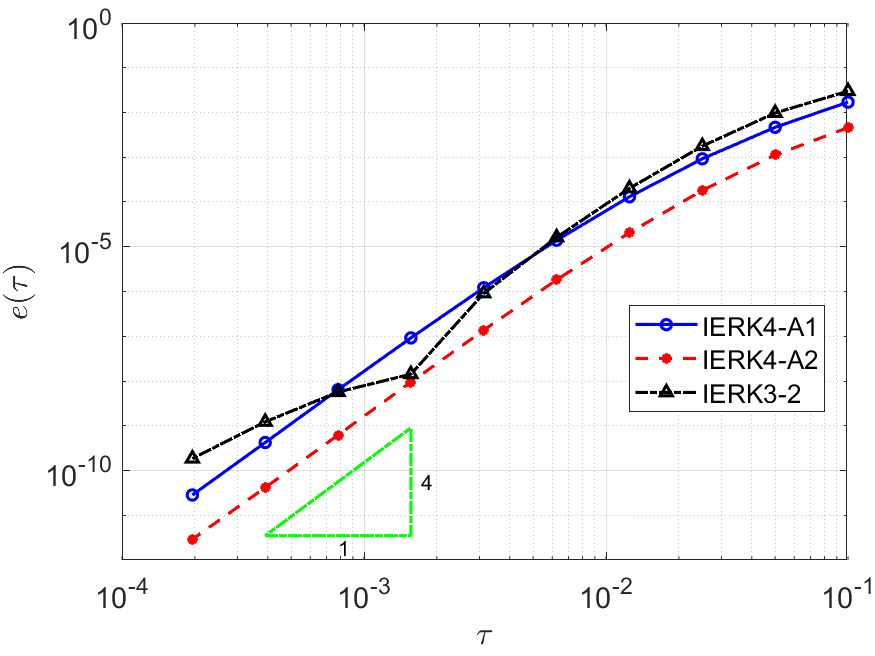}}
	\caption{Solution errors of IERK4 methods.}
	\label{fig: IERK4 convergence}
\end{figure}

Example \ref{ex: with source term} is also used to test the convergence of IERK4-A1 \eqref{Scheme: IERK4-A1} and IERK4-A2 \eqref{Scheme: IERK4-A2} methods by choosing the final time $T=1$ with different stabilized parameters $\kappa=2, 3$. The $L^\infty$ norm errors $e(\tau)$ of the two IERK4-A methods on halving time steps $\tau=2^{-k}/10$ for $0\le k\le9$ are shown in Figure \ref{fig: IERK4 convergence}. We see that the two IERK4-A methods can achieve fourth-order accuracy and the  error of IERK4-A2 method is always smaller than the error of IERK4-A1 method. 

\begin{figure}[htb!]
	\centering
	\subfigure[$\tau=0.01,\kappa=2$]{
		\includegraphics[width=2in]{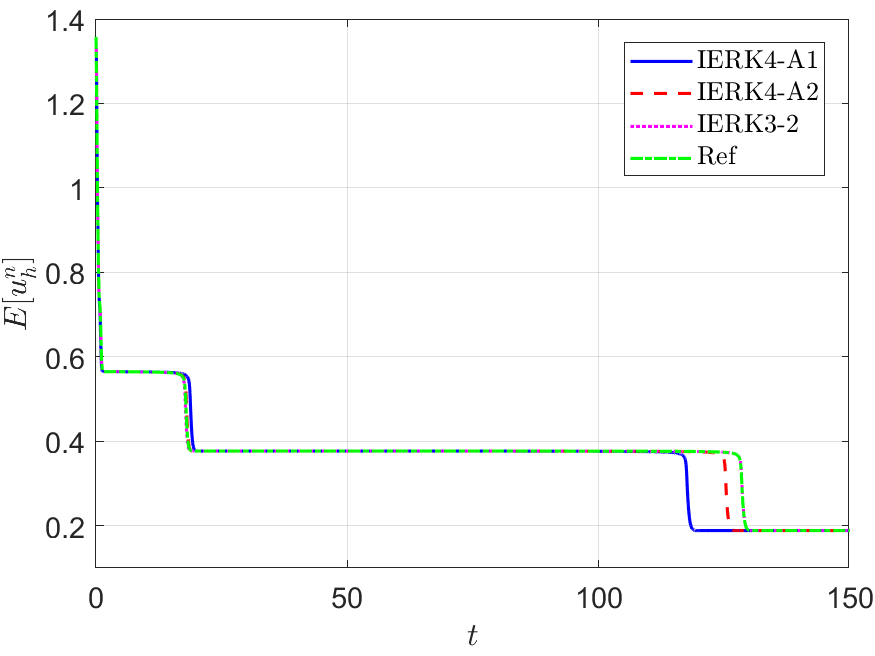}} 
	\subfigure[$\tau=0.05,\kappa=2$]{
		\includegraphics[width=2in]{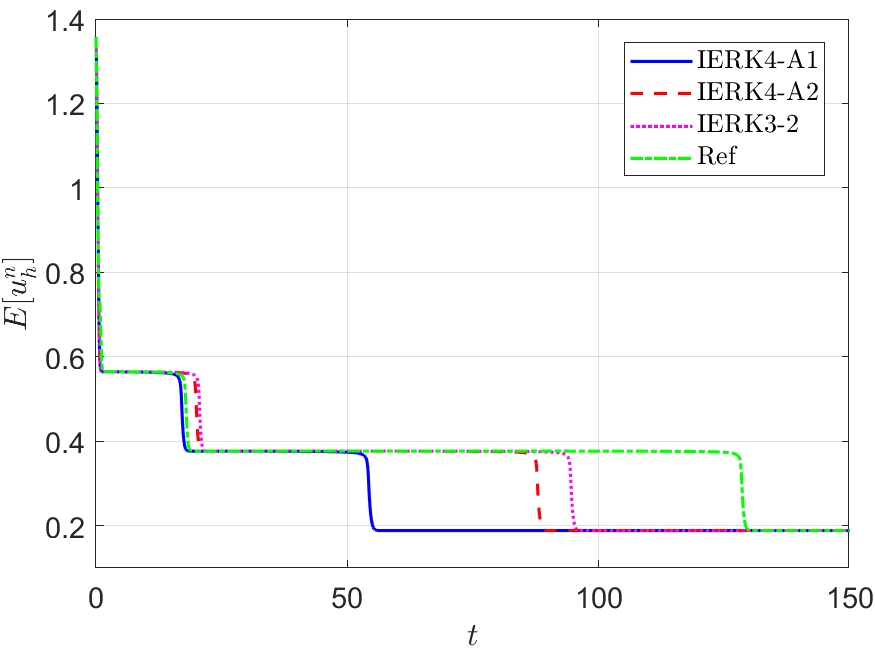}}
	\subfigure[$\tau=0.05,\kappa=3$]{
		\includegraphics[width=2in]{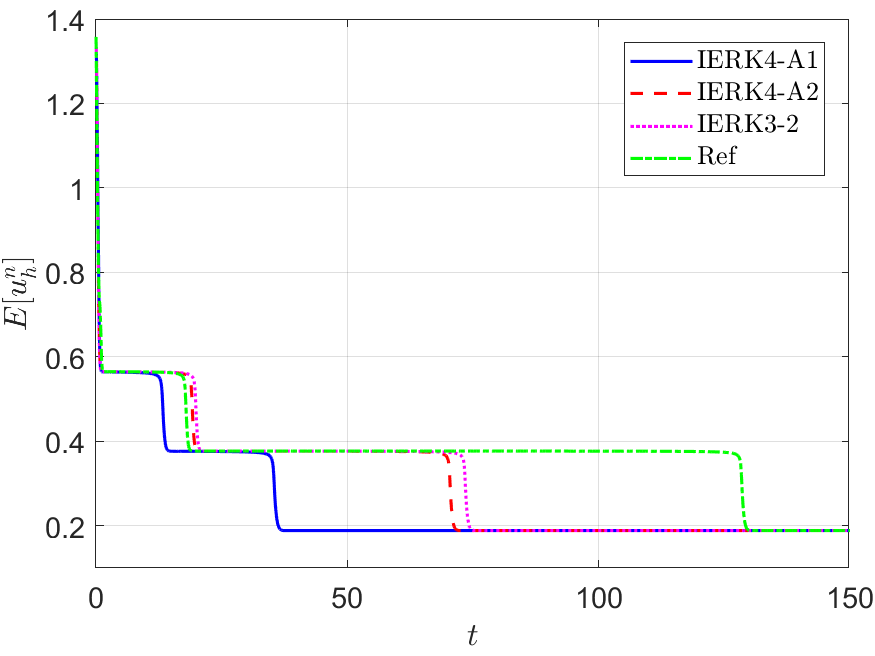}}
	\caption{Energy behaviors of IERK4-A1 and IERK4-A2 methods.}
	\label{fig: IERK4 decay tau_kappa}
\end{figure}

The discrete energy behaviors of the two IERK4-A methods are examined by Example \ref{example: energy_2}, in which the reference solution and energy are computed by the IERK4-A1 method with a small time-step size $\tau=10^{-3}$. Figure \ref{fig: IERK4 decay tau_kappa} depicts the discrete energies generated by the IERK4-A1 and IERK4-A2 methods, for three different scenes: (i) $\tau=0.01,\kappa=2$; (ii) $\tau=0.05,\kappa=2$ and (iii) $\tau=0.05,\kappa=3$. As predicted by Theorem \ref{thm: IERK4-A} , the discrete energies in all cases always monotonically decrease. Also, the discrete energy obtained from the IERK4-A2 method is always closer to the reference energy than that of the IERK4-A1 method. It is predictable by the numerical errors in Figure \ref{fig: IERK4 convergence}, which shows that the IERK4-A2 method is always more accurate than the IERK4-A1 method, and the associated average dissipation rates in \eqref{AverRate: seven-stage IERK4}, that is, $\abs{\mathcal{R}^{(4,2)}-1}<\abs{\mathcal{R}^{(4,1)}-1}$.

In Figure \ref{fig: IERK4 decay tau_kappa} (a)-(c), we also include the discrete energy generated by the IERK3-2  method \eqref{Scheme: IERK3-5stage-a43} with $a_{43}=-0.6$, which has the average dissipation rate  $\mathcal{R}^{(3,2)}=\tfrac{5}{4} + \tfrac{2}{5}\tau\overline{\lambda}_{\mathrm{ML}}$ and is regarded as the ``best" IERK3 method in this paper. More surprisingly, the energy curve generated by the IERK3-2 method \eqref{Scheme: IERK3-5stage-a43} with $a_{43}=-0.6$ is always closer to the reference one than the two IERK4-A schemes. Actually, this phenomenon can not be explained by the precision of numerical solutions. As seen in Figure \ref{fig: IERK4 convergence}, the solution errors of IERK3-2 method are larger than the errors of IERK4-A1 and IERK4-A2 methods for the time-step sizes $\tau=0.01$ and $0.05$. Obviously, the well preservation of the original energy dissipation law does not entirely depend on the numerical precision of solution  
although we are unsure of the complete mechanism behind it. One possible reason is that the average dissipation rates in \eqref{AverRate: seven-stage IERK4} of the two IERK4-A schemes are relatively large. 

  This raises an important issue: finding some IERK4-A method or formally fourth-order IERK method that preserves the original energy dissipation law \eqref{problem: energy dissipation law} and perform better than the ``best'' third-order method \eqref{Scheme: IERK3-5stage-a43} with $a_{43}=-0.6$. We have planned to continue exploring this issue and will present the relevant results in a forthcoming report.

\section{Concluding remarks}
\label{section: conclusions}\setcounter{equation}{0}

We construct some parameterized IERK methods by adopting the stiffly-accurate Lobatto-type DIRK methods in the implicit part, and present a unified theoretical framework to examine the energy dissipation properties at all stages of IERK methods up to fourth-order accurate for gradient flow problems.
The novel framework contains two main parts: one is the differential form and the associated differentiation  matrix of an IERK method by using the difference coefficients of method and the so-called discrete orthogonal convolution kernels; the other is the average energy dissipation rate defined via the average eigenvalue of the differentiation matrix. The rough but simple concept of average energy dissipation rate seems very useful to evaluate the overall energy dissipation behaviors of the IERK methods, including choosing proper parameters in some parameterized IERK methods or comparing different IERK methods with the same accuracy. 
Our main contributions include:
\begin{itemize}
  \item[(i)] Among the second-order IERK2-1 \eqref{Scheme: IERK2-two parameters}, IERK2-2 \eqref{Scheme: IERK2-one parameter} and IERK2-Radau \eqref{Scheme: IERK2-Radau-one parameter} methods preserving the original energy dissipation law \eqref{problem: energy dissipation law}, the one-parameter IERK2-2 method \eqref{Scheme: IERK2-one parameter} with  $a_{33}=\tfrac{1+\sqrt{2}}{4}$ would be the best one from both the average energy dissipation rate $\mathcal{R}^{(2,2)}(\tfrac{\sqrt{2}}{2},\tfrac{1+\sqrt{2}}{4})
  =\sqrt{2}+\tfrac{\sqrt{2}}{4}\tau\overline{\lambda}_{\mathrm{ML}}$, and the numerical experiments in Subsection 3.2.
  \item[(ii)] Among the third-order IERK3-1 \eqref{Scheme: IERK3-5stage-a55}, IERK3-2 \eqref{Scheme: IERK3-5stage-a43} and IERK3-Radau \eqref{Scheme: IERK3-5stage-Radau} methods preserving the original energy dissipation law \eqref{problem: energy dissipation law} unconditionally, the one-parameter IERK3-2 method  \eqref{Scheme: IERK3-5stage-a43} with the parameter $a_{43}=-0.6$ would be the best one from both the average energy dissipation rate $\mathcal{R}^{(3,2)}
  =\frac{5}4+\tfrac{2}{5}\tau\overline{\lambda}_{\mathrm{ML}}$, and the numerical tests in Subsection 4.4. 
  \item[(iii)] Two approximately fourth-order IERK methods with the  order-condition tolerance $\varepsilon_{\text{tol}} = 10^{-6}$ are constructed and shown to  preserve the original energy dissipation law \eqref{problem: energy dissipation law} unconditionally. From the average dissipation rates and the numerical experiments in Subsection 5.2, the IERK4-A2 method \eqref{Scheme: IERK4-A2} is better than the IERK4-A1 method \eqref{Scheme: IERK4-A1}. 
  \item[(iv)]  From the perspective of maintaining the original energy dissipation law \eqref{problem: energy dissipation law},  numerical tests show that the IERK1 method \eqref{Scheme: IERK1} with $\theta=\tfrac12$ performs better than some of IERK2 methods; the IERK2-2 method \eqref{Scheme: IERK2-one parameter} with $a_{33}=\tfrac{1+\sqrt{2}}{4}$ performs better than some of IERK3 methods and the IERK3-2 method  \eqref{Scheme: IERK3-5stage-a43} with $a_{43}=-0.6$ performs better than the two IERK4-A methods. They suggest that the time accuracy of an IERK method cannot fully characterize the degree to which it maintains the original energy dissipation law. In practice, the selection of method parameters in IERK methods is at least as important as the selection of different IERK methods, if not more important.
\end{itemize}
At the same time, our theory is far away from complete. There are many interesting issues that we have not yet addressed. Some of them are listed as follows:
\begin{itemize}
  \item[(a)] We address two approximately fourth-order IERK methods that preserve the energy dissipation law \eqref{problem: energy dissipation law} unconditionally; however, up to now, we can not construct any formally fourth-order IERK methods that preserve the energy dissipation law \eqref{problem: energy dissipation law}. 
  \item[(b)] In order to maintain the continuous dissipation rate of the original energy as much as possible, the average dissipation rate, $\mathcal{R}=\frac1{s_{\mathrm{I}}}\mathrm{tr}(D_{\mathrm{E}})+\frac1{s_{\mathrm{I}}}\mathrm{tr}(D_{\mathrm{EI}})\tau\overline{\lambda}_{\mathrm{ML}}$, is hoped to be as close to 1 as possible.   
  It is seen from \eqref{AverRate: IERK1} that the value of $\mathrm{tr}(D_{\mathrm{EI}})$ can reach zero, while we always have $\frac1{s_{\mathrm{I}}}\mathrm{tr}(D_{\mathrm{E}})\ge1$ for the presented eight IERK methods, cf. Tables \ref{table: comparison of 2nd-order methods}-\ref{table: comparison of 3rd-order methods}. 
  We wonder what the minimum value of $\frac1{s_{\mathrm{I}}}\mathrm{tr}(D_{\mathrm{E}})$ may take, and how to construct the corresponding IERK method. 
  \item[(c)] As we have repeatedly stated, the average dissipation rate can only be used to compare different IERK algorithms with the same accuracy. Is it possible to reasonably integrate the accuracy information of the algorithm to form a more effective indicator so that the discrete energy dissipation behaviors of IERK methods with different accuracies can be compared?
  \item[(d)] Last but not least,  the global Lipschitz continuity assumption in Lemma \ref{lemma: origional energy derivation} of the nonlinear function $g_h$ is limited. For specific gradient flow problems, weakening the above assumption is theoretically important, and as the closely related issue, determining the minimum stabilized parameter $\kappa$ in \eqref{def: stabilized parameter} is also practically useful.
\end{itemize}


\begin{thebibliography}{99}
  
\bibitem{AkrivisLiLi:2019}
  {\sc G. Akrivis, B. Li and D. Li},
  {Energy-decaying extrapolated RK-SAV methods
    for the Allen-Cahn and Cahn-Hilliard equations},
  {\em SIAM J. Sci. Comput.}, 41:6 (2019), pp. A3703-A3727.
  
  



  
\bibitem{AscherRuuthSpiteri:1997}
{\sc U. Ascher, S. Ruuth and R.J. Spiteri},
Implicit-explicit Runge-Kutta methods for time dependent partial differential equations,
{\em Appl. Numer. Math.}, 25 (1997), pp. 151-167.



\bibitem{BoscarinoPareschiRusso:2017}
{\sc S. Boscarino, L. Pareschi and G. Russo},
A unified IMEX Runge-Kutta approach for hyperbolic systems with multiscale relaxation,
{\em SIAM J. Numer. Anal.}, 55:4 (2017), pp. 2085-2109.


\bibitem{BoscarinoRusso:2009}
{\sc S. Boscarino and G. Russo},
On a class of uniformly accurate IMEX Runge-Kutta schemes and applications to hyperbolic systems with relaxation,
{\em SIAM J. Sci. Comput.}, 31:3 (2009), pp. 1926-1945.


\bibitem{CardoneJackiewiczSanduZhang:2014MMA}
{\sc A. Cardone, Z. Jackiewicz, A. Sandu and H. Zhang},
Extrapolated implicit-explicit Runge-Kutta methods,
{\em Math. Model. Anal.}, 19:1 (2014), pp. 18-43.

\bibitem{CooperSayfy:1983}
{\sc G. J. Cooper and A. Sayfy},
Additive Runge-Kutta methods for stiff ordinary differential equations,
{\em Math. Comp.}, 40:161 (1983), pp. 207-218.



%
%
%
%
\bibitem{DimarcoPareschi:2013}
{\sc G. Dimarco and L. Pareschi},
Asymptotic preserving implicit-explicit Runge-Kutta methods for nonlinear kinetic equations,
{\em SIAM J. Numer. Anal.}, 51:2 (2013), pp. 1064-1087.





\bibitem{DuJuLiQiao:2021SIREV}
{\sc Q. Du, L. Ju, X. Li and Z. Qiao},
Maximum bound principles for  a class of semilinear parabolic equations and exponential time-differencing schemes,
{\em SIAM Rev.}, 63 (2021), pp. 317-359.



 
\bibitem{FasiGaudreaultLundSchweitzer:2024} 
{\sc M. Fasi, S. Gaudreault, K. Lund and M. Schweitzer},
Challenges in computing matrix functions,
{\em arXiv:2401.16132v1}, 2024.

\bibitem{FuYang:2022JCP}
{\sc Z. Fu and J. Yang},
Energy-decreasing exponential time differencing Runge-Kutta methods for phase-field models, 
{\em J. Comput. Phys.}, 454 (2022), 110943.

\bibitem{FuTangYang:2024}
{\sc Z. Fu, T. Tang and J. Yang}, 
Energy diminishing implicit-explicit Runge-Kutta methods for gradient flows, 
{\em Math. Comput.}, 2024, doi: 10.1090/mcom/3950.

  
\bibitem{GongZhaoWang:2020HEQRK}
{\sc Y. Gong, J. Zhao and Q. Wang},
{Arbitrarily high-order unconditionally energy stable schemes for thermodynamically consistent gradient flow models},
{\em SIAM J. Sci. Comput.}, 42:1 (2020), pp. B135-B156.

\bibitem{GuiWangChen:2023}
{\sc Y. Gui, C. Wang and J. Chen},
IMEX-RK methods for Landau-Lifshitz equation with arbitrary damping,
{\em arXiv: 2312.15654}, 2023.

  
  
  




\bibitem{HundsdorferVerwer:2003book}
{\sc W. Hundsdorfer and J.G. Verwer},
{Numerical Solution of Time-Dependent Advection Diffusion-Reaction Equations},
Springer Ser. Comput. Math. 33, Springer-Verlag, Berlin, 2003.

\bibitem{HochbruckOstermann:2010ActaNu}
{\sc M. Hochbruck and A. Ostermann}, 
\newblock Exponential integrators,
\newblock {\em Acta Numerica.}, 19 (2010), pp. 209-286.

\bibitem{IzzoJackiewicz:2017}
{\sc G. Izzo and Z. Jackiewicz},
Highly stable implicit-explicit Runge-Kutta methods,
{\em Appl. Numer. Math.}, 113 (2017), pp. 71-92.



\bibitem{JiangZhangZhao:2022}
{\sc M. Jiang, Z. Zhang and J. Zhao},
Improving the accuracy and consistency of the scalar auxiliary variable (SAV) method with relaxation,
{\em J. Comput. Phys.}, 456 (2022), 110954.

\bibitem{Jin:1995}
{\sc S. Jin},
Runge-Kutta methods for hyperbolic conservation laws with stiff relaxation terms,
{\em J. Comput. Phys.}, 122:1 (1995), pp. 51-67.

\bibitem{KanevskyCarpenterGottliebHesthaven:2007}
{\sc A. Kanevsky, M.H. Carpenter, D. Gottlieb and J.S. Hesthaven},
Application of implicit-explicit high order Runge-Kutta methods to discontinuous-Galerkin schemes,
{\em J. Comput. Phys.}, 225:2 (2007), pp. 1753-1781.


%

%
%
%
%


\bibitem{KennedyCarpenter:2003}
{\sc C.A. Kennedy and M.H. Carpenter}, Additive Runge-Kutta schemes for convection-diffusion-reaction
equations,
{\em Appl. Numer. Math.}, 44 (2003), pp. 139-181.

\bibitem{KennedyCarpenter:2016}
{\sc C.A. Kennedy and M.H. Carpenter},
{Diagonally implicit Runge-Kutta methods for ordinary differential equations: a review},
Technical Memorandum: NASA/TM 2016-219173, 2016.



%

\bibitem{LiLiao:2022}
{\sc Z. Li and H.-L. Liao},
{Stability of variable-step BDF2 and BDF3 methods},
{\em SIAM J. Numer. Anal.}, 60:4 (2022), pp. 2253-2272.


  
\bibitem{LiaoTangZhou:2024}
{\sc H.-L. Liao, T. Tang and T. Zhou},
{Positive definiteness of real quadratic forms resulting from variable-step L1-type approximations of convolution operators},
{\em Sci. China. Math.}, 67:2 (2024), pp. 237-252.
  
\bibitem{LiaoWang:2024arxiv}
{\sc H.-L. Liao and X. Wang},
{Average energy dissipation rates of explicit exponential 
	Runge-Kutta methods for gradient flow problems}, 
	{\em Math. Comput.}, 2024, doi: 10.1090/mcom/4015.

\bibitem{LiaoZhang:2021}
{\sc H.-L. Liao and Z. Zhang},
{Analysis of adaptive BDF2 scheme for diffusion equations},
{\em Math. Comput.}, 90 (2021), pp. 1207-1226.

\bibitem{LiuZou:2006}
{\sc H. Liu and J. Zou},
{Some new additive Runge-Kutta methods and their applications},
{\em J. Comput. Appl. Math.}, 190 (2006), pp. 74-98.



\bibitem{MolerVanLoan:2003}
{\sc C. Moler and C. Van Loan},
{Nineteen dubious ways to compute the exponential of a matrix, twenty-five years later}, 
{\em SIAM Rev.}, 45:1 (2003), pp. 3-49.

%



\bibitem{OlssonSoderlind:2000}
{\sc H. Olsson and G. S\"{o}derlind},
The approximate Runge-Kutta computational process,
{\em BIT}, 40:2 (2000), pp. 351-373.




\bibitem{ShenXuYang:2018}
{\sc J. Shen, J. Xu and J. Yang},
{The scalar auxiliary variable (SAV) approach for gradient flows},
{\em J. Comput. Phys.}, 353 (2018), pp. 407-416.

\bibitem{ShinLeeLee:2017CMA}
{\sc J. Shin, H.-G. Lee and J.-Y. Lee},
{Convex splitting Runge-Kutta methods for phase-field models},
{\em Comput. Math. Appl.}, 73:11 (2017), pp. 2388-2403.


\bibitem{ShinLeeLee:2017}
{\sc J. Shin, H.-G. Lee and J.-Y. Lee},
{Unconditionally stable methods for gradient flow using convex splitting Runge-Kutta scheme},
{\em J. Comput. Phys.}, 347 (2017), pp. 367-381.


\bibitem{StuartHumphries:1998}
{\sc A. M. Stuart and A. R. Humphries},
{Dynamical systems and numerical analysis},
Cambridge University Press, New York, 1998.


\bibitem{Zhong:1995}
{\sc X. Zhong},
{Additive Semi-implicit Runge-Kutta methods for computing high-speed nonequilibrium reactive flows},
{\em J. Comput. Phys.}, 128 (1996), pp. 19-31. 

%
%
%


%
%
%

\end{thebibliography}
\end{document}